\def\ArxivBuild{}
\documentclass[aos]{imsart}

\newif\ifArxiv
\ifdefined\ArxivBuild
    \Arxivtrue
\fi

\ifArxiv\else

\makeatletter\@load@xr@or@xr@hyperfalse\makeatother
\fi

\RequirePackage{amsthm,amsmath,amsfonts,amssymb}
\RequirePackage[authoryear]{natbib}
\ifArxiv\else\RequirePackage{xr-hyper}\fi
\RequirePackage[colorlinks,citecolor=blue,urlcolor=blue]{hyperref}

\newcommand{\revision}[1]{\ifArxiv {#1}\else\textcolor{red}{#1}\fi}

\usepackage{dsfont}
\usepackage{mathtools}
\mathtoolsset{showonlyrefs=true,showmanualtags=true}
\usepackage[only,llfloor,rrfloor]{stmaryrd}
\usepackage{xspace}
\usepackage{graphicx}

\makeatletter
\def\journal@name{Working paper}
\def\journal@url{}
\makeatother
\hypersetup{pdfsubject={Working paper}}

\numberwithin{equation}{section}
\theoremstyle{plain}
\newtheorem{theorem}{Theorem}[section]
\newtheorem{proposition}[theorem]{Proposition}
\newtheorem{lemma}[theorem]{Lemma}
\newtheorem{corollary}[theorem]{Corollary}

\newtheorem*{theorem*}{Theorem}
\newtheorem*{proposition*}{Proposition}
\theoremstyle{definition}
\newtheorem{definition}[theorem]{Definition}

\newtheorem{remark}[theorem]{Remark}

\newcommand{\N}{\mathbb{N}}
\newcommand{\R}{\mathbb{R}}
\newcommand{\PP}{\mathsf{P}}
\newcommand{\EE}{\mathsf{E}}
\newcommand{\Var}{\mathsf{Var}}

\newcommand{\1}{\mathds{1}}
\newcommand{\rd}{\mathrm{d}}

\newcommand{\ConditionSymbol}[2][]{%
  \ensuremath{\boldsymbol{\mathrm{#2}}%
  \if\relax\detokenize{#1}\relax\else_{#1}\fi}\xspace}
\newcommand{\ATC}[1][]{\ConditionSymbol[#1]{ATC}}
\newcommand{\PBM}[1][]{\ConditionSymbol[#1]{PBM}}
\newcommand{\BM}[1][]{\ConditionSymbol[#1]{BM}}
\newcommand{\PL}[1][]{\ConditionSymbol[#1]{PL}}

\newcommand{\domain}{\mathcal{D}}
\newcommand{\Dk}{\mathcal{D}_k}
\newcommand{\Dexp}{\mathcal{D}_{\mathrm{exp}}}
\newcommand{\Dlad}{\mathcal{D}_{\mathrm{lad}}}

\newcommand{\bh}{\mathbf{h}}
\newcommand{\hmax}{{h^*}}
\newcommand{\Escale}[1]{\mathcal{V}_{#1}}
\newcommand{\Eh}{\Escale{h}}
\newcommand{\Ehmax}{\Escale{\hmax}}

\newcommand{\Pm}{\mathcal{P}_{\!m}}
\newcommand{\Vscale}[1]{V_{#1}}
\newcommand{\Vh}{\Vscale{h}}
\newcommand{\Wscale}[1]{W_{#1}}
\newcommand{\Wh}{\Wscale{h}}
\newcommand{\Bg}{\mathcal{B}_\gamma}

\newcommand{\Bgiso}{\mathcal{B}_\gamma^{\mathrm{iso}}}
\newcommand{\lamiso}{\lambda_\gamma^{\mathrm{iso}}}

\newcommand{\kg}{\revision{k_\gamma(0,0)}}

\newcommand{\Qg}{Q_\gamma}
\newcommand{\Zm}{\mathcal{Z}_m}
\newcommand{\deltam}{\delta_m(h)}
\newcommand{\tmu}{\revision{\mu_h}}
\newcommand{\Gn}{\Gamma_n}

\DeclareMathOperator*{\argmin}{arg\,min}

\newcommand{\pwrisk}{R_n}
\newcommand{\finfty}{\|f(t+\cdot)\|_{\infty,\Ehmax}}
\newcommand{\FsL}{\mathcal{F}_{\!s,L}}

\ifArxiv\else\fi

\newcommand{\stmtlabel}[1]{\label{#1}}
\newenvironment{stmtequation}{\equation}{\endequation}
\newenvironment{stmtalign}{\align}{\endalign}

\newcommand{\StmtBiasFactors}{%
    For $\gamma=(m,{\bh})$ and $r>0$, define
    \begin{stmtalign}
        \Lambda_m({\bh})
        &:=\int_{\Escale{\bh}}|\revision{k}_\gamma(0,u)|\revision{\mu}_{\bh}(u)\,\rd u,
        \stmtlabel{eq:Lambda.def}\\
        \mathcal A_{m,r}({\bh})
        &:=\int_{\Escale{\bh}}|\revision{k}_\gamma(0,u)|\|u\|_\infty^r \revision{\mu}_{\bh}(u)\,\rd u.
        \stmtlabel{eq:bias.moment.def}
    \end{stmtalign}%
}

\newcommand{\StmtRiskBound}{%
\revision{Let $\gamma=(m,\bh)$ with $m\in\N_0$, ${\bh}\in(0,\hmax]^d$.
    \begin{enumerate}
    \item[(a)] We have for any $r>0$,
    \begin{stmtequation}\stmtlabel{eq:Lambda.bounds}
        1\leq\Lambda_m({\bh})\leq\delta_m({\bh})^{-1/2},
        \qquad
        \mathcal A_{m,r}({\bh})\leq \|{\bh}\|_\infty^r\Lambda_m({\bh}).
    \end{stmtequation}
    \item[(b)] The variance of $\hat f_\gamma$ satisfies
    \begin{stmtequation}
     \Var\bigl(\hat f_\gamma(t)\bigr)\le \frac{\finfty \kg}{n\revision{V_{\bh}}}=\frac{\finfty}{n\Vscale{\bh}\,\delta_m(\bh)}.
    \end{stmtequation}
   \item[(c)] Let $s,L>0$ and $f\in\Sigma(s,L)$.
    Set $\beta_m(s):=\min(m+1,s)\in(0,s]$.
    There exists $\mathfrak{L}_{m,s,L}>0$, independent of ${\bh}$ and of $f\in\Sigma(s,L)$, with $\mathfrak{L}_{m,s,L}=L$ whenever $m\ge\llfloor s\rrfloor$, such that, for ${\bh}\in(0,\hmax]^d$ and $\gamma=(m,{\bh})$, the pointwise bias $b_\gamma:=\EE\hat f_\gamma(t)-f(t)$ satisfies
    \begin{stmtequation}\stmtlabel{eq:bias.bound.general}
        |b_\gamma|
        \leq\mathfrak L_{m,s,L}\mathcal A_{m,\beta_m(s)}({\bh})
        \leq\mathfrak L_{m,s,L}\|{\bh}\|_\infty^{\beta_m(s)}\Lambda_m({\bh}).
    \end{stmtequation}
\item[(d)] Finally we have
    \begin{stmtalign}
        \pwrisk\bigl(\hat f_\gamma,f\bigr)^2
        &\;\le\;
        \bigl(\mathfrak L_{m,s,L}\Lambda_m({\bh})\|{\bh}\|_\infty^{\beta_m(s)}\bigr)^2
        +\frac{\finfty}{n\Vscale{\bh}\,\delta_m(\bh)}
        \stmtlabel{eq:risk.bound}
    \end{stmtalign}%
\end{enumerate}  } 
}

\newcommand{\StmtSupBound}{%
    Let $s,L>0$. There exists a finite constant $\FsL$, depending only
    on $s,L,d,\domain,t$ and $\hmax$, such that every
    $f\in\Sigma(s,L)$ satisfies $\finfty\leq\FsL$.%
}

\newcommand{\StmtPBM}{%
    Fix $m\in\N_0$. The following conditions are equivalent:
    \begin{enumerate}
        \item[\textup{(i)}] \PBM[m] holds;
        \item[\textup{(ii)}]
        $\displaystyle\inf_{0<h\leq\hmax}\delta_m(h)>0$;
        \item[\textup{(iii)}]
        $\displaystyle\sup_{0<h\leq\hmax}
        \revision{k_{m,h}(0,0)}<\infty$.
    \end{enumerate}%
}

\newcommand{\StmtBM}{%
    Fix $m\in\N_0$. The following conditions are equivalent:
    \begin{enumerate}
        \item[\textup{(i)}] \BM[m] holds;
        \item[\textup{(ii)}]
        $\displaystyle\sup_{0<h\leq\hmax}Q_{m,h}<\infty$.
    \end{enumerate}
    In particular, \BM[m] implies \PBM[m], and one may take
    $C_m^{\mathrm{PBM}}\leq C_m^{\mathrm{BM}}$ since $0\in\Eh$.%
}

\newcommand{\StmtRescaledSetup}{%
    Let ${\bh}\in(0,\hmax]^d$ and let \revision{ $A_{\bh}$ be an invertible $d\times d$ matrix.}
     In this section, limits ($\longrightarrow$) are taken along a specified family with $\|{\bh}\|_\infty\to0$ that will be chosen differently in Sections~\ref{sec:results} to \ref{sec:Dlad}  according to $\mathcal{D}$.
    Define \revision{$\sigma_{\bh}(u)=A_{\bh}^{-1}u$ for $u\in\R^d$ and let $\varrho_{\bh}=(\sigma_{\bh})_\#\mu_{\bh}$ be the pushforward of $\mu_{\bh}$ by $\sigma_{\bh}$} and 
    \[       \qquad G_{m,{\bh}}=\int\Phi(v)\Phi(v)^\top\,\rd\varrho_{\bh}(v).
    \]
    The choice $A_{\bh}=\operatorname{diag}(h_1,\ldots,h_d)$ normalises
    $\Escale{\bh}$ coordinate-wise.
    Other invertible linear rescalings are allowed.
}

\newcommand{\StmtMomentStability}{%
    Fix $m\in\N_0$ and suppose that we have $G_{m,{\bh}}\longrightarrow G_{m,0}$ where $G_{m,0}$ is positive definite.
    Define
    \[\revision{\tilde{k}}_{m,0}(v,w)=\Phi(v)^\top G_{m,0}^{-1}\Phi(w).\]
    Then, uniformly on each compact subset of $\R^d\times\R^d$,
    \begin{stmtequation}\stmtlabel{eq:moment.kernel.limit}
        k_{m,{\bh}}(A_{{\bh}}v,A_{{\bh}}w)\longrightarrow\revision{\tilde k}_{m,0}(v,w).
    \end{stmtequation}
    In particular,
    \begin{stmtequation}\stmtlabel{eq:moment.delta.limit}
       \delta_m({\bh})\longrightarrow
        \delta_{m,0}:=(e_1^\top G_{m,0}^{-1}e_1)^{-1}>0,
        \qquad \limsup_{{\bh}\to0}\Lambda_m({\bh})<\infty.
    \end{stmtequation}
    The following additional conclusions hold.
    \begin{enumerate}
        \item[\textup{(i)}] If $\varrho_{\bh}\Rightarrow\varrho_0$ and
        $G_{m,0}=\int\Phi\Phi^\top\,\rd\varrho_0$, then
        \begin{stmtequation}\stmtlabel{eq:moment.Lambda.limit}
           \Lambda_m({\bh})\longrightarrow
            \int|\tilde{k}_{m,0}(0,v)|\,\rd\varrho_0(v).
        \end{stmtequation}
        \item[\textup{(ii)}] If $\sigma_{\bh}(\Escale{\bh})$ is contained in a fixed
        compact set for all sufficiently small bandwidths, then
        \begin{stmtequation}\stmtlabel{eq:moment.uniform}
            \limsup_{{\bh}\to0}\sup_{u\in\Escale{\bh}}k_{m,{\bh}}(u,u)<\infty.
        \end{stmtequation}
    \end{enumerate}%
}

\newcommand{\StmtGeometricStability}{%
    Suppose that $F_{\bh}:=\sigma_{\bh}(\Escale{\bh})$ is contained in a fixed compact
    set $B$, and that there is a measurable set $F_0\subset B$ with
    $\mathrm{Leb}(F_0)>0$ such that
    \[
        \varepsilon_{\bh}:=\mathrm{Leb}(F_{\bh}\triangle F_0)\longrightarrow0.
    \]
    Then all conclusions of Proposition~\ref{prop:moment.stability}
    hold at every fixed degree, with $\varrho_0$ the uniform probability
    measure on $F_0$ and $G_{m,0}=\int\Phi\Phi^\top\,\rd\varrho_0$.
    Moreover,
    \begin{stmtequation}\stmtlabel{eq:geometric.volume}
        \Vscale{\bh}=|\det A_{\bh}|\bigl(\mathrm{Leb}(F_0)+O(\varepsilon_{\bh})\bigr).
    \end{stmtequation}
    The errors in~\eqref{eq:moment.kernel.limit},
    \eqref{eq:moment.delta.limit} for $\delta_m$, and
    \eqref{eq:moment.Lambda.limit} are $O(\varepsilon_{\bh})$.
    No inclusion of the origin in $F_0$ is required.%
}

\newcommand{\StmtMinimaxUpper}{%
    Let $s,L>0$ and $\nu\geq0$, and set $m=\llfloor s\rrfloor$.
    Suppose that \PBM[m] and $\PL[-](\nu)$ hold at $t$.
    Set $h_n=\hmax\wedge n^{-1/(2s+d+\nu)}$.
    Then, for every $n\geq3$,
    \[
        \sup_{f\in\Sigma(s,L)}
        \pwrisk\bigl(\hat f_{(m,h_n)},f\bigr)
        \;\lesssim\;
        n^{-s/(2s+d+\nu)}.
    \]%
}

\newcommand{\StmtMinimaxLower}{%
    Let $s,L>0$ and $\nu\geq0$.
    If $\PL[+](\nu)$ holds at $t$, then, for every $n\geq3$,
    \[
        \inf_{\hat f}\sup_{f\in\Sigma(s,L)}
        \pwrisk(\hat f,f)
        \;\gtrsim\;
        n^{-s/(2s+d+\nu)}.
    \]%
}

\newcommand{\StmtMinimaxRate}{%
    Let $s,L>0$ and $\nu\geq0$.
    Suppose that $\domain$ is regular of degree $\llfloor s\rrfloor$ at $t$, with exponent $\nu$.
    Then, for every $n\geq3$,
    \[
        \inf_{\hat{f}}\,\sup_{f\in\Sigma(s,L)}\pwrisk(\hat{f},f)
        \;\asymp\;
        n^{-s/(2s+d+\nu)}.
    \]
    The estimator $\hat f_{(\llfloor s\rrfloor,h_n)}$, with $h_n=\hmax\wedge n^{-1/(2s+d+\nu)}$, attains this minimax rate.%
}

\newcommand{\StmtOracle}{%
    Suppose $\PL[-](\nu)$ holds at $t$, let $f$ be a density on $\domain$ with $\finfty<\infty$, and fix $\vartheta>0$ and $A>d+\nu$.
    There exist constants $C_1,C_2>0$, depending only on $\vartheta,A,d,\nu,\hmax$, the constant in $\PL[-](\nu)$, and on $f$ only through $\finfty$, such that, for every $n\geq3$ for which the grid $\Gn$ is nonempty,
    \begin{stmtequation}\stmtlabel{eq:oracle.adaptive}
        \pwrisk(\hat f,f)
        \;\leq\;
        C_1\min_{\gamma\in\Gn}
        \Bigl\{B(\gamma) + r_\gamma\bigl(v_\gamma^\star+2\varepsilon_\gamma,\,\mathrm{pen}(\gamma)\bigr)\Bigr\}
        \;+\;\frac{C_2}{\sqrt n}.
    \end{stmtequation}%
}

\newcommand{\StmtAdaptive}{%
    Suppose that $\domain$ is uniformly regular at $t$, with exponent $\nu$.
    Fix $\vartheta>0$ and $A>d+\nu$.
    Then, for every $s,L>0$, there exists a constant $C_{\mathrm{ad}}>0$ such that, for every $n\geq3$,
    \begin{stmtequation}\stmtlabel{eq:rate.adaptive}
        \sup_{f\in\Sigma(s,L)}\pwrisk(\hat f,f)
        \;\leq\;
        C_{\mathrm{ad}}
        \left(\frac{\log n}{n}\right)^{s/(2s+d+\nu)}.
    \end{stmtequation}
    The constant $C_{\mathrm{ad}}$ depends only on
    $s,L,d,\nu,\hmax,\vartheta,A$, the constants in $\PL(\nu)$, and
    the constants in $\BM[m]$ for finitely many degrees $m$ determined
    by $s,d,\nu$.%
}

\newcommand{\StmtATCPBM}{%
    Suppose that condition~\textup{(a)} of $\ATC(\alpha)$ holds at $t$, and set $\nu=|\alpha|-d$.
    Then $\domain$ is regular at $t$, with exponent $\nu$.
    More precisely, for every $m\in\N_0$, there exists $\delta_{m,0}\in(0,1]$ such that
    \[
        \delta_m(h)\longrightarrow\delta_{m,0}
        \qquad\text{as }h\to0.
    \]%
}

\newcommand{\StmtATCBM}{%
    Suppose that $\domain$ satisfies $\ATC(\alpha)$ at $t$, and set $\nu=|\alpha|-d$.
    Then $\domain$ is uniformly regular  at $t$, with exponent $\nu$.

    More precisely, for every fixed $m\in\N_0$, there exist $C_m^{\mathrm{ATC}}<\infty$ such that, for $0<h\leq\hmax$,
  \revision{  \[
       1
        \;\leq\; k_{m,h}(0,0)
        \;\leq\; Q_{m,h}
        \;\leq\; C_m^{\mathrm{ATC}}.
    \]%
}
}

\newcommand{\StmtDladVolume}{%
    There exist $c,C>0$ such that, for all sufficiently small $h$,
    \revision{\begin{stmtequation}\stmtlabel{eq:volume.Dlad}
        c\,h^{2+\nu}\,\theta_{k(h)}^{1+\nu}\;\le\;V_h\;\le\;C\,h^{2+\nu}.
    \end{stmtequation}}
    
    In particular, after adjusting the constants if necessary,
    \begin{stmtequation}\stmtlabel{eq:volume.Dlad22}
        c\,\exp\left\{-2(1+\nu)\sqrt{\log (1/h)}\right\}\;\le\;\frac{V_h}{h^{2+\nu}}\;\le\;C.
    \end{stmtequation}
    Moreover,
    \[
        V_h=h^{2+\nu+o(1)},
        \qquad
        \liminf_{h\downarrow0}h^{-(2+\nu)}V_h=0,
        \qquad
        0<\limsup_{h\downarrow0}h^{-(2+\nu)}V_h<\infty.
    \]
    Consequently, \PL does not hold at $t$.

}

\newcommand{\StmtDladPBM}{%
    Let $h_k:=2\zeta_k$.
    There exists $C>0$ such that, for every $m\ge1$ and all sufficiently
    large $k$,
    \begin{stmtequation}\stmtlabel{eq:delta.Dlad}
        \delta_m(h_k)\;\le\;\delta_1(h_k)\;\le\;
        C\,\theta_{k+1}^{1+\nu}\;\xrightarrow[k\to\infty]{}\;0.
    \end{stmtequation}
    Consequently, \PBM[m] fails for every $m\ge1$, while
    $\Lambda_1(h_k)$ remains bounded along the same sequence.%
}

\newcommand{\StmtDladLow}{%
    Let $0<s\le1$ and $L>0$.
    For every $n\ge3$, define $h_n$ as follows.
    If $n\hmax^{2s}\Vscale{\hmax}\ge1$, let $h_n$ be the unique solution in
    $(0,\hmax]$ of
    \begin{stmtequation}\stmtlabel{eq:hn.implicit.Dlad}
        \revision{n h_n^{2s}\Vscale{h_n}=1}.
    \end{stmtequation}
    Otherwise, set $h_n=\hmax$.
    Then
    \begin{stmtequation}\stmtlabel{eq:rate.implicit.Dlad}
        \mathcal R_n(s,L)\asymp h_n^s,
    \end{stmtequation}
    where the upper bound is attained by the local-constant estimator $\hat f_{(0,h_n)}$.
    Moreover, this minimax rate oscillates.
    More precisely,
    \begin{stmtalign}
        \liminf_{n\to\infty}
        \frac{\log\mathcal R_n(s,L)
        +\frac{s}{2s+2+\nu}\log n}{\sqrt{\log n}}
        &=0,
        \stmtlabel{eq:oscillations.Dlad.inf}\\
        \limsup_{n\to\infty}
        \frac{\log\mathcal R_n(s,L)
        +\frac{s}{2s+2+\nu}\log n}{\sqrt{\log n}}
        &=\frac{2s(1+\nu)}{(2s+2+\nu)^{3/2}} \eqcolon c_{s,\nu}^*.
        \stmtlabel{eq:oscillations.Dlad.sup}
    \end{stmtalign}
    Along the subsequence
    \[
        n_k^-:=\left\lceil
        \bigl\{(2\zeta_k)^{2s}\Vscale{2\zeta_k}\bigr\}^{-1}
        \right\rceil,
    \]
    the baseline order is attained:
    \[
        \mathcal R_{n_k^-}(s,L)
        \asymp (n_k^-)^{-s/(2s+2+\nu)}.
    \]%
}

\newcommand{\StmtDladAnisotropic}{%
    Let $1<s\leq2$ and $L>0$, and use the definitions
    in~\eqref{eq:Dlad.high.windows}--\eqref{eq:Dlad.high.estimator}.
    For all sufficiently large $n$, the bandwidth $\bh_{k_n,n}$ belongs to
    $(0,\hmax]^2$, and
    \begin{stmtequation}\stmtlabel{eq:Dlad.high.minimax}
        \mathcal R_n(s,L)
        \asymp\sup_{f\in\Sigma(s,L)}\pwrisk(\hat f_n^{\rm lad},f)
        \asymp\mathfrak r_n(s).
    \end{stmtequation}
    Writing $a=2s+2+\nu$, we have
    \begin{stmtalign}
        \liminf_{n\to\infty}
        \frac{\log\mathcal R_n(s,L)+(s/a)\log n}{\sqrt{\log n}}
        &=\frac{2(1+\nu)(s-1)}{a^{3/2}},
        \stmtlabel{eq:Dlad.high.liminf}\\
        \limsup_{n\to\infty}
        \frac{\log\mathcal R_n(s,L)+(s/a)\log n}{\sqrt{\log n}}
        &=\frac{2(1+\nu)(2s-1)}{a^{3/2}}.
        \stmtlabel{eq:Dlad.high.limsup}
    \end{stmtalign}
    These limits are attained, respectively, along
    \begin{stmtequation}\stmtlabel{eq:Dlad.high.subsequences}
        n_k^-=
        \left\lceil B_k(s)^{-(2s+1)/s}\omega_k^{-1}\right\rceil,
        \qquad
        n_k^+=
        \left\lceil B_k(s)^{-(2s+1)/s}\omega_{k+1}^{-1}\right\rceil.
    \end{stmtequation}
    The same limits hold with $\mathcal R_n(s,L)$ replaced by
    $\sup_{f\in\Sigma(s,L)}\pwrisk(\hat f_n^{\rm lad},f)$.
    The comparison constants may depend on $s,L,\nu,\hmax$ and $k_0$.%
}

\newcommand{\StmtDexpBest}{%
    Let $s,L>0$, and let $M\ge\llfloor s\rrfloor$ be a fixed integer.
    Set
    \begin{equation*}
        \alpha_\star(s):=
        \max_{0\le m\le\llfloor s\rrfloor}\min\{m+1,s-m\}.
    \end{equation*}
    Then
    \begin{stmtequation}\stmtlabel{eq:rate.Dexp.best.fixed.degree}
        \inf_{\hat f_\gamma\in\mathcal S_M}
        \sup_{f\in\Sigma(s,L)}\pwrisk(\hat f_\gamma,f)
        =\min_{0\le m\le M}\;\inf_{0<h\le\hmax}\;
        \sup_{f\in\Sigma(s,L)}
        \pwrisk\bigl(\hat f_{(m,h)},f\bigr)
        \asymp(\log n)^{-\alpha_\star(s)}.
    \end{stmtequation}
    The maximising degree is one of the integers closest to
    $(s-1)/2$. In particular, $\alpha_\star(s)<s$ for every $s>1$.%
}

\newcommand{\StmtDexpMinimax}{%
    For every $s,L>0$,
    \begin{stmtequation}\stmtlabel{eq:rate.Dexp}
        \inf_{\hat f}\sup_{f\in\Sigma(s,L)}\pwrisk(\hat f,f)
        \;\asymp\;(\log n)^{-s}.
    \end{stmtequation}
    For all sufficiently large $n$, the estimator
    $\hat f_{(\llfloor s\rrfloor,\bh(r_n))}$ attains this rate.
    For $0<s\leq1$ and all sufficiently large $n$, the scalar-bandwidth
    estimator $\hat f_{(0,h_n)}$ also attains it.%
}

\newcommand{\StmtDexpAnisotropic}{%
    Fix $m\in\N_0$, and let $0<s\leq m+1$ and $L>0$.
    For $\bh(r)=(2r,e^{-1/r})$, there exist $r_0>0$ and $C<\infty$
    such that, for $0<r\leq r_0$ and $n\geq3$,
    \begin{stmtequation}\stmtlabel{eq:Dexp.anisotropic.risk}
        \sup_{f\in\Sigma(s,L)}
        \pwrisk\bigl(\hat f_{(m,\bh(r))},f\bigr)^2
        \leq C\left(r^{2s}+\frac{e^{1/r}}{nr}\right).
    \end{stmtequation}
    Set
    \begin{stmtequation}\stmtlabel{eq:Dexp.anisotropic.bandwidth}
        r_n=\frac{2}{\log n},
        \qquad \bh(r_n)=\left(\frac{4}{\log n},n^{-1/2}\right).
    \end{stmtequation}
    Then, for all sufficiently large $n$,
    \begin{stmtequation}\stmtlabel{eq:Dexp.anisotropic.rate}
        \sup_{f\in\Sigma(s,L)}
        \pwrisk\bigl(\hat f_{(m,\bh(r_n))},f\bigr)
        \leq C'(\log n)^{-s}.
    \end{stmtequation}
    The constants may depend on $m,s,L$ and $\hmax$.%
}

\newcommand{\StmtDefSigma}{%
    Let $(s,L)\in(0,\infty)^2$.
    A probability density $f$ on $\domain$ belongs to $\Sigma(s,L)$ if there exists $q\in\mathcal{P}_{\llfloor s\rrfloor}$ such that
    \begin{stmtequation}\stmtlabel{eq:holder}
        |f(t+u)-q(u)|\leq L\|u\|_\infty^s,
        \qquad\forall\,u\in \Ehmax.
    \end{stmtequation}
    In particular, taking $u=0$ gives $q(0)=f(t)$.
    The target $f(t)$ is thus the value at the origin of the local approximating polynomial.
    We write
    \[
        \Sigma:=\bigcup_{s>0}\bigcup_{L>0}\Sigma(s,L).
    \]%
}

\newcommand{\StmtDefDelta}{%
    For $\gamma=(m,{\bh})$ with $m\in\N_0$ and ${\bh}\in(0,\hmax]^d$,
    \[
        \delta_m(\bh)
        := \frac{1}{\revision{k_\gamma(0,0)}}
        = \inf_{\substack{p\in\Pm\\ p(0)=1}}\|p\|_{L^2(\mu_{\bh})}^2
        = \inf_{p\in\Zm}\|1-p\|_{L^2(\mu_{\bh})}^2 ,
    \]%
    where $\Zm=\{p\in\Pm:p(0)=0\}=\operatorname{Span}\{\varphi_\alpha:1\leq|\alpha|\leq m\}$.
}

\newcommand{\StmtDefPL}{%
    Fix $\nu\geq0$. The domain satisfies $\PL[-](\nu)$ at $t$ if there
    exists $c_\nu^{\mathrm{PL}}>0$ such that
    $\revision{V_h}\geq c_\nu^{\mathrm{PL}}h^{\revision{d+\nu}}$ for every
    $h\in(0,\hmax]$. It satisfies $\PL[+](\nu)$ if there exists
    $C_\nu^{\mathrm{PL}}<\infty$ such that
    $\revision{V_h}\leq C_\nu^{\mathrm{PL}}h^{\revision{d+\nu}}$ on the same interval.
    We write $\PL(\nu)$ when both conditions hold, equivalently when
    \begin{stmtequation}\stmtlabel{eq:PL.def}
        c_\nu^{\mathrm{PL}}\,h^{\revision{d+\nu}}
        \;\leq\; \revision{V_h}\;\leq\;
        C_\nu^{\mathrm{PL}}\,h^{\revision{d+\nu}},
        \qquad h\in(0,\hmax].
    \end{stmtequation}
    Finally, \PL means that $\PL(\nu)$ holds for some $\nu\geq0$.%
}

\newcommand{\StmtDefPBM}{%
    We say that \PBM[m] holds at $t$ if there exists
    $C_m^{\mathrm{PBM}}<\infty$ such that
    \[
        |p(0)|^2
        \;\leq\; C_m^{\mathrm{PBM}}\,\|p\|_{L^2(\revision{\mu_h})}^2,
        \qquad
        p\in\Pm,\quad h\in(0,\hmax].
    \]%
}

\newcommand{\StmtDefBM}{%
    We say that \BM[m] holds at $t$ if there exists $C_m^{\mathrm{BM}}<\infty$ such that
    \[
        \|p\|_{\infty,\Eh}^2
        \;\leq\; C_m^{\mathrm{BM}}\,\|p\|_{L^2(\revision{\mu_h})}^2,
        \qquad
        p\in\Pm,\quad h\in(0,\hmax].
    \]%
}

\newcommand{\StmtDefRegular}{%
    Fix $\nu\geq0$.
    For $m\in\N_0$, the domain $\domain$ is \emph{regular of degree $m$ at $t$}, with exponent $\nu$, if it satisfies $\PL(\nu)$ and \PBM[m].
    It is \emph{uniformly regular of degree $m$ at $t$}, with exponent $\nu$, if it satisfies $\PL(\nu)$ and \BM[m].
    The domain is \emph{regular at $t$}, with exponent $\nu$, if it is regular of every fixed degree.
    Equivalently, it satisfies $\PL(\nu)$ and \PBM.
    It is \emph{uniformly regular at $t$}, with exponent $\nu$, if it is uniformly regular of every fixed degree.
    Equivalently, it satisfies $\PL(\nu)$ and \BM.
    Here uniformity refers to polynomial evaluation over the local neighbourhood, not to the degree $m$.
    The constants may still depend on $m$.%
}

\newcommand{\StmtDefATC}{%
    Let $\alpha=(\alpha_1,\ldots,\alpha_d)\in(0,\infty)^d$, with associated
    coordinate-wise rescaling maps
    $\sigma_h(u)=(u_1/h^{\alpha_1},\ldots,u_d/h^{\alpha_d})$, $h>0$.
    We say that $\domain$ satisfies $\ATC(\alpha)$ if the following two
    conditions hold:
    \begin{enumerate}
        \item[\textup{(a)}] for every $\beta\in\N_0^d$, there exists
        a constant $c_\beta\in\R$, with $c_0>0$, such that, as $h\to 0$,
        \revision{\begin{stmtequation}\stmtlabel{eq:ATC}
            \int_{\Eh}u^\beta\,\rd u
            = h^{\langle\alpha,\beta\rangle+d+\nu}\bigl(c_\beta+o(1)\bigr),
        \end{stmtequation}}
        where $\nu=|\alpha| - d$ is the volume exponent (necessarily
        $\nu\geq0$: the case $\beta=0$ gives \revision{$V_h\sim c_0h^{d+\nu}$} and $V_h\leq 2^d h^d$);
        \item[\textup{(b)}] there exists $R>0$ such that
        $\sigma_h(\Eh)\subseteq[-R,R]^d$ for all $h\in(0,\hmax]$.
    \end{enumerate}%
}

\AtBeginDocument{\noeqref{%
eq:ATC,%
eq:Agamma.new,%
eq:Dexp.anisotropic.bandwidth,%
eq:Dexp.anisotropic.geometry,%
eq:Dexp.anisotropic.rate,%
eq:Dexp.anisotropic.risk,%
eq:Dlad.high.estimator,%
eq:Dlad.high.liminf,%
eq:Dlad.high.limsup,%
eq:Dlad.high.minimax,%
eq:Dlad.high.profile,%
eq:Dlad.high.subsequences,%
eq:Lambda.bounds,%
eq:Lambda.def,%
eq:Ln.def,%
eq:PL.def,%
eq:Qg,%
eq:Qg.def,%
eq:bias.bound.general,%
eq:bias.moment.def,%
eq:cgamma.def,%
eq:delta.Dlad,%
eq:eps.gamma,%
eq:gammahat.new,%
eq:geometric.volume,%
eq:grid.new,%
eq:hatv.new,%
eq:hl.ml,%
eq:hn.Dexp,%
eq:hn.implicit.Dlad,%
eq:holder,%
eq:moment.Lambda.limit,%
eq:moment.delta.limit,%
eq:moment.kernel.limit,%
eq:moment.uniform,%
eq:oracle.adaptive,%
eq:oscillations.Dlad.inf,%
eq:oscillations.Dlad.sup,%
eq:rate.Dexp,%
eq:rate.Dexp.best.fixed.degree,%
eq:rate.adaptive,%
eq:rate.implicit.Dlad,%
eq:risk.bound,%
eq:volume.Dexp,%
eq:volume.Dlad,%
eq:volume.Dlad22}}

\begin{document}


\begin{frontmatter}
\title{A minimax theory for pointwise density estimation on geometrically complex domains}
\runtitle{Pointwise density estimation on complex domains}

\begin{aug}
\author[a1]{
    \fnms{Karine}~\snm{Bertin}%
    \ead[label=e1]{karine.bertin@uv.cl}%
    \orcid{0000-0002-8167-137X}%
}

\author[a2]{
    \fnms{Nicolas}~\snm{Klutchnikoff}%
    \ead[label=e2]{nicolas.klutchnikoff@univ-rennes2.fr}%
    \orcid{0000-0001-9349-0771}%
}

\address[a1]{
    CIMFAV-INGEMAT, Universidad de Valpara\'{\i}so,\space
    Valpara\'{\i}so, Chile%
    \printead[presep={,\ }]{e1}%
}
\address[a2]{
    Univ Rennes, CNRS, IRMAR-UMR 6625,\space
    F-35000 Rennes, France%
    \printead[presep={,\ }]{e2}%
}

\end{aug}

\begin{abstract}
We study pointwise density estimation on a known domain whose boundary may have geometric singularities.
Our main contribution is a variance bound for local polynomial estimators expressed through a local Christoffel function.
This quantity separates the effect of local volume from the stability of polynomial evaluation at the estimation point.
It can yield substantially sharper bounds than spectral estimates based on the smallest eigenvalue of the local Gram matrix.
We combine this variance control with bounds on weighted absolute moments of the equivalent kernel to control the bias.
Under polynomial local volume growth and a pointwise Bernstein--Markov condition, we establish minimax rates for every positive smoothness level over isotropic H\"older-type classes.
A uniform Bernstein--Markov condition yields adaptation to unknown smoothness up to a logarithmic factor.
A purely geometric anisotropic tangent cone condition, readily verifiable in concrete examples, implies these assumptions and covers polynomial sectors.
Beyond this setting, an exponential cusp has logarithmic minimax rates, while a connected ladder exhibits oscillations for smoothness up to order two, with exact lower and upper envelopes.
On these domains, anisotropic bandwidths attain the minimax rates over isotropic classes.
On the cusp, deterministic scalar bandwidths with bounded polynomial degrees are suboptimal for smoothness greater than one.
These examples show how support geometry affects both the attainable rates and the choice of bandwidth.
\end{abstract}



\end{frontmatter}


\section{Introduction}\label{sec:introduction}

\paragraph{Density estimation on a known domain}
We study estimation of a density $f$ at a fixed point $t$ of a known domain $\domain\subset\R^d$, from $n$ independent observations.
At an interior point, the minimax pointwise $L^2$ risk over an isotropic smoothness class of order $s$ is $n^{-s/(2s+d)}$.
Adaptation to unknown smoothness is possible at a logarithmic cost \citep{Lepski1991adaptive,GoldenshlugerLepski2011}.
At a boundary point, estimation neighbourhoods are truncated by the domain.
This truncation changes the bias of standard kernel estimators and reduces the amount of local information.
At a geometric singularity, it can also make polynomial evaluation at the target unstable.
We study how these effects determine the minimax risk and the performance of local polynomial estimators.

Boundary correction has a substantial literature, including reflection methods \citep{Schuster1985}, boundary kernels \citep{GasserMuller1979,Muller1991}, Bernstein and asymmetric kernels \citep{Vitale1975,Chen1999}, transformations \citep{MarronRuppert1994} and local polynomial fitting \citep{FanGijbels1996}.
On intervals, boxes and smooth domains, such corrections can preserve the interior rate.
More complicated observation domains arise, for instance, in spatial density estimation \citep{DaviesMarshallHazelton2018}.
Methods designed for these settings include diffusion estimators \citep{BotevGrotowskiKroese2010,BaddeleyEtAl2022} and PDE-regularised maximum likelihood \citep{Ferraccioli_et_al2021}.
Related approaches address intensity estimation on linear networks \citep{McSwigganBaddeleyNair2017} and irregular domains \citep{BiscioLavancier2025}.
Our concern is the statistical difficulty imposed by the local geometry of a given domain.

\citet{BertinKlutchnikoffOuimet2025} introduced local polynomial density estimators on arbitrary known domains, together with a Goldenshluger--Lepski selection rule.
They established minimax optimality on simple domains for every positive smoothness level.
On polynomial sectors, their bounds match the minimax lower bounds for $0<s\leq1$, but deteriorate for $s>1$.
These planar domains are pinched at the origin between the horizontal axis and a power curve.
The deterioration comes from controlling the local Gram matrix through its smallest eigenvalue.
We replace this spectral bound by a pointwise quantity and obtain matching rates on polynomial sectors for every $s>0$.

\paragraph{Christoffel functions and pointwise risk}
Our main contribution is a variance bound expressed through the Christoffel function of the normalised local Lebesgue measure.
Its reciprocal is the diagonal value of the corresponding reproducing kernel at the estimation point.
For a scalar bandwidth $h$, let $V_h$ denote the volume of $\domain\cap(t+[-h,h]^d)$ and let $\delta_m(h)$ denote this Christoffel function at degree $m$.
The variance bound has geometric factor
\[
    \frac{1}{nV_h\delta_m(h)}.
\]
The local volume determines the amount of information available near the target.
The Christoffel function measures the additional cost of polynomial evaluation there.
This distinction matters even on polynomial sectors, where a spectral bound can overestimate the variance by a diverging power of the bandwidth.

The bias is controlled by weighted absolute moments of the equivalent kernel.
We give a general risk bound for vector bandwidths and sufficient conditions for controlling the kernel through convergence of rescaled moments.
These results allow the variance and bias to be studied separately.
They connect local polynomial density estimation with the classical theory of Christoffel functions \citep{Nevai1986,Totik2000,Simon2008}, which also has statistical applications in support inference and data analysis \citep{LasserrePauwels2019,LasserrePauwelsPutinar2022}.

\paragraph{Minimax rates and adaptation}
We formulate two intrinsic conditions for scalar bandwidths.
The power-law volume condition (\PL) requires $V_h\asymp h^{d+\nu}$ for some $\nu\geq0$.
The pointwise Bernstein--Markov condition (\PBM) keeps $\delta_m(h)$ bounded away from zero at each fixed degree.
Together, they yield the minimax rate
\[
    n^{-s/(2s+d+\nu)}
\]
for every $s>0$ over the pointwise polynomial approximation classes $\Sigma(s,L)$ defined in Section~\ref{sec:framework}.
The smoothness requirement is isotropic.
The volume exponent $\nu$ acts as an excess dimension, while \PBM prevents an additional loss from polynomial evaluation.
Polynomial volume growth alone does not ensure this stability.

The uniform Bernstein--Markov condition (\BM) controls polynomial evaluation throughout the local neighbourhood.
Under \PL and \BM, a Goldenshluger--Lepski procedure attains the upper bound
\[
    (\log n/n)^{s/(2s+d+\nu)}
\]
without knowledge of $s$ or $L$.
Its thresholds use pointwise and uniform kernel diagonals.
For each fixed smoothness level, the analysis requires \BM only at finitely many fixed degrees.
The procedure uses the known domain and a bound on its volume exponent.

An anisotropic tangent cone condition (\ATC) provides a sufficient geometric condition for both the minimax and adaptive results.
It describes convergence of local moments after a coordinate rescaling, together with boundedness of the rescaled neighbourhoods.
All polynomial sectors satisfy this condition.
For a sector bounded by a power curve of order $k$, the volume exponent is $k-1$.
The role of this exponent is related to that of a vanishing design density in pointwise regression \citep{Gaiffas2005,Gaiffas2007,guillou2011minimax}.
It differs from estimation on a submanifold, where the intrinsic dimension can replace the ambient dimension \citep{BickelLi2007,BerenfeldHoffmann2021}.

\paragraph{Singular domains and anisotropic bandwidths}
Two examples describe phenomena beyond these intrinsic conditions.
On an exponential cusp, the minimax rate is $(\log n)^{-s}$ for every $s>0$.
An explicit anisotropic bandwidth attains this rate.
For $s>1$, local polynomial estimators with scalar bandwidths and uniformly bounded degrees are suboptimal under every deterministic choice of degree and bandwidth.
Thus isotropic regularity does not make scalar bandwidths sufficient on a singular domain.
For each fixed degree, the proposed anisotropic bandwidth is independent of smoothness and attains the minimax rate throughout the corresponding smoothness range.

On a connected ladder, we identify the minimax rate for $0<s\leq2$ and give the exact lower and upper fluctuation envelopes of its logarithm.
For $0<s\leq1$, a local-constant estimator with a scalar bandwidth is minimax.
For $1<s\leq2$, an anisotropic bandwidth attains the full minimax profile.
In this latter range, the risk exceeds the reference polynomial rate by a factor tending to infinity along every sequence of sample sizes.
The general adaptive theorem does not cover either example.
The ladder bandwidth depends on the smoothness.

\paragraph{Organisation of the paper}
Section~\ref{sec:framework} introduces the framework, the estimators and the general variance and bias bounds, together with sufficient conditions based on rescaled moments.
Section~\ref{sec:conditions} formulates the intrinsic conditions for scalar bandwidths.
Section~\ref{sec:results} gives the minimax and adaptive results.
\revision{Sections~\ref{sec:polynomial.sectors}, \ref{sec:Dexp} and~\ref{sec:Dlad} study polynomial sectors, the exponential cusp and the connected ladder, respectively.}
\ifArxiv
Proofs and auxiliary results are collected in Appendices~\ref{sec:proofs}--\ref{sec:proof.Dlad}.
\else
Proofs and technical results are collected in the supplementary material.
\fi


\section{Framework}\label{sec:framework}

\subsection{Domain, neighbourhood and risk}\label{sec:domain}

Endow $\R^d$ with the sup norm $\|\cdot\|_\infty$.
Let $\domain$ be the closure of a nonempty open subset of $\R^d$, let $t\in\domain$ be a fixed estimation point, and fix $\hmax>0$.
Set $\Delta=[-1,1]^d$.
For a bandwidth vector $\bh=(h_1,\ldots,h_d)\in(0,\hmax]^d$, write
\[
    \pi({\bh}):=\prod_{j=1}^d h_j,
    \qquad {\bh}\Delta:=\prod_{j=1}^d[-h_j,h_j].
\]
Define the local neighbourhood of $t$ at bandwidth ${\bh}$ by
\begin{equation}\label{eq:Eh}
    \Escale{\bh}=(\domain-t)\cap {\bh}\Delta.
\end{equation}
For a scalar bandwidth $h$, we write $\bh=h\mathbf{1}=(h,\ldots,h)$
and use the subscript or argument $h$ in place of $\bh$.
In particular, $\Ehmax=(\domain-t)\cap\hmax\Delta$ remains the fixed
neighbourhood used to define local smoothness.
Every $\Escale{\bh}$ has positive volume and is contained in $\Ehmax$.
Define the local volume 
\begin{equation}\label{eq:local.volumes}
    \Vscale{\bh}:=\mathrm{Leb}(\Escale{\bh}).
\end{equation}
The quantity $n\Vscale{\bh}$ gives the geometric scale of the local sample size.
The domain is known, so all geometric quantities entering the procedures
are deterministic.
We suppress the dependence on the fixed objects $\domain,d,t$ and $\hmax$.
The symbols $\lesssim$, $\gtrsim$ and $\asymp$ denote comparisons up to
positive constants independent of the sample size and bandwidth.
Unless stated otherwise, these constants may depend on the fixed domain,
target, maximal bandwidth, polynomial degree and smoothness parameters.
Comparisons of risks over $\Sigma(s,L)$ are uniform in $f$.

Throughout, we assume $n\ge3$, so that $\log n>1$.
Let $X_1,\ldots,X_n$ be i.i.d.\ copies of a random variable $X$ supported on $\domain$ with density $f$ with respect to Lebesgue measure on $\R^d$.
Our goal is to estimate the value $f(t)$ at the fixed target $t$ from these observations.
The accuracy of an estimator $\tilde{f}$ of $f$ at the point $t$ is measured by the pointwise $L^2$ risk
\begin{equation}\label{eq:risk}
    \pwrisk(\tilde{f},f)
    = \Bigl(\EE_f^n\bigl[\{\tilde{f}(t)-f(t)\}^2\bigr]\Bigr)^{1/2}.
\end{equation}
Here $\EE_f^n$ denotes expectation under the joint distribution of the sample when the common density is $f$.
We write simply $\EE$ whenever no ambiguity can arise.
We use the usual $L^p(\mu)$ norms, writing $L^p(w)$ when
$\rd\mu=w\,\rd u$ and $L^p(A)$ for Lebesgue measure restricted to $A$.
The notation $\|g\|_{\infty,A}:=\sup_{u\in A}|g(u)|$ denotes
the pointwise supremum, including boundary points.
In particular, set $\finfty=\sup_{u\in\Ehmax}|f(t+u)|$.

\subsection{Smoothness classes}\label{sec:holder.class}

We encode local smoothness at $t$ through polynomial approximation on $\Ehmax$.
This formulation does not require $\Ehmax$ to contain a full neighbourhood of the origin and is therefore well suited to irregular domains.

For $s>0$, let $\llfloor s\rrfloor$ denote the greatest integer strictly less than $s$.
Thus $\llfloor s\rrfloor=s-1$ when $s\in\N$.
For $m\in\N_0=\N\cup\{0\}$, let $\Pm$ be the space of polynomials on $\R^d$ of total degree at most $m$.

    \begin{definition}\label{def:Sigma}
        \StmtDefSigma
    \end{definition}

We study densities $f\in\Sigma$, with risk bounds uniform over each
class $\Sigma(s,L)$. The defining condition concerns approximation at $t$.
It does not impose H\"older regularity between arbitrary pairs of nearby
points. If $\Ehmax$ is star-shaped and $f$ is classically H\"older-smooth
on $t+\Ehmax$, a Taylor polynomial at $t$ gives this approximation
\citep{BertinKlutchnikoffOuimet2025}.

Positivity and unit mass of $f$ control its approximating polynomial on
$\Ehmax$, and hence bound $\finfty$ uniformly over $\Sigma(s,L)$.
The bound depends on the geometry of this fixed neighbourhood.

\begin{lemma}[Uniform local boundedness]\label{lem:sup.bound}
    \StmtSupBound
\end{lemma}

    \subsection{Local polynomial estimators}\label{sec:lp.estimator}

For $m\in\N_0$, we use the monomial basis $\{\varphi_\alpha:\alpha\in\N_0^d,\,|\alpha|\leq m\}$ of $\Pm$, where $\varphi_\alpha(u)=u^\alpha=\prod_{j=1}^d u_j^{\alpha_j}$ and $|\alpha|=\sum_j\alpha_j$.
Its dimension is $D_m=\binom{m+d}{d}$.
We organise these monomials into a column vector $\Phi(u)\in\R^{D_m}$ under a fixed total order on multi-indices, placing $\varphi_0\equiv1$ first so that $\Phi(0)=e_1$.
Unlike \citet{BertinKlutchnikoffOuimet2025}, we use the unrescaled monomials $\Phi(u)$.
This change of basis leaves the projection estimator unchanged and simplifies the notation.

\revision{For ${\bh}\in(0,\hmax]^d$, let $\mu_{\bh}$ be the uniform probability measure on $\Escale{\bh}$. We also write $\mu_{\bh}(u)$ for its Lebesgue density
\[
    \mu_{\bh}(u)
    :=V_{\bh}^{-1}\1_{\Escale{\bh}}(u)
    =V_{\bh}^{-1}\1_\Delta(u/{\bh})\1_\domain(t+u),
\]
where $u/{\bh}=(u_1/h_1,\ldots,u_d/h_d)$.
It indeed satisfies $\int_{\R^d}\mu_{\bh}(u)\,\rd u=1$.}

For $\gamma=(m,{\bh})$ with $m\in\N_0$ and ${\bh}\in(0,\hmax]^d$, define the $D_m\times D_m$ Gram matrix
\revision{
\[
    \Bg = \int_{\R^d}\Phi(u)\Phi^\top(u)\,\mu_{\bh}(u)\,\rd u.
\]}
Since $\Vscale{\bh}>0$ and no non-trivial polynomial vanishes Lebesgue-a.e.\ on $\Escale{\bh}$, the matrix $\Bg$ is symmetric positive definite. The Christoffel--Darboux reproducing kernel of $\Pm$ in \revision{$L^2(\mu_{\bh})$} is
\begin{equation}\label{eq:repro.kernel}
   \revision{ k}_\gamma(u,v)=\Phi^\top(u)\Bg^{-1}\Phi(v),
    \qquad u,v\in\R^d.
\end{equation}
It satisfies \revision{the reproducing property} $\langle \revision{ k}_\gamma(\cdot,v),p\rangle_{L^2(\revision{ \mu}_{\bh})}=p(v)$
for every $p\in\Pm$. Note that when $\gamma=(m,{\bh})$, we write interchangeably $k_\gamma$ or $k_{m,{\bh}}$.

\revision{Proceeding as in \citet{BertinKlutchnikoffOuimet2025} by estimating the projection coefficients by empirical averages, we obtain the estimator}
\begin{equation}\label{eq:estimator.kernel}
    \hat f_\gamma(t)=\frac1n\sum_{i=1}^n
    \revision{ k}_\gamma(0,X_i-t)\mu_{\bh}(X_i-t).
\end{equation}

\noindent\revision{The equivalent kernel has value $k_\gamma(0,0)/V_{\bh}$ at the origin. The next section relates this diagonal to the Christoffel function.}

\subsection{Kernel identities}

\revision{We recall some identities for the Christoffel--Darboux reproducing kernel $k_\gamma$ for $\gamma=(m,\bh)$ with $\bh\in(0,h^*]^d$ and $m\in\N_0$. 
By the extremal characterisation of a reproducing kernel \citep[][]{Nevai1986,Simon2008}, for every $v\in\R^d$,
    \begin{equation*}
        k_{m,\bh}(v,v)
        =\sup_{0\neq P\in\Pm}\frac{P(v)^2}{\|P\|_{L^2(\mu_{\bh})}^2},
    \end{equation*}
    with equality attained at $P=k_{m,\bh}(\cdot,v)$.
    Indeed, this corresponds to equality case in the Cauchy--Schwarz inequality.
    In fact, we have for any $P\in \Pm$ and using that $k_{m,\bh}$ is the reproducing kernel of $\Pm$ in $L^2(\mu_{\bh})$,
    \[
        P(v)^2=\left|\langle P,k_{m,\bh}(\cdot,v)\rangle_{L^2(\mu_{\bh})}\right|^2\le \|P\|_{L^2(\mu_{\bh})}^2 \|k_{m,\bh}(\cdot,v)\|_{L^2(\mu_{\bh})}^2=\|P\|_{L^2(\mu_{\bh})}^2 k_{m,\bh}(v,v).
    \]
    Setting $v=0$ and taking reciprocals,
    \begin{equation*}
        \frac{1}{k_{m,\bh}(0,0)}
        =\inf_{\substack{P\in\Pm\\ P(0)=1}}\|P\|_{L^2(\mu_{\bh})}^2
    \end{equation*}
    with the infimum attained by the kernel section normalised to equal one at the origin.}
   
    The reproducing property also gives
    \begin{equation}
        \int_{\R^d}k_\gamma^2(0,u)\mu_{\bh}(u)\,\rd u
        =
        \bigl\langle k_\gamma(0,\cdot),k_\gamma(0,\cdot)\bigr\rangle_{L^2(\mu_{\bh})}
        =
        \kg.
        \label{eq:reproo}
    \end{equation}
    
    Additionally, by the reproducing property and Cauchy--Schwarz inequality for any $x\in\R^d$
     \begin{align}
       k^2_{\gamma}(0,x) =&\left(\int_{\R^d} k_{\gamma}(0,u) k_{\gamma}(u,x) \mu_{\bh}(u)\,\rd u\right)^2\nonumber\\
        \le &\left(\int_{\R^d} k_{\gamma}^2(0,u)  \mu_{\bh}(u)\,\rd u \right) \left(\int_{\R^d}  k_{\gamma}^2(u,x) \mu_{\bh}(u)\,\rd u\right)\nonumber\\
       \le & \ \kg k_{\gamma}(x,x).\label{eq:iden}
    \end{align}

\subsection{Christoffel function}\label{sec:christoffel}

\noindent\revision{The Christoffel function at the origin is the reciprocal of the kernel diagonal $k_\gamma(0,0)$.}
\begin{definition}[Christoffel function]\label{def:delta}
\StmtDefDelta
\end{definition}
    
\revision{The variational identities follow from the extremal characterisation of the reproducing kernel $k_\gamma$, see \citep{Nevai1986,Simon2008}. The finite-dimensional variational problem gives $0<\delta_m(\bh)\leq1$.
In particular, $\delta_0({\bh})=1$ for every ${\bh}$.}

In fact, the Christoffel function corresponds to the squared distance from the constant polynomial $1$ to $\Zm$ in $L^2(\mu_{\bh})$. The quantity $\delta_m(\bh)$ measures how well a polynomial vanishing at the origin can approximate the constant $1$ under the local measure.
If it stays bounded away from zero, polynomial evaluation at the target is stable.
By contrast, $\delta_m(\bh)\to0$ means that the local mass allows a polynomial to be close to $1$ on most of the local neighbourhood while vanishing at the target, thereby forcing extrapolation to the origin.

\subsection{General risk bound}\label{sec:risk.bound}

The variance is controlled by the local volume and the Christoffel function.
The bias depends on the absolute moments of the kernel $k_\gamma$.

\StmtBiasFactors

\begin{proposition}[General risk bound]\label{prop:risk.bound}
    \StmtRiskBound
\end{proposition}

If $f(t+u)=f(t)$ throughout $\Escale{\bh}$, reproduction of constants and
the squared-norm identity~\eqref{eq:reproo} give the exact variance
\[
    \Var\bigl(\hat f_\gamma(t)\bigr)
    =\frac1n\left\{
        \frac{f(t)}{V_{\bh}\,\delta_m(\bh)}-f(t)^2
    \right\}.
\]
For a fixed positive value of $f(t)$ and shrinking bandwidths, this
attains the order of the Christoffel variance bound.

\subsection{Moment conditions under linear rescaling}\label{sec:moment.stability}

The preceding risk bound involves the Christoffel function and the
absolute moments of the kernel. We give sufficient conditions for their
control through convergence of the local moments after an invertible
linear change of coordinates.
Since $\Pm$ consists of polynomials of total degree at most $m$,
composition with an invertible linear map gives another basis of the
same polynomial space. We use this change of basis to analyse the local
moment matrix.

\StmtRescaledSetup

\begin{proposition}[Stability through moments]\label{prop:moment.stability}
    \StmtMomentStability
\end{proposition}

These conclusions concern small bandwidths along the specified family.
Only moments up to order $2m$ enter the Gram matrix.
Their normalisation imposes no volume law.

\begin{corollary}[Convergence of rescaled domains]\label{cor:geometric.stability}
    \StmtGeometricStability
\end{corollary}

The limit domain may be separated from the origin.
Positive volume still makes polynomial evaluation at the origin a bounded
functional at every fixed degree.
Neither the matrices $A_{\bh}$ nor their determinants need follow power laws.


\section{Intrinsic geometric conditions}\label{sec:conditions}

We now restrict attention to scalar bandwidths $h$.
For the isotropic classes $\Sigma(s,L)$, we identify geometric conditions
under which scalar bandwidths yield the minimax rate.

\revision{For $m=\llfloor s\rrfloor$, we have
$\beta_m(s)=s$ and $\mathfrak L_{m,s,L}=L$ in Proposition~\ref{prop:risk.bound} and}
\begin{equation}\label{eq:scalar.risk.bound}
    \sup_{f\in\Sigma(s,L)}\pwrisk(\hat f_{(m,h)},f)
    \leq Lh^s\Lambda_m(h)
    +\left\{\frac{\FsL}{n\Vh\delta_m(h)}\right\}^{1/2}.
\end{equation}

The variance term separates local volume from the cost of polynomial
evaluation. We control these through \PL and \PBM, respectively.
The stronger condition \BM supplies the uniform control used by the
adaptive procedure. 
All conditions in this section concern scalar bandwidths.

\subsection{Power-law volume}\label{sec:PL}

The volume condition describes the amount of domain available in a
local neighbourhood through a power law.

\begin{definition}[Power-law volume conditions]\label{def:PL}
    \StmtDefPL
\end{definition}

The exponent $\nu$ in $\PL(\nu)$ is unique and is called the volume exponent. 
Its two sides have different statistical roles:
$\PL[-](\nu)$ provides enough observations for an upper bound,
whereas $\PL[+](\nu)$ makes local alternatives sufficiently
difficult for the minimax lower bound.

\subsection{Pointwise Bernstein--Markov property}\label{sec:PBM}

The pointwise Bernstein--Markov property controls polynomial evaluation
at the target uniformly over scalar bandwidths.

\begin{definition}[Pointwise Bernstein--Markov {\PBM[m]}]\label{def:PBM}
    \StmtDefPBM
\end{definition}

We write \PBM when \PBM[m] holds for every fixed $m\in\N_0$.

The Christoffel representation gives the following characterisation.

\begin{proposition}[Characterisation of {\PBM[m]}]\label{prop:PBM}
    \StmtPBM
\end{proposition}

\subsection{Uniform Bernstein--Markov property}\label{sec:BM}

The adaptive procedure also requires control over the whole local
neighbourhood. The uniform Bernstein--Markov property provides this control.
Such inequalities are classical in approximation theory
\citep{Bloom1997}. Here the bounds are uniform over bandwidths at each
fixed degree. We impose no growth condition on the constants as the
degree increases.

For $\gamma=(m,h)$, define the uniform kernel diagonal by
\begin{equation}\label{eq:Qg.def}
    \revision{\Qg:=\sup_{u\in\Eh}k_\gamma(u,u).}
\end{equation}

\begin{definition}[Bernstein--Markov {\BM[m]}]\label{def:BM}
    \StmtDefBM
\end{definition}

We write \BM when \BM[m] holds for every fixed $m\in\N_0$.

The uniform kernel diagonal characterises this property.

\begin{proposition}[Characterisation of {\BM[m]}]\label{prop:BM}
    \StmtBM
\end{proposition}

We now combine the volume and polynomial-stability conditions into the two regimes used in the main results.

\begin{definition}[Regular and uniformly regular points]\label{def:regular}
    \StmtDefRegular
\end{definition}

\subsection{The Anisotropic Tangent Cone condition}\label{sec:ATC}

The Anisotropic Tangent Cone (\ATC) condition gives sufficient
conditions to obtain regular and uniformly regular domains.
The anisotropy in the condition concerns the
coordinate-wise rescaling used to analyse the local domain, but the bandwidth remains scalar.
\begin{definition}[Anisotropic tangent cone condition $\ATC(\alpha)$]\label{def:ATC}
    \StmtDefATC
\end{definition}

We write \ATC when $\ATC(\alpha)$ holds for some exponent vector
$\alpha\in(0,\infty)^d$.

As we will see in the following, condition~\textup{(a)} implies \PL and \PBM.
The additional support condition~\textup{(b)} yields \BM.
It is automatic in directions with $\alpha_i\leq1$ and constrains
only the pinched directions $\alpha_i>1$.


\begin{proposition}[Moment condition implies regularity]\label{prop:ATC.PBM}
    \StmtATCPBM
\end{proposition}

\revision{To prove Proposition~\ref{prop:ATC.PBM}, we apply Proposition~\ref{prop:moment.stability} with $A_h=\operatorname{diag}(h^{\alpha_1},\ldots,h^{\alpha_d})$ and use condition~\textup{(a)} to obtain the limit of the rescaled Gram matrix $G_{m,h}$.}
The full \ATC condition provides the uniform control required by our adaptive result.

\begin{proposition}[\ATC implies uniform regularity]\label{prop:ATC.BM}
    \StmtATCBM
\end{proposition}


\revision{We say that $\ATC$ holds in exact form when  the transformed domain $\sigma_h(\Eh)=F$ where $F$ is a fixed
bounded set of positive volume.} Then $\ATC(\alpha)$ holds with
$c_\beta=\int_F v^\beta\,\rd v$.
This relates the condition to the model domains used in the study of
Christoffel functions on balls, simplices and planar corners
\citep{Xu1996,Xu1999,Kroo2015,Totik2010,PrymakUsoltseva2019}.

More generally, the measures
\revision{\[
    \rho_h(B):=h^{-(d+\nu)}\int_{\Eh}\1_B(\sigma_h(u))\,\rd u
\]}
are Lebesgue measure restricted to $\sigma_h(\Eh)$.
Condition~\textup{(a)} specifies the limits of their moments.
Under \ATC, they converge weakly to an absolutely continuous measure
$\rho_0$ of mass $c_0$.
The limit $\delta_{m,0}$ is the Christoffel function at the origin
for the probability measure $c_0^{-1}\rho_0$.

\section{Main results}\label{sec:results}

All bandwidths in this section are scalar.

\subsection{Known smoothness: minimax bounds}

The two bounds rely on different sides of the local volume law.
The estimator upper bound requires $\PL[-](\nu)$ and pointwise polynomial stability at degree $\llfloor s\rrfloor$, whereas the minimax lower bound requires only $\PL[+](\nu)$.

\begin{proposition}[Estimator upper bound]\label{prop:minimax.upper}
    \StmtMinimaxUpper
\end{proposition}

\begin{proposition}[Minimax lower bound]\label{prop:minimax.lower}
    \StmtMinimaxLower
\end{proposition}

Together, these bounds give the minimax rate under regularity.

\begin{theorem}[Minimax rate]\label{thm:minimax}
    \StmtMinimaxRate
\end{theorem}

Consequently, if $\domain$ is regular at $t$, with exponent $\nu$, the minimax rate is $n^{-s/(2s+d+\nu)}$ for every $s,L>0$.


\subsection{Adaptive procedure}\label{sec:adaptive}

We use a Goldenshluger--Lepski selector \citep{GoldenshlugerLepski2011,GoldenshlugerLepski2014} over a discrete family of the estimators~\eqref{eq:estimator.kernel}.
The procedure depends on the known domain and the data, but not on $s$ or $L$.
Its tuning uses an exponent $\nu$ for which $\PL[-](\nu)$ holds,
or a known upper bound on such an exponent.
Its thresholds are calibrated by the pointwise and uniform kernel diagonals $\kg$ and $\Qg$ introduced in Sections~\ref{sec:christoffel} and~\ref{sec:BM}.

\paragraph{Candidate grid}
For $\ell\in\N$, set
    \begin{equation}\label{eq:hl.ml}
        h_\ell = \hmax\,e^{-\ell},
        \qquad
        m_\ell = \left\lfloor\frac{\log n}{2\ell}\right\rfloor,
        \qquad
        \gamma_\ell=(m_\ell,h_\ell).
    \end{equation}
This coupling provides a sufficient degree at the scale associated
with each fixed smoothness level. For candidates with
$m<\llfloor s\rrfloor$, the bias is controlled through the exponent
$\beta_m(s)$ in Proposition~\ref{prop:risk.bound}.
We first retain scales with sufficient local volume:
    \begin{equation}\label{eq:Ln.def}
        \mathcal{L}_n
        = \bigl\{1\leq \ell \leq\lfloor\log n\rfloor
          \;:\; n\,\Vscale{h_\ell} \geq (\log n)^3\bigr\}.
    \end{equation}
To control geometric conditioning, recall 
   \revision{ \begin{equation}\label{eq:Qg}
        \Qg \;\geq\; \kg \;=\; \frac{1}{\deltam} \;\geq\; 1,
    \end{equation}}
for every $\gamma=(m,h)$. The candidate grid is
    \begin{equation}\label{eq:grid.new}
        \Gn = \bigl\{\gamma_\ell : \ell\in\mathcal{L}_n,\ Q_{\gamma_\ell}\leq n\bigr\}.
    \end{equation}

\paragraph{Thresholds}
For $\gamma=\gamma_\ell=(m,h)\in\Gn$, first define the summand envelope
   \revision{ \begin{equation}\label{eq:cgamma.def}
        c_\gamma = \frac{\sqrt{\kg\,\Qg}}{n\,V_h}.
    \end{equation}}
The reproducing-kernel Cauchy--Schwarz inequality \revision{\eqref{eq:iden} gives $k_\gamma(0,u)^2\leq\kg\Qg$} on $\Eh$.
Hence $c_\gamma$ bounds the absolute value of each summand, including
the factor $1/n$, in $\hat f_\gamma(t)$.

The associated empirical quadratic variance proxy is
    \begin{equation}\label{eq:hatv.new}
        \hat v_\gamma
        = \frac{1}{n^2}\sum_{i=1}^n
          \bigl\{\revision{k}_\gamma(0,X_i-t)\,\revision{\mu}_h(X_i-t)\bigr\}^2,
        \qquad
        v_\gamma = \EE[\hat v_\gamma].
    \end{equation}
Its expectation $v_\gamma$ upper-bounds the variance of $\hat f_\gamma(t)$.

Fix $\vartheta>0$ and set
    \begin{equation}\label{eq:eps.gamma}
        \varepsilon_\gamma := \vartheta\,c_\gamma^2\,n\Vh.
    \end{equation}
This tolerance controls the deviation of $\hat v_\gamma$ from $v_\gamma$ uniformly over the retained grid.

To account for geometric conditioning, set
    \[
        \Pi_\gamma := 2\log(e\vee Q_\gamma).
    \]
Fix $A>d+\nu$ using the available geometric bound, and define
    \begin{equation}\label{eq:pen.new}
        \mathrm{pen}(\gamma) := A\,(m+1)\,\ell + \Pi_\gamma.
    \end{equation}
The first term ensures summability over the candidate grid, while $\Pi_\gamma$ compensates for the geometric conditioning $Q_\gamma$.

With $r_\gamma(v,x)=\sqrt{2vx}+c_\gamma x$, define the data-driven threshold
    \[
        \hat{\mathbb{U}}_\gamma
        = r_\gamma\bigl(\hat v_\gamma+\varepsilon_\gamma,\ \mathrm{pen}(\gamma)\bigr).
    \]
\begin{remark}[Christoffel quantities in the thresholds]\label{rem:currency.adaptive}
\revision{The pointwise diagonal $\kg$ controls the variance, while
$\sqrt{\kg\Qg}$ controls the summands needed for concentration.
Under $\BM[m]$, both quantities and $\Pi_\gamma$ are uniformly bounded over bandwidths.}
\end{remark}

\paragraph{Comparisons and selection}
Suppose $\Gn$ is nonempty. We order candidates by their bandwidths:
$\gamma\preceq\gamma'$ if $h\leq h'$, equivalently $\ell\geq\ell'$.
Write $\gamma\vee\gamma'$ for the candidate with the larger bandwidth.
Define
    \begin{equation}\label{eq:Agamma.new}
        \hat A_\gamma
        = \max_{\gamma'\in\Gn}
          \Bigl\{\bigl|\hat f_{\gamma\vee\gamma'}(t)-\hat f_{\gamma'}(t)\bigr|
                 -\hat{\mathbb{U}}_{\gamma\vee\gamma'}-\hat{\mathbb{U}}_{\gamma'}\Bigr\}_+.
    \end{equation}
The quantity $\hat A_\gamma$ compares estimators after subtracting their
stochastic thresholds. We select
    \begin{equation}\label{eq:gammahat.new}
        \hat\gamma \in \argmin_{\gamma\in\Gn}\bigl\{\hat A_\gamma+\hat{\mathbb{U}}_\gamma\bigr\},
        \qquad
        \hat f(t) := \hat f_{\hat\gamma}(t).
    \end{equation}
Ties are resolved according to a fixed deterministic rule.
If $\Gn$ is empty, set $\hat f(t)=0$.

\subsection{Oracle inequality and adaptive upper bound}\label{sec:adaptive.results}

The oracle inequality requires $\PL[-](\nu)$ and local boundedness
of the density. Uniform regularity then yields the adaptive bound.

To state the oracle inequality, for $\gamma\in\Gn$ set
\begin{equation}\label{eq:Bgamma.def}
    B(\gamma):=\max\{|b_{\gamma'}|:\gamma'\in\Gn,\
    \gamma'\preceq\gamma\}.
\end{equation}
Thus $B(\gamma)$ is the largest absolute pointwise bias among candidates at least as fine as $\gamma$.

For the analysis, we also use the deterministic variance envelope
\revision{\[
    v_\gamma^\star:=\frac{\finfty\kg}{nV_h}\;\geq\;v_\gamma,
\]}
where the inequality follows from the definition of $v_\gamma$ and
 the reproduction identity~\eqref{eq:reproo}.

\begin{theorem}[Oracle inequality]\label{thm:oracle}
    \StmtOracle
\end{theorem}

Under uniform regularity, the grid contains, for each fixed smoothness, a well-conditioned candidate at the appropriate bias--variance scale.

\begin{theorem}[Adaptive upper bound]\label{thm:adaptive}
    \StmtAdaptive
\end{theorem}

The same estimator satisfies this bound for every $s,L>0$.
For each fixed $s$, the proof uses $\BM[m]$ only for the comparison
candidate and the finer candidates, whose degrees are bounded by a
constant depending on $s,d$ and $\nu$.
The logarithmic loss relative to Theorem~\ref{thm:minimax} is familiar
from pointwise adaptation \citep{Lepski1991adaptive,K-MMS-2014}.
We do not establish a matching adaptive lower bound.


\section{Polynomial sectors}\label{sec:polynomial.sectors}

Polynomial sectors provide a basic family of uniformly regular domains.
For $k>1$, let
    \begin{equation*}
        \Dk := \bigl\{(x,y)\in\R^2 : 0\leq x\leq 1,\;0\leq y\leq x^k\bigr\},
    \end{equation*}
as introduced by \citet{BertinKlutchnikoffOuimet2025}, and consider estimation at $t=(0,0)$.
For $\hmax\leq1$ and $\alpha=(1,k)$, one has
$\sigma_h(\Eh)=\Dk$ for every $h\in(0,\hmax]$.
Thus \ATC holds in exact form and $\nu=k-1$.
The dilation transforms $\tmu$ into the uniform probability measure
on $\Dk$. Composition with it preserves total polynomial degree.
Hence, for each fixed $m$, both $\delta_m(h)$ and $\Lambda_m(h)$ are
independent of $h$.
Consequently, Theorem~\ref{thm:minimax} shows that $\hat f_{(\llfloor s\rrfloor,h_n)}$, with $h_n=\hmax\wedge n^{-1/(2s+k+1)}$, attains the minimax rate
    \begin{equation*}
        \inf_{\hat{f}}\,\sup_{f\in\Sigma(s,L)}\pwrisk(\hat{f},f)
        \;\asymp\;
        n^{-s/(2s+k+1)},
        \qquad s,L>0.
    \end{equation*}
Since the regularity is uniform, Theorem~\ref{thm:adaptive} also gives the adaptive upper bound
\[
    \sup_{f\in\Sigma(s,L)}\pwrisk(\hat f,f)
    \;\lesssim\;
    \left(\frac{\log n}{n}\right)^{s/(2s+k+1)},
    \qquad s,L>0,
\]
without knowledge of $s$ or $L$.

The spectral analysis of \citet{BertinKlutchnikoffOuimet2025}
established minimax optimality only for $0<s\leq1$.
The Christoffel bound establishes it for every $s>0$ using the same
family of local polynomial estimators. The difference is seen by
comparing the geometric quantities entering the two risk bounds.
\revision{To obtain the spectral order, fix} $\gamma=(m,h)$ and define
\revision{\[
    \Bgiso:=\int_{\R^d}\frac{1}{h^d}\Phi(u/h)\Phi^\top(u/h)\1_{\Escale{\bh}}\,\rd u,
    \qquad \lamiso:=\lambda_{\min}(\Bgiso).
\]}
Let $G_m$ be the Gram matrix of the monomials under the uniform
probability measure on $\Dk$, and set
\[
    T_h=\operatorname{diag}\bigl(h^{(k-1)\beta_2}\bigr)_{|\beta|\leq m}.
\]
The change of variables $(x,y)=(hv_1,h^kv_2)$ gives
\[
    \Bgiso=\Wh T_hG_mT_h,
    \qquad \Wh=\frac{V_h}{h^2}=\frac{h^{k-1}}{k+1}.
\]
Since $h\leq1$, the smallest diagonal entry of $T_h$ is
$h^{m(k-1)}$. A lower bound using $\lambda_{\min}(G_m)>0$ and the
Rayleigh quotient for the monomial $y^m$ give
\[
    \Wh\lambda_{\min}(G_m)h^{2m(k-1)}
    \leq\lamiso
    \leq\Wh(G_m)_{(0,m),(0,m)}h^{2m(k-1)}.
\]
\revision{In \citet{BertinKlutchnikoffOuimet2025}, the risk of the polynomial estimator is bounded, up to a multiplicative constant, by $(W_h/\lamiso) (h^{2\beta_m(s)}+(nV_h)^{-1})^{1/2}$. The Christoffel bound, up to a multiplicative constant, is $(k_{m,h}(0,0))^{1/2}(h^{2\beta_m(s)}+(nV_h)^{-1})^{1/2}$.
By Proposition~\ref{prop:ATC.BM}, $k_{m,h}(0,0)$ and $Q_{m,h}$ are uniformly bounded over $0<h\leq h^*$ at each fixed degree. Moreover
        $\frac{W_h}{\lambda_{m,h}^{\mathrm{iso}}}
    \asymp h^{-2m(k-1)}.$}
For $m=0$, the two \revision{multiplicative terms} coincide, whereas for $m\geq1$ the spectral proxy is larger by the diverging factor $h^{-2m(k-1)}$.
Thus the spectral envelope produces an inflated upper bound for the risk of the corresponding estimator, whereas the Christoffel bound matches the minimax lower bound for $m\geq\llfloor s\rrfloor$ with a suitably chosen bandwidth.

\section{The exponential cusp: scalar and anisotropic bandwidths}\label{sec:Dexp}

On the exponential cusp, deterministic choices of scalar bandwidths
and bounded polynomial degrees are suboptimal for $s>1$.
An explicit anisotropic bandwidth attains the minimax rate over the
isotropic classes $\Sigma(s,L)$.

\subsection{Geometry and reference scale}

    Consider the domain
    \begin{equation}\label{eq:Dexp}
        \Dexp=\overline{\bigl\{(x,y)\in\R^2:0<x<1,\;0<y<e^{-1/x}\bigr\}},
        \qquad t:=(0,0).
    \end{equation}
    Fix $\hmax<1/4$.

    For scalar bandwidths $h$, up to boundary sets of Lebesgue measure zero,
    $\Eh=\{(x,y):0<x<h,\,0<y<e^{-1/x}\}$.
    Using $e^{-1/x}=x^2(e^{-1/x})'$, integration by parts gives
    $V_h=h^2e^{-1/h}-2\int_0^h xe^{-1/x}\,\rd x$. The remainder is at most $2hV_h$, so
    \begin{equation}\label{eq:volume.Dexp}
        \Vh=\int_0^h e^{-1/x}\,\rd x\sim h^2e^{-1/h}, \qquad \text{as $h\to0$}.
    \end{equation}
    Thus \PL fails for every exponent $\nu$, \ATC fails already at the
    zeroth moment, and the regular-domain rates do not apply.

    For $s>0$ and all sufficiently large $n$, define $h_n$ by
    \begin{equation}\label{eq:hn.Dexp}
        n\,h_n^{2s+2}\,e^{-1/h_n}=1.
    \end{equation}
    This equation has a unique positive solution, with $h_n\sim1/\log n$.
    It gives the scale used in the minimax lower bound below.

\subsection{Limits of bounded-degree estimators with scalar bandwidths}

    To understand what changes when $s>1$, consider first a fixed
    polynomial degree $m$.
    A local polynomial fit with a scalar bandwidth encounters an extrapolation
    problem.
    Indeed, the local measure concentrates near the relative abscissa
    $v_1=1$, away from the estimation point.
    For every fixed $m$, this concentration yields
    \begin{equation}\label{eq:Dexp.kernel.orders}
       \revision{ \delta^{-1}_m(h)= k_{m,h}(0,0)}\asymp h^{-2m},
        \qquad
        \Lambda_m(h)\asymp h^{-m}.
    \end{equation}

    The two orders in~\eqref{eq:Dexp.kernel.orders} allow to control the variance and the bias terms.
    The latter has two different sources.
    If $m<\llfloor s\rrfloor$, the horizontal monomial $x^{m+1}$
    is not reproduced by the estimator and contributes at order $h^{m+1}$.
    The approximation remainder contributes at most order $h^{s-m}$.
    Thus, for $m\le\llfloor s\rrfloor$, the relevant bias exponent is
    \begin{equation}\label{eq:alpha.m.Dexp}
        \alpha_m(s):=\min\{m+1,s-m\}.
    \end{equation}

    The standard deviation is bounded uniformly over $\Sigma(s,L)$
    by a constant times
    \[
        \left\{\frac{e^{1/h}}{nh^{2m+2}}\right\}^{1/2}.
    \]
    This order is attained under the uniform density on $\Dexp$.
    Optimising over the bandwidth therefore gives the best risk at degree
    $m$, of order $(\log n)^{-\alpha_m(s)}$.
    If $m>\llfloor s\rrfloor$, the worst-case approximation bias no longer
    tends to zero, and the best risk is of constant order.

    For a fixed integer $M\ge0$, define the family with bounded degrees and scalar bandwidths
    \begin{equation}\label{eq:single.scale.family.Dexp}
        \mathcal S_M
        :=\bigl\{\hat f_{(m,h)}:0\le m\le M,\ 0<h\le\hmax\bigr\}.
    \end{equation}
    The infimum below allows a deterministic pair $(m_n,h_n)$ to be
    chosen for each $n$, but not from the data.

\begin{theorem}[Best risk with bounded degrees and scalar bandwidths]\label{thm:Dexp.best.fixed.degree}
    \StmtDexpBest
\end{theorem}

    For every fixed $M\geq0$, the lower bound
    $(\log n)^{-\alpha_\star(s)}$ remains valid.
    Indeed, if $M<\llfloor s\rrfloor$, then
    $\mathcal S_M\subseteq\mathcal S_{\llfloor s\rrfloor}$.

    This lower bound concerns the pointwise approximation class
    $\Sigma(s,L)$ of Definition~\ref{def:Sigma}.
    Its proof uses remainders aligned with the sign of the equivalent
    kernel, subject to the envelope $L\|u\|_\infty^s$.
    Such remainders need not belong to a uniform H\"older class on a
    neighbourhood, so the theorem does not assert the same limitation
    for that smaller class.

\subsection{Anisotropic bandwidths and the minimax rate}

    For $0<r\leq\hmax/2$, choose the bandwidth vector and rescaling
    \begin{equation}\label{eq:Dexp.anisotropic.window}
        \bh(r)=(2r,e^{-1/r}),
        \qquad
        \sigma_{\bh(r)}(x,y)
        =\left(\frac{x}{2r},\frac{y}{e^{-1/r}}\right).
    \end{equation}
    Since $e^{-1/r}\leq2r$, one has
    $\Escale{\bh(r)}\subseteq\Ehmax$.
    Up to null boundary sets, its transformed domain is
    \[
        F_r:=\sigma_{\bh(r)}(\Escale{\bh(r)})
        =\left\{0<v_1<1,\quad
        0<v_2<\min\left(1,
        \exp\left\{\frac{1-(2v_1)^{-1}}{r}\right\}\right)\right\}.
    \]
    It is contained in $[0,1]^2$ and converges in measure to
    $F_0=[1/2,1]\times[0,1]$.
    More precisely,
    \begin{equation}\label{eq:Dexp.anisotropic.set.limit}
        \mathrm{Leb}(F_r\triangle F_0)
        \leq r/2.
    \end{equation}
    Corollary~\ref{cor:geometric.stability} therefore yields, for each
    fixed $m$,
    \begin{equation}\label{eq:Dexp.anisotropic.geometry}
        \begin{aligned}
        V_{\bh(r)}&=r e^{-1/r}(1+O(r)),\\
        \delta_m(\bh(r))&\longrightarrow\delta_{m,0}>0,
        \qquad \Lambda_m(\bh(r))\longrightarrow\Lambda_{m,0}<\infty.
        \end{aligned}
    \end{equation}
    The limits are determined by the uniform probability measure on
    $F_0$. Although the origin lies outside $F_0$, polynomial evaluation
    at the origin remains bounded at every fixed degree.

\begin{proposition}[Anisotropic upper bound]\label{prop:Dexp.anisotropic}
    \StmtDexpAnisotropic
\end{proposition}

\begin{theorem}[Minimax rate on $\Dexp$]\label{thm:Dexp}
    \StmtDexpMinimax
\end{theorem}

\begin{remark}[Adaptation on the exponential cusp]
    Fix $M\in\N_0$.
    The estimator $\hat f_{(M,\bh(r_n))}$ uses neither $s$ nor $L$.
    By Proposition~\ref{prop:Dexp.anisotropic} and
    Theorem~\ref{thm:Dexp}, it attains the minimax rate
    $(\log n)^{-s}$ for every fixed $0<s\leq M+1$ and $L>0$.
    Adaptation over this range incurs no additional rate loss.
\end{remark}

\begin{remark}[Isotropic smoothness and anisotropic bandwidths]
    For $s>1$, deterministic selection among scalar bandwidths with bounded
    degrees is suboptimal. An anisotropic bandwidth makes the same
    polynomial projection minimax over $\Sigma(s,L)$.
    The anisotropy of the bandwidth is imposed by the domain geometry.
\end{remark}


\section{An oscillating geometry: the connected ladder
\texorpdfstring{$\Dlad$}{D\_lad}}\label{sec:Dlad}

The connected ladder $\Dlad$ exhibits oscillations in the minimax rate.
Its normalised local volume fluctuates, while $\delta_m(h)$ tends to zero
along a sequence of scalar bandwidths for every $m\geq1$.

\subsection{Construction and geometric properties}

    Fix $\nu>0$.
    For $j\ge2$, set $\zeta_j=e^{-j^2}$,
    $\omega_j=\zeta_j^{1+\nu}$ and
    $\theta_j=\zeta_j/\zeta_{j-1}=e^{-(2j-1)}$.

    Let $I_j=[\zeta_j,\zeta_j+\omega_j]\times[0,1]$ be a vertical bar
    of width $\omega_j$.
    The bars $(I_j)_{j\ge2}$ are pairwise disjoint and accumulate on
    $\{0\}\times[0,1]$.
    We connect their lower edges by the exponential cusp $\Dexp$
    defined in Section~\ref{sec:Dexp}.
    The resulting domain is
    \begin{equation}\label{eq:Dlad}
        \Dlad
        :=\overline{\bigcup_{j\ge2}\mathring I_j\cup\mathring\Dexp},
        \qquad t:=(0,0).
    \end{equation}
    The cusp keeps the local neighbourhoods connected but is asymptotically
    negligible for the geometric quantities below.

    Write $k(h):=\min\{j\ge2:\zeta_j+\omega_j\le h\}$.

\begin{proposition}[Local-volume oscillations]\label{prop:Dlad.volume}
    \StmtDladVolume 
\end{proposition}

    Since \ATC implies \PL, the connected ladder cannot satisfy
    \ATC for any exponent vector.

\begin{proposition}[Failure of the pointwise Bernstein--Markov
condition]\label{prop:Dlad.PBM}
    \StmtDladPBM
\end{proposition}

    At the scales $h_k$, most of the local mass lies in bar $k$,
    near $u_1/h_k=1/2$.
    The polynomial $p_k(u)=2u_1/h_k$ therefore approaches $1$ in
   \revision{ $L^2(\mu_{h_k})$}, while $p_k(0)=0$.
    Since $p_k\in\Zm$ for every $m\geq1$, this explains the
    degeneration of $\delta_m(h_k)$ in Proposition~\ref{prop:Dlad.PBM}.
    It contrasts with the polynomial sectors $\Dk$, where the volume may
    decay with an arbitrarily large power while $\delta_m(h)$ remains
    bounded away from zero.

\subsection{Low smoothness: an implicit minimax rate}

    For the connected ladder, denote the minimax risk by
    \[
        \mathcal R_n(s,L):=\inf_{\hat f}\sup_{f\in\Sigma(s,L)}
        \pwrisk(\hat f,f).
    \]
    When $0<s\le1$, the local-constant estimator requires no polynomial
    extrapolation, and $\delta_0(h)=1$ identically.
    The local volume $\Vh$ determines the minimax scale through the
    implicit equation below.

\begin{theorem}[Implicit minimax rate and oscillations]\label{thm:Dlad.low}
    \StmtDladLow
\end{theorem}

    Along two subsequences of sample sizes, the minimax risk has,
    respectively, the baseline order $n^{-s/(2s+2+\nu)}$ and the inflated
    order
    $n^{-s/(2s+2+\nu)}\allowbreak
    \exp\{c_{s,\nu}^*
    \sqrt{\log n}(1+o(1))\}$.
    Consequently,
    the logarithmic volume exponent alone does not determine a
    non-oscillating minimax rate.

\subsection{Anisotropic bandwidths and full minimax oscillations}

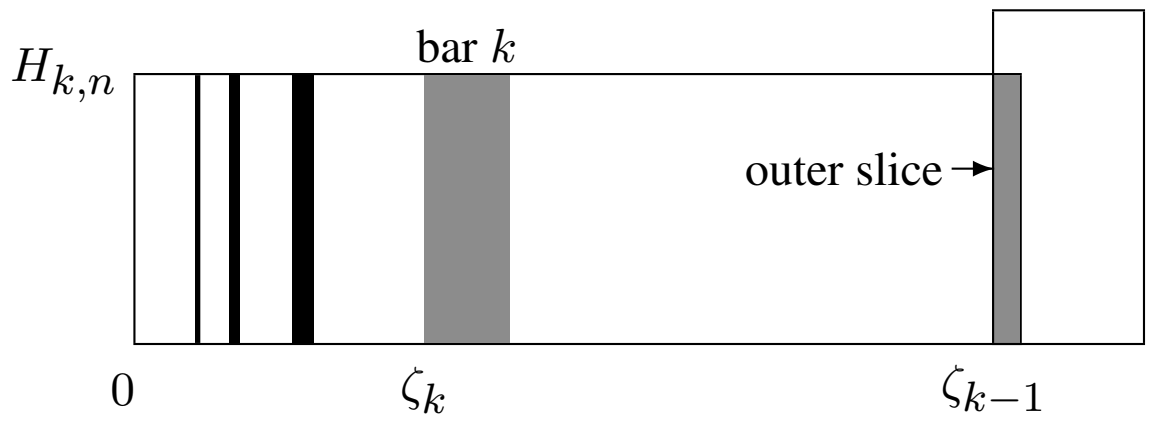
\begin{figure}[t]
    \centering
    \setlength{\unitlength}{0.8pt}
    \begin{picture}(290,116)
        \put(26,25){\framebox(208,63){}}
        \put(40,25){\rule{1pt}{50.4pt}}
        \put(48,25){\rule{2pt}{50.4pt}}
        \put(63,25){\rule{4pt}{50.4pt}}
        \put(94,25){\color[gray]{0.55}\rule{16pt}{50.4pt}}
        \put(228,25){\color[gray]{0.55}\rule{4.8pt}{50.4pt}}
        \put(228,25){\framebox(35,78){}}
        \put(23,14){\makebox(0,0){$0$}}
        \put(94,14){\makebox(0,0){$\zeta_k$}}
        \put(228,14){\makebox(0,0){$\zeta_{k-1}$}}
        \put(22,88){\makebox(0,0)[r]{$H_{k,n}$}}
        \put(104,95){\makebox(0,0){bar $k$}}
        \put(216,66){\makebox(0,0)[r]{outer slice}}
        \put(218,66){\vector(1,0){10}}
    \end{picture}
    \caption{Schematic of the anisotropic window for $H_{k,n}<1$.
    The outer slice has width $\ell_k$, compared with $\omega_k$ for the
    dominant bar. Its relative mass is of order
    $\ell_k/\omega_k=\zeta_k/\zeta_{k-1}$.
    The connecting cusp is omitted.}
    \label{fig:ladder.window}
\end{figure}

    For $1<s\leq2$, we use the local polynomial estimator of degree one
    with an anisotropic bandwidth. Fix a sufficiently large integer $k_0\geq3$.
    For $k\geq k_0$, set
    \begin{equation}\label{eq:Dlad.high.windows}
        \ell_k=\omega_k\frac{\zeta_k}{\zeta_{k-1}},
        \qquad H_{k,n}=(n\omega_k)^{-1/(2s+1)},
        \qquad \bh_{k,n}=(\zeta_{k-1}+\ell_k,H_{k,n}).
    \end{equation}
    Intersecting the bars with
    $[0,\zeta_{k-1}+\ell_k]\times[0,H_{k,n}]$ retains those with index
    at least $k$ and a slice of width $\ell_k$ from bar $k-1$,
    truncated at height $H_{k,n}\wedge1$.
    Figure~\ref{fig:ladder.window} shows the retained bars and the outer
    slice. The slice has relative mass of order $\zeta_k/\zeta_{k-1}$.

    The horizontal bias and the balance between vertical bias and variance
    lead to the quantities
    \begin{equation}\label{eq:Dlad.high.profile}
        \begin{aligned}
            B_k(s)&=\zeta_k\zeta_{k-1}^{s-1},
            &v_k(n,s)&=(n\omega_k)^{-s/(2s+1)},\\
            \mathfrak r_n(s)
            &=\inf_{k\geq k_0}\max\{B_k(s),v_k(n,s)\}.
        \end{aligned}
    \end{equation}
    Let $k_n$ be the smallest index attaining this infimum. We take
    \begin{equation}\label{eq:Dlad.high.estimator}
        \hat f_n^{\rm lad}:=\hat f_{(1,\bh_{k_n,n})}.
    \end{equation}

\begin{theorem}[Minimax rate and full oscillations on the ladder]
\label{thm:Dlad.anisotropic}
    \StmtDladAnisotropic
\end{theorem}

    At the selected index $k=k_n$, the slice of the outer bar has relative
    mass of order $\zeta_k/\zeta_{k-1}$.
    It stabilises affine evaluation at the origin and contributes at most
    order $B_k(s)$ to the weighted kernel moment. More precisely,
    \[
        \mathcal A_{1,s}(\bh_{k_n,n})
        \lesssim B_{k_n}(s)+H_{k_n,n}^s.
    \]
    This bound enters Proposition~\ref{prop:risk.bound} directly.
    A bound using only the maximum component of the bandwidth would
    lose this gain.

    The lower fluctuation envelope is strictly positive for $s>1$.
    Consequently,
    \[
        n^{s/(2s+2+\nu)}\mathcal R_n(s,L)\longrightarrow\infty.
    \]
    The baseline polynomial rate is therefore not attained along any
    subsequence, although it remains the logarithmic power of the risk.
    The choices of $k_n$ and $H_{k_n,n}$ depend on $s$.
    Theorem~\ref{thm:Dlad.anisotropic} does not assert adaptation to
    unknown smoothness on the ladder.


\ifArxiv
    \appendix
    \newcommand{\Restate}[3]{}
    \newcounter{proofpart}
\newcounter{proofstep}[proofpart]
\renewcommand{\theproofpart}{\alph{proofpart}}
\renewcommand{\theproofstep}{\roman{proofstep}}
\AtBeginEnvironment{proof}{%
    \setcounter{proofpart}{0}%
    \setcounter{proofstep}{0}%
}
\newcommand{\proofpart}[1]{%
    \stepcounter{proofpart}%
    \par\medskip\noindent
    \textbf{Part (\theproofpart). #1.}%
}
\newcommand{\proofstep}[1]{%
    \stepcounter{proofstep}%
    \par\medskip\noindent
    \textit{Step (\theproofstep). #1.}%
}

\section{Proofs of Section \emph{Framework}}\label{sec:proofs}

\subsection{Proof of Lemma~\ref{lem:sup.bound}}\label{sec:proof.sup.bound}

\Restate{relem}{lem:sup.bound}{\StmtSupBound}

\begin{proof}
    Let $f\in\Sigma(s,L)$ and let
    $q\in\mathcal P_{\llfloor s\rrfloor}$ be an associated polynomial.
    Since $\|u\|_\infty\le\hmax$ on $\Ehmax$, Definition~\ref{def:Sigma}
    and the non-negativity of $f$ give
    \[
        |q(u)|\le f(t+u)+L\hmax^s,
        \qquad u\in\Ehmax.
    \]
    Since $t+\Ehmax\subseteq\domain$ and $f$ is a probability density,
    \[
        \|q\|_{L^1(\Ehmax)}
        \le 1+L\hmax^s\,\Vscale{\hmax}.
    \]

    Because $\Vscale{\hmax}>0$, a polynomial that vanishes
    Lebesgue-almost everywhere on $\Ehmax$ is identically zero.
    Hence $\|\cdot\|_{L^1(\Ehmax)}$ is a norm on the finite-dimensional
    space $\mathcal P_{\llfloor s\rrfloor}$. Since
    the pointwise supremum over the compact set $\Ehmax$ is also a finite
    norm on this space,
    there exists $\Theta_{s,\hmax}<\infty$, depending only on
    $\llfloor s\rrfloor$, $d$ and $\Ehmax$, such that
    \[
        \sup_{u\in\Ehmax}|q(u)|
        \le\Theta_{s,\hmax}\|q\|_{L^1(\Ehmax)}.
    \]

    Finally, Definition~\ref{def:Sigma} also gives
    \[
        \finfty
        \le\sup_{u\in\Ehmax}|q(u)|+L\hmax^s
        \le
        \Theta_{s,\hmax}
        \bigl(1+L\hmax^s\Vscale{\hmax}\bigr)+L\hmax^s.
    \]
    The result follows by defining the last quantity as $\FsL$.
\end{proof}

\subsection{Proof of Proposition~\ref{prop:risk.bound}}\label{sec:proof.risk.bound}

\ifArxiv\else\StmtBiasFactors\fi

\Restate{reprop}{prop:risk.bound}{\StmtRiskBound}

\begin{proof}
   
\revision{Let  $m\in\N_0$, ${\bh}\in(0,\hmax]^d$ and $\gamma=(m,\bh)$.}
The bias--variance decomposition gives
    \[
        \pwrisk\bigl(\hat{f}_\gamma,f\bigr)^2
        =
        b_\gamma^2+\Var\bigl(\hat{f}_\gamma(t)\bigr).
    \]

     \proofpart{Absolute kernel moments}
    Reproduction of constants  gives
    \[
        1=\left|\int_{\Escale{\bh}}\revision{k}_\gamma(0,u)\revision{\mu}_{\bh}(u)\,\rd u\right|
        \leq\Lambda_m({\bh}).
    \]
    
    Moreover, Cauchy--Schwarz and \eqref{eq:reproo} implies
    
    \[ \Lambda_m({\bh})\leq \left(\int_{\Escale{\bh}}\revision{k}^2_\gamma(0,u)\revision{\mu}_{\bh}(u)\,\rd u\right)^{1/2}
        = \left(\revision{k}_\gamma(0,0)\right)^{1/2} =\delta_m({\bh})^{-1/2}.\]
    Since $\|u\|_\infty\leq\|{\bh}\|_\infty$ on $\Escale{\bh}$, this also proves
    $\mathcal A_{m,r}({\bh})\leq \|{\bh}\|_\infty^r\Lambda_m({\bh})$ for $r>0$.
    \proofpart{Variance bound}
    By the i.i.d.\ nature of the sample, and \revision{using \eqref{eq:reproo}, we obtain the following bounds for the variance term}
    \begin{align*}
        \Var\bigl(\hat f_\gamma(t)\bigr)
        &\leq \frac{1}{n}\,\EE\Bigl[\bigl\{\revision{k}_\gamma(0,X-t)\,\revision{\mu}_{\bh}(X-t)\bigr\}^2\Bigr] \\
        &= \frac{1}{n}\int_{\R^d} \revision{k}_\gamma^2(0,u) \revision{\mu}_{\bh}^2(u) f(t+u)\,\rd u \\
        &\leq \frac{\finfty}{n\revision{V_{\bh}}}\int_{\R^d}\revision{k}_\gamma^2(0,u) \revision{\mu}_{\bh}(u)\,\rd u \\
        &= \frac{\finfty \kg}{n\revision{V_{\bh}}}\\
        &=\frac{\finfty}{n\Vscale{\bh}\,\delta_m(\bh)}.
    \end{align*}

    \proofpart{Pointwise bias}
    Since $f\in\Sigma(s,L)$, there exists
    \[
        q=\sum_{|\alpha|\leq\llfloor s\rrfloor}q_\alpha\varphi_\alpha
        \in\mathcal P_{\llfloor s\rrfloor}
    \]
    such that $q(0)=f(t)$ and the remainder $r(u):=f(t+u)-q(u)$ satisfies
    \[
        |r(u)|\leq L\|u\|_\infty^s,
        \qquad u\in\Ehmax.
    \]
We define
    \[
        q_m^-
        :=
        \sum_{|\alpha|\leq\min(m,\llfloor s\rrfloor)}
        q_\alpha\varphi_\alpha,
        \qquad
        q_m^+:=q-q_m^-.
    \]
    Then $q_m^-\in\Pm$, $q_m^-(0)=q(0)=f(t)$, and $q_m^+\equiv0$ whenever $m\geq\llfloor s\rrfloor$.
    Define
    \begin{equation*}
        \widetilde r(u)
        :=
        f(t+u)-q_m^-(u)
        =
        r(u)+q_m^+(u).
    \end{equation*}

    We now bound $\widetilde r$ on $\Escale{\bh}$.
    If $m\geq\llfloor s\rrfloor$, then $q_m^+\equiv0$ and $\beta_m(s)=s$, so
    \[
        |\widetilde r(u)|=|r(u)|\leq L\|u\|_\infty^s,
        \qquad u\in\Escale{\bh}.
    \]
    Suppose instead that $m<\llfloor s\rrfloor$, so that $\beta_m(s)=m+1$.
For $u\in\Escale{\bh}$ and $m<|\alpha|\leq\llfloor s\rrfloor$,
    \[
        |\varphi_\alpha(u)|
        \leq
        \|u\|_\infty^{|\alpha|}
        \leq
        \max(1,\hmax)^{|\alpha|-(m+1)}\|u\|_\infty^{m+1}.
    \]
    Moreover, Lemma~\ref{lem:sup.bound} gives
    \[
        \|q\|_{\infty,\Ehmax}
        \leq
        \FsL+L\hmax^s.
    \]
    Since $\Vscale{\hmax}>0$, equivalence of norms on the finite-dimensional space $\mathcal P_{\llfloor s\rrfloor}$ therefore allows to bound all coefficients $q_\alpha$ uniformly over $f\in\Sigma(s,L)$.
    It follows that there exists $C_{m,s,L}>0$, independent of ${\bh}$ and of $f\in\Sigma(s,L)$, such that
    \[
        |q_m^+(u)|
        \leq
        C_{m,s,L}\|u\|_\infty^{m+1},
        \qquad u\in\Escale{\bh}.
    \]
    Since $|r(u)|\leq L\|u\|_\infty^s\leq L\max(1,\hmax)^{s-(m+1)}\|u\|_\infty^{m+1}$ on $\Escale{\bh}$, define
    \[
        \mathfrak L_{m,s,L}
        :=
        \begin{cases}
            L,
            & m\geq\llfloor s\rrfloor,\\
            L\max(1,\hmax)^{s-(m+1)}+C_{m,s,L},
            & m<\llfloor s\rrfloor.
        \end{cases}
    \]
    In both cases,
    \begin{equation}\label{eq:proof.r.tilde.bound}
        |\widetilde r(u)|
        \leq
        \mathfrak L_{m,s,L}\|u\|_\infty^{\beta_m(s)},
        \qquad u\in\Escale{\bh}.
    \end{equation}

    Since the kernel reproduces polynomials in $\Pm$,
    \[
        \int_{\Escale{\bh}}\revision{k}_\gamma(0,u)q_m^-(u)\revision{\mu}_{\bh}(u)\,\rd u
        =
        q_m^-(0)
        =
        f(t).
    \]
    Hence
    \begin{align}
        b_\gamma
        &= \int_{\Escale{\bh}} \revision{k}_\gamma(0,u) f(t+u) \revision{\mu}_{\bh}(u)\,\rd u - f(t) \nonumber \\
        &= \int_{\Escale{\bh}} \revision{k}_\gamma(0,u) \widetilde r(u) \revision{\mu}_{\bh}(u)\,\rd u.
        \label{eq:proof.bias.kernel}
    \end{align}
    Taking absolute values and applying~\eqref{eq:proof.r.tilde.bound} gives
    \begin{align*}
        |b_\gamma|
        &\leq\int_{\Escale{\bh}}|\revision{k}_\gamma(0,u)|\,|\widetilde r(u)|\revision{\mu}_{\bh}(u)\,\rd u\\
        &\leq\mathfrak L_{m,s,L}\mathcal A_{m,\beta_m(s)}({\bh})\\
        &\leq\mathfrak L_{m,s,L}\|{\bh}\|_\infty^{\beta_m(s)}\Lambda_m({\bh}).
    \end{align*}
    This proves~\eqref{eq:bias.bound.general}.
    
     \proofpart{Pointwise risk}
    Combining the variance bound and the bias bound yields
    \begin{align*}
        \pwrisk\bigl(\hat f_\gamma,f\bigr)^2
        &\leq b_\gamma^2 + \frac{\finfty \kg}{nV_{\bh}}\\
        &\leq\mathfrak L_{m,s,L}^2\mathcal A_{m,\beta_m(s)}({\bh})^2
        +\frac{\finfty}{n\Vscale{\bh}\,\delta_m(\bh)}.
    \end{align*}
    Finally, $\mathcal A_{m,\beta_m(s)}({\bh})\leq
    \|{\bh}\|_\infty^{\beta_m(s)}\Lambda_m({\bh})$ gives
    \begin{equation*}
        \pwrisk\bigl(\hat f_\gamma,f\bigr)^2
        \leq\bigl(\mathfrak L_{m,s,L}\Lambda_m({\bh})\|{\bh}\|_\infty^{\beta_m(s)}\bigr)^2
        +\frac{\finfty}{n\Vscale{\bh}\,\delta_m(\bh)},
    \end{equation*}
    which completes the proof.
\end{proof}

\subsection{Proof of Proposition~\ref{prop:moment.stability}}

\StmtRescaledSetup

\Restate{reprop}{prop:moment.stability}{\StmtMomentStability}

\begin{proof}
    \proofpart{Kernel and Christoffel limit}
    Composition with $A_{\bh}$ maps $\Pm$ bijectively onto itself, and
    $A_{\bh}$ sends $\varrho_{\bh}$ to \revision{$\mu_{\bh}$.}
    The reproducing property therefore gives that the  Christoffel--Darboux reproducing kernel of $\Pm$ in $L^2(\varrho_\bh)$ is 
    \begin{equation}\label{eq:general.kernel.rescaling}
        \Phi(v)^\top G_{m,{\bh}}^{-1}\Phi(w)=k_{m,{\bh}}(A_{{\bh}}v,A_{{\bh}}w).
    \end{equation}
    This identity holds for all $v,w\in\R^d$, including points outside
    the transformed domain.
    Since $G_{m,{\bh}}\to G_{m,0}$ as $\bh$ tends to $0$, and both $G_{m,\bh}$ and $G_{m,0}$ are positive definite their inverses converge.
    The monomials  $\varphi_\alpha$ are bounded on each compact set, which proves
    \eqref{eq:moment.kernel.limit}.
    Taking $v=w=0$ gives the limit of $\delta_m({\bh})$.
    Using that
    \[
        1\leq\Lambda_m({\bh})\leq\delta_m({\bh})^{-1/2}
    \]
     proves the remaining assertion of~\eqref{eq:moment.delta.limit}.

    \proofpart{Limit of the bias factor}

    \revision{We have 
    \begin{align*}
    \Lambda_m({\bh})
        &=\int_{\Escale{\bh}}|{k}_\gamma(0,u)|{\mu}_{\bh}(u)\,\rd u\\
        &=\int_{A_{\bh}^{-1}\Escale{\bh}}|{k}_\gamma(0,A_{\bh}v)|\rd \varrho_{\bh}(v)\\    
            &=\int_{A_{\bh}^{-1}\Escale{\bh}}|a_{\bh}(v)|\,\rd \varrho_{\bh}(v).
        \end{align*}
    where $a_{\bh}(v)=\Phi(0)^\top G_{m,{\bh}}^{-1}\Phi(v)=e_1^\top G_{m,{\bh}}^{-1}\Phi(v)$.}
If $m=0$, then $a_{\bh}=1$ and $\Lambda_0({\bh})=1$.
Suppose now $m\geq1$.
    The convergence of $G_{m,{\bh}}$ bounds its diagonal entries associated
    with $v_j^m$. Consequently,
    \[
        \sup_{{\bh}\text{ small}}\int\|v\|_\infty^{2m}\,\rd\varrho_{\bh}(v)
        \leq\sup_{{\bh}\text{ small}}\sum_{j=1}^d
        \int v_j^{2m}\,\rd\varrho_{\bh}(v)<\infty.
    \]
    The coefficients of $a_{\bh}$ are bounded, so
    $|a_{\bh}(v)|\leq C(1+\|v\|_\infty^m)$.
    For $R\geq1$ it follows that
    \[
        \sup_{{\bh}\text{ small}}
        \int_{\|v\|_\infty>R}|a_{\bh}(v)|\,\rd\varrho_{\bh}(v)
        \leq C R^{-m}.
    \]
    The same tail bound holds for $a_0$ under $\varrho_0$.
    On each compact set $a_{\bh}\to a_0$ uniformly.
    Using a continuous compactly supported cutoff, weak convergence
    therefore gives convergence of the truncated absolute integrals.
    Letting the cutoff radius tend to infinity proves
    \[
        \Lambda_m({\bh})=\int|a_{\bh}|\,\rd\varrho_{\bh}
        \longrightarrow\int|a_0|\,\rd\varrho_0,
    \]
    as asserted in~\eqref{eq:moment.Lambda.limit}.

    \proofpart{Uniform kernel control}
    Under part~\textup{(ii)}, every $u\in\Escale{\bh}$ has the form $u=A_{{\bh}}v$
    with $v$ in a fixed compact set $B$.
    Equation~\eqref{eq:general.kernel.rescaling} gives
    \[
        \sup_{u\in\Escale{\bh}}k_{m,{\bh}}(u,u)
        \leq\|G_{m,{\bh}}^{-1}\|_2\sup_{v\in B}\|\Phi(v)\|_2^2.
    \]
    The right-hand side is bounded for small bandwidths, proving
    \eqref{eq:moment.uniform}.
\end{proof}

\subsection{Proof of Corollary~\ref{cor:geometric.stability}}

\Restate{recor}{cor:geometric.stability}{\StmtGeometricStability}

\begin{proof}
    Put $v_{\bh}=\mathrm{Leb}(F_{\bh})$ and $v_0=\mathrm{Leb}(F_0)$.
    The change of variables $u=A_{{\bh}}v$ shows that $\varrho_{\bh}$ is uniform on
    $F_{\bh}$ and that $\Vscale{\bh}=|\det A_{\bh}|v_{\bh}$.
    Since $|v_{\bh}-v_0|\leq\varepsilon_{\bh}$, this proves
    \eqref{eq:geometric.volume}.
    Moreover,
    \[
        \int\left|\frac{\1_{F_{\bh}}(v)}{v_{\bh}}
        -\frac{\1_{F_0}(v)}{v_0}\right|\,\rd v
        =O(\varepsilon_{\bh}).
    \]
    Thus $\varrho_{\bh}$ converges to $\varrho_0$ in total variation.
    Every monomial is bounded on $B$, so
    $\|G_{m,{\bh}}-G_{m,0}\|_2=O(\varepsilon_{\bh})$.
    For each nonzero coefficient vector $c$, the polynomial
    $c^\top\Phi$ is nonzero and its zero set has Lebesgue measure zero.
    Hence
    \[
        c^\top G_{m,0}c
        =\frac{1}{v_0}\int_{F_0}|c^\top\Phi(v)|^2\,\rd v>0.
    \]
    All assumptions of Proposition~\ref{prop:moment.stability} follow.

    The resolvent identity gives
    \[
        G_{m,{\bh}}^{-1}-G_{m,0}^{-1}
        =G_{m,{\bh}}^{-1}(G_{m,0}-G_{m,{\bh}})G_{m,0}^{-1}
        =O(\varepsilon_{\bh}).
    \]
    This proves the claimed rates for the kernel and for $\delta_m$.
    With $a_{\bh}(v)=e_1^\top G_{m,{\bh}}^{-1}\Phi(v)$, one also has
    \begin{align*}
        \left|\int|a_{\bh}|\,\rd\varrho_{\bh}-\int|a_0|\,\rd\varrho_0\right|
        &\leq\sup_B|a_{\bh}-a_0|
        +\sup_B|a_0|\int|\rd\varrho_{\bh}-\rd\varrho_0|\\
        &=O(\varepsilon_{\bh}).
    \end{align*}
    This proves the rate for $\Lambda_m$ and completes the proof.
\end{proof}

\section{Proofs of Section \emph{Intrinsic geometric conditions}}\label{sec:proof.BM}

The Bernstein--Markov characterisations follow from the extremal representation of the reproducing kernel. The same representation expresses the bias factor through the extremal polynomial.

\subsection{Proof of Proposition~\ref{prop:PBM}}\label{sec:proof.PBM}

\Restate{reprop}{prop:PBM}{\StmtPBM}

\begin{proof}
   \revision{Since by definition $k_{m,h}(0,0)=\delta_m(h)^{-1}$ conditions~\textup{(ii)} and~\textup{(iii)} are equivalent.
    Moreover, the extremal characterisation gives
    \[
        k_{m,h}(0,0)
        =
        \sup_{0\ne p\in\mathcal P_m}
        \frac{|p(0)|^2}
             {\|p\|_{L^2(\mu_h)}^2}.
    \]
    Uniform boundedness in $h$ is exactly \PBM[m], proving the equivalence with condition~\textup{(i)}.}
\end{proof}

\subsection{Proof of Proposition~\ref{prop:BM}}\label{sec:proof.uniform.BM}

\Restate{reprop}{prop:BM}{\StmtBM}

\begin{proof}
    \revision{By the extremal characterisation we have, for every $u\in \Eh$,
    \[
        k_{m,h}(u,u)
        =\sup_{0\ne p\in\mathcal P_m}
          \frac{|p(u)|^2}{\|p\|_{L^2(\mu_h)}^2}.
    \]
    Taking the supremum over $u\in \Eh$ yields
    \[
        Q_{m,h} =\sup_{0\ne p\in\mathcal P_m}
          \frac{\|p\|_{\infty,\Eh}^2}
               {\|p\|_{L^2(\mu_h)}^2}.
    \]
    Uniform boundedness in $h$ is therefore equivalent to \BM[m].}
    Finally, $0\in \Eh$ implies $|p(0)|\leq\|p\|_{\infty,\Eh}$.
    Thus \BM[m] implies \PBM[m], and the same constant may be used: $C_m^{\mathrm{PBM}}\leq C_m^{\mathrm{BM}}$.
\end{proof}

\subsection{Proof of Proposition~\ref{prop:ATC.PBM}}\label{sec:proof.ATC.BM}

We first identify the limiting measures, then apply
Proposition~\ref{prop:moment.stability} to their normalised moments.

\begin{lemma}[Anisotropic tangent measure]\label{lem:tangent}
    Assume condition \textup{(a)} of $\ATC(\alpha)$, set
    $\nu=|\alpha|-d$, and define
  \revision{  $\rho_h(B):=h^{-(d+\nu)}\int\1_B(\sigma_h(u))\1_{\Eh}(u)\,\rd u$.}
    Then the family $(\rho_h)_{h\in(0,\hmax]}$ is tight, and every weak
    limit point $\rho_\star$ as $h\to0$ is absolutely continuous with
    respect to Lebesgue measure, with a density $g_\star\in[0,1]$, total
    mass $\rho_\star(\R^d)=c_0$, and moments
    $\int v^\beta\,\rd\rho_\star=c_\beta$ for every
    $\beta\in\N_0^d$. If condition \textup{(b)} also holds, the limit
    point is unique and $\rho_h\to\rho_\star$ weakly. In the exact case $\sigma_h(\Eh)=F$,
    $g_\star=\1_F$.
\end{lemma}

\begin{proof}
    \proofpart{Tightness}
    The change of variables $v=\sigma_h(u)$, with $\rd u=h^{|\alpha|}\rd v$ shows that $\rho_h$ has density
    $\1_{\sigma_h(\Eh)}$ with respect to the Lebesgue measure, since
    $\nu=|\alpha|-d$.
    For each $i\in\{1,\ldots,d\}$ and $q\in\N$, applying
    condition~\textup{(a)} to the multi-index $\beta$ with $2q$ in position
    $i$ and zeros elsewhere, we deduce that
    the even moments $\int v_i^{2q}\,\rd\rho_h$ converge as
    $h\to0$.
    In particular, there exist $h_1\in(0,\hmax]$ and constants
    $C_q<\infty$ such that
    \[
        \sup_{h\leq h_1}
        \max_{i\leq d}
        \int v_i^{2q}\,\rd\rho_h
        \leq
        C_q.
    \]
    By the Markov inequality,
    \begin{equation*}
        \sup_{h\leq h_1}\rho_h\bigl(\{v:\|v\|_\infty>\tau\}\bigr)
        \;\leq\;\sup_{h\leq h_1}\sum_{i=1}^d\tau^{-2q}\int v_i^{2q}\,\rd\rho_h
        \;\leq\;d\,C_q\,\tau^{-2q}
        \;\xrightarrow[\tau\to\infty]{}\;0,
    \end{equation*}
    so the family $(\rho_h)_{h\leq h_1}$ is tight.
    On the complementary range $h\in[h_1,\hmax]$, the sets
    $\sigma_h(\Eh)\subset\sigma_h(h\Delta)$ lie in a common compact set,
    and the masses
    \[
        \revision{\rho_h(\R^d)=h^{-(d+\nu)}V_h}
    \]
    are uniformly bounded.
    Hence the full family $(\rho_h)_{h\in(0,\hmax]}$ is tight.

    \proofpart{Limit moments}
    Moreover, since for each $\beta$ the integrals
    $\int v^{2\beta}\,\rd\rho_h$ are uniformly bounded for $h\leq h_1$,
    the monomials $v^\beta$ are uniformly integrable under
    $(\rho_h)_{h\leq h_1}$.
    By Prokhorov's theorem, weak limit points exist.
    If $\rho_{h_j}\to\rho_\star$ weakly along a sequence $h_j\to0$,
    uniform integrability allows us to pass to the limit in the moments:
    \[
        \int v^\beta\,\rd\rho_\star
        =
        \lim_{j\to\infty}\int v^\beta\,\rd\rho_{h_j}
        =
        c_\beta.
    \]
    This holds for every $\beta$, and the case $\beta=0$ gives
    $\rho_\star(\R^d)=c_0$.

    \proofpart{Absolute continuity}
    For absolute continuity, the Portmanteau theorem gives
    \[
        \rho_\star(U)
        \leq
        \liminf_{j\to\infty}\rho_{h_j}(U)
        \leq
        \mathrm{Leb}(U)
    \]
    for every open set $U$.
    By outer regularity, the same inequality holds for every Borel set.
    Thus $\rho_\star\ll\mathrm{Leb}$, with a density
    $g_\star\in[0,1]$ almost everywhere.

    \proofpart{Uniqueness and the exact case}
    Under condition~\textup{(b)}, every $\rho_h$, and hence every weak
    limit point, is supported in the compact set $[-R,R]^d$.
    A finite measure with compact support is determined by its moments.
    All weak limit points therefore coincide, since they share the moments
    $\{c_\beta\}$, and $\rho_h\to\rho_\star$ weakly as $h\to0$.
    Finally, in the exact case, $\sigma_h(\Eh)=F$ for every $h$.
    Hence $\rho_h=\mathrm{Leb}|_{F}$ does not depend on $h$, and
    $g_\star=\1_{F}$.
\end{proof}

\Restate{reprop}{prop:ATC.PBM}{\StmtATCPBM}

\begin{proof}
    Fix $m\in\N_0$ and take
    $A_h=\operatorname{diag}(h^{\alpha_1},\ldots,h^{\alpha_d})$
    in Proposition~\ref{prop:moment.stability}.
    Condition~\textup{(a)} gives \revision{$V_h\sim c_0h^{d+\nu}$} and
    \[
        G_{m,h}\longrightarrow G_{m,0}
        :=\bigl(c_{\beta+\beta'}/c_0\bigr)_{|\beta|,|\beta'|\leq m}, \qquad \text{as $h\rightarrow 0$}.
    \]
    By Lemma~\ref{lem:tangent}, this is the Gram matrix of the
    probability measure $\mu_\star=c_0^{-1}\rho_\star$, for any
    weak limit point $\rho_\star$ of the rescaled measures.
    This measure is absolutely continuous. The zero set of a nonzero
    polynomial has Lebesgue measure zero, so $G_{m,0}$ is positive
    definite.
    Proposition~\ref{prop:moment.stability} therefore yields
    \[
        \delta_m(h)\longrightarrow
        \delta_{m,0}:=(e_1^\top G_{m,0}^{-1}e_1)^{-1}>0,  \qquad \text{as $h\rightarrow 0$}.
    \]
    Since $\delta_m(h)\leq1$, the limit belongs to $(0,1]$.

     \revision{Then Condition \eqref{eq:PL.def} of $\PL(\nu)$ and (ii) of Proposition~\ref{prop:PBM} are satisfied for $0<h\le h_1$ for some fixed $h_1<h^*$.} To obtain the conditions on all of $(0,\hmax]$, note that
    \revision{$h^{-(d+\nu)}V_h$} and $\delta_m(h)$ are continuous and strictly
    positive for $h>0$.
    Indeed, the local moment integrals are continuous by dominated
    convergence, since the boundary of $h\Delta$ is Lebesgue-null,
    and the local Gram matrices are positive definite.
    Compactness extends the bounds near zero to every remaining
    bandwidth. Thus $\PL(\nu)$ and $\PBM[m]$ hold.
    Since $m$ was arbitrary, the domain is regular at $t$.
\end{proof}

\subsection{Proof of Proposition~\ref{prop:ATC.BM}}\label{sec:proof.ATC.uniform}

\Restate{reprop}{prop:ATC.BM}{\StmtATCBM}

\begin{proof}
    Fix $m\in\N_0$ and use the rescaling from the preceding proof.
    Condition~\textup{(b)} places $\sigma_h(\Eh)$ in a fixed compact
    set. Part~\textup{(ii)} of Proposition~\ref{prop:moment.stability}
    therefore gives
    \[
        Q_{m,h}\leq C_m^{\mathrm{ATC}}
    \]
    for all sufficiently small $h$.
    Together with $\delta_m(h)\leq1$ and $0\in\Eh$, this proves
   \revision{ \[
        1\leq k_{m,h}(0,0)
        \leq Q_{m,h}\leq C_m^{\mathrm{ATC}}.
    \]}
    On any compact bandwidth interval bounded away from zero,
    the inverse local Gram matrices are uniformly bounded.
    Since $\Eh\subseteq\Ehmax$, the monomials are uniformly bounded
    there, which gives a uniform bound for $Q_{m,h}$ after increasing $C_m^{\mathrm{ATC}}$ if needed.
    Hence $\sup_{0<h\leq\hmax}Q_{m,h}<\infty$.
    Proposition~\ref{prop:BM} yields $\BM[m]$, while
    Proposition~\ref{prop:ATC.PBM} gives $\PL(\nu)$.
    Since $m$ was arbitrary, the domain is uniformly regular at $t$.
\end{proof}

\section{Proofs of Section \emph{Main results}}\label{sec:supp.results}

\subsection{Proof of Proposition~\ref{prop:minimax.upper}}
\label{sec:proof.minimax.upper}

\Restate{reprop}{prop:minimax.upper}{\StmtMinimaxUpper}

\begin{proof}
    Set
    \[
        \bar h_n:=n^{-1/(2s+d+\nu)},
        \qquad
        h_n=\hmax\wedge\bar h_n,
        \qquad
        \gamma_n:=(m,h_n).
    \]
    By construction,
    \[
        \bar h_n^{2s}
        =
        \frac{1}{n\bar h_n^{d+\nu}}
        =
        n^{-2s/(2s+d+\nu)}.
    \]
    \revision{Proposition~\ref{prop:PBM} gives
    \[
        \inf_{0<h\leq\hmax}\delta_m(h)>0,
    \]
    while $\PL[-](\nu)$ gives $V_h\gtrsim h^{d+\nu}$ uniformly over $h\in(0,\hmax]$.
    Combining these facts with Lemma~\ref{lem:sup.bound} and  Proposition~\ref{prop:risk.bound} with $m=\llfloor s\rrfloor$, $\beta_m(s)=s$ and $\mathfrak L_{m,s,L}=L$, we obtain}
    \[
        \sup_{f\in\Sigma(s,L)}
        \pwrisk\bigl(\hat f_{\gamma_n},f\bigr)^2
        \lesssim
        h_n^{2s}
        +
        \frac{1}{n h_n^{d+\nu}}.
    \]
    Since $h_n\leq\bar h_n$, the bias term is bounded by $\bar h_n^{2s}$.
    For the variance term,
    \begin{align*}
        \frac{1}{n h_n^{d+\nu}}
        &=
        \max\left\{
            \frac{1}{n\hmax^{d+\nu}},
            \frac{1}{n\bar h_n^{d+\nu}}
        \right\} \\
        &\leq
        \max\left\{(\hmax)^{-(d+\nu)},1\right\}
        n^{-2s/(2s+d+\nu)},
    \end{align*}
    where we used $n^{-1}\leq n^{-2s/(2s+d+\nu)}$.
    Therefore,
    \[
        \sup_{f\in\Sigma(s,L)}
        \pwrisk\bigl(\hat f_{\gamma_n},f\bigr)^2
        \lesssim
        n^{-2s/(2s+d+\nu)}.
    \]
    Taking square roots proves the result for every $n\geq3$.
\end{proof}

\subsection{Proof of Proposition~\ref{prop:minimax.lower}}
\label{sec:proof.minimax}

\Restate{reprop}{prop:minimax.lower}{\StmtMinimaxLower}

\begin{proof}
    We use Le~Cam's two-point method with a bump function perturbation.
    By $\PL[+](\nu)$,
    \[
        \revision{\Vh}
        \leq C_\nu^{\mathrm{PL}}h^{d+\nu},
        \qquad 0<h\leq\hmax.
    \]
    The perturbation constructed below has amplitude of order $h^s$,
    so its $n$-sample Kullback--Leibler divergence will be bounded by a
    constant times $nh^{2s}\Vh$, hence by a constant times
    $nh^{2s+d+\nu}$. This motivates setting
    \[
        \bar h_n:=n^{-1/(2s+d+\nu)},
    \]
    for which $n\bar h_n^{2s+d+\nu}=1$.
    In particular, $\Vh\to0$ as $h\to0$.
    Fix $h_0\in(0,\hmax]$ such that
    \[
        \Vscale{h_0}\leq\frac{\Vscale{\hmax}}{2},
    \]
    and set
    \[
        h_n:=h_0\wedge\bar h_n.
    \]
    In particular, $nh_n^{2s+d+\nu}\leq1$. Moreover, since
    $\bar h_n\leq\bar h_3$ for every $n\geq3$,
    \begin{equation}\label{eq:proof.lower.hn}
        \kappa_0\bar h_n\leq h_n\leq\bar h_n,
        \qquad
        \kappa_0:=1\wedge\frac{h_0}{\bar h_3}>0.
    \end{equation}

    \proofpart{Construction}
    Choose $\psi\in C_c^\infty(\R^d)$ such that $\psi\geq0$ and
    $\psi=1$ on $t+\hmax\Delta$, and define
    \[
        f_0(x):=\frac{\psi(x)}{\int_\domain\psi(y)\,\rd y},
        \qquad x\in\domain.
    \]
    The denominator is finite and strictly positive, and $f_0$ is
    constant on $t+\Ehmax$, with $f_0(x)=f_0(t)>0$ there.

    Choose $\phi_0\in C_c^\infty(\Delta)$ such that
    $0\leq\phi_0\leq1$, $\phi_0(0)=1$, and $\phi_0$ is constant in a
    neighbourhood of the origin. There exists $C_\phi>0$ such that
    \[
        |\phi_0(v)-1|
        \leq C_\phi\|v\|_\infty^s,
        \qquad v\in\Delta.
    \]
    Set $\phi_{h_n}(u):=\phi_0(u/h_n)$ and
    \[
        c_n
        :=
        \frac{1}{\Vscale{\hmax}}
        \int_{\Escale{h_n}}\phi_{h_n}(u)\,\rd u.
    \]
    Since $0\leq\phi_0\leq1$ and $h_n\leq h_0$,
    \begin{equation}\label{eq:proof.lower.cn}
        0\leq c_n
        \leq\frac{\Vscale{h_n}}{\Vscale{\hmax}}
        \leq\frac12.
    \end{equation}

    Let
    \[
        0<c_1\leq
        \min\left\{
            \frac{L}{C_\phi},\,
            L,\,
            \frac{f_0(t)}{h_0^s},\,
            \left(
                \frac{f_0(t)}
                     {2C_\nu^{\mathrm{PL}}\|\phi_0\|_\infty^2}
            \right)^{1/2}
        \right\},
        \qquad
        \rho_n:=c_1h_n^s,
    \]
    and define, for $x\in\domain$,
    \begin{equation}\label{eq:f1}
        f_1(x)
        :=
        f_0(x)
        +
        \rho_n\left\{
            \phi_{h_n}(x-t)\1_{t+\Escale{h_n}}(x)
            -
            c_n\1_{t+\Ehmax}(x)
        \right\}.
    \end{equation}

    \proofstep{Normalisation}
    By a change of variables and the definition of $c_n$,
    \begin{align*}
        \int_\domain f_1(x)\,\rd x
        &=
        \int_\domain f_0(x)\,\rd x
        +
        \rho_n\left\{
            \int_{\Escale{h_n}}\phi_{h_n}(u)\,\rd u
            -
            c_n\Vscale{\hmax}
        \right\} \\
        &=1.
    \end{align*}

    \proofstep{Non-negativity}
    On $t+\Ehmax$, the non-negativity of $\phi_0$ and the choice of
    $c_1$ give
    \[
        f_1(x)
        \geq
        f_0(t)-\rho_nc_n
        \geq
        f_0(t)-\frac{c_1h_0^s}{2}
        \geq
        \frac{f_0(t)}{2}>0.
    \]
    Outside $t+\Ehmax$, $f_1=f_0\geq0$. Thus $f_1$ is a probability
    density on $\domain$.

    \proofstep{Signal}
    Since $\phi_{h_n}(0)=1$,~\eqref{eq:proof.lower.cn} gives
    \begin{equation}\label{eq:proof.lower.signal}
        f_1(t)-f_0(t)
        =
        \rho_n(1-c_n)
        \geq
        \frac{\rho_n}{2}.
    \end{equation}

    \proofstep{Membership in $\Sigma(s,L)$}
    Since $f_0$ is constant on $t+\Ehmax$, it belongs to
    $\Sigma(s,L)$ with approximating polynomial $q_0\equiv f_0(t)$.
    For $f_1$, take the constant polynomial
    \[
        q(u)\equiv f_0(t)+\rho_n(1-c_n).
    \]
    If $u\in\Escale{h_n}$, then
    \[
        |f_1(t+u)-q(u)|
        =
        \rho_n|\phi_0(u/h_n)-1|
        \leq
        c_1C_\phi\|u\|_\infty^s
        \leq
        L\|u\|_\infty^s.
    \]
    If $u\in\Ehmax\setminus\Escale{h_n}$, then
    $\|u\|_\infty>h_n$ and
    \[
        |f_1(t+u)-q(u)|
        =
        \rho_n
        \leq
        c_1h_n^s
        \leq
        L\|u\|_\infty^s.
    \]
    Hence $f_0,f_1\in\Sigma(s,L)$.

    \proofstep{Kullback--Leibler divergence}
    Set
    \[
        g_n(u)
        :=
        \phi_{h_n}(u)\1_{\Escale{h_n}}(u)-c_n,
        \qquad u\in\Ehmax.
    \]
    Since $c_n$ is the mean of
    $\phi_{h_n}\1_{\Escale{h_n}}$ over $\Ehmax$,
    \begin{align*}
        \int_{\Ehmax}g_n(u)^2\,\rd u
        &=
        \int_{\Escale{h_n}}\phi_{h_n}(u)^2\,\rd u
        -
        c_n^2\Vscale{\hmax} \\
        &\leq
        \|\phi_0\|_\infty^2\Vscale{h_n}.
    \end{align*}
    Since $f_0=f_0(t)>0$ on $t+\Ehmax$ and $f_1=f_0$ outside this
    set, the standard $\chi^2$ upper bound yields
    \begin{align}
        D_{\mathrm{KL}}(f_1\,\|\,f_0)
        &\leq
        \int_{t+\Ehmax}\frac{(f_1-f_0)^2}{f_0}\,\rd x \notag\\
        &\leq
        \frac{\rho_n^2}{f_0(t)}
        \|\phi_0\|_\infty^2 \Vscale{h_n}.
        \label{eq:KL.bound}
    \end{align}
    By $\PL[+](\nu)$, the choice of $c_1$, and
    $nh_n^{2s+d+\nu}\leq1$,
    \[
        nD_{\mathrm{KL}}(f_1\,\|\,f_0)
        \leq
        \frac{C_\nu^{\mathrm{PL}}\|\phi_0\|_\infty^2}{f_0(t)}
        c_1^2nh_n^{2s+d+\nu}
        \leq\frac12.
    \]

    \proofpart{Conclusion}
    By Pinsker's inequality, Le~Cam's two-point lemma,
    and~\eqref{eq:proof.lower.signal},
    \begin{align*}
        \inf_{\hat f}\sup_{f\in\Sigma(s,L)}
        \pwrisk(\hat f,f)^2
        &\geq
        \frac{\rho_n^2}{16}
        \left(
            1-\sqrt{\frac n2D_{\mathrm{KL}}(f_1\,\|\,f_0)}
        \right)^2 \\
        &\geq
        \frac{\rho_n^2}{64}
        \geq
        \frac{c_1^2\kappa_0^{2s}}{64}
        n^{-2s/(2s+d+\nu)},
    \end{align*}
    where the last inequality follows from~\eqref{eq:proof.lower.hn}.
    Taking square roots proves the result for every $n\geq3$.
\end{proof}

\subsection{Proof of Theorem~\ref{thm:minimax}}
\label{sec:proof.minimax.rate}

\Restate{rethm}{thm:minimax}{\StmtMinimaxRate}

\begin{proof}
    Regularity of degree $\llfloor s\rrfloor$ means that
    \PBM[\llfloor s\rrfloor] and $\PL(\nu)$ hold at $t$.
    In particular, the assumptions of both
    Propositions~\ref{prop:minimax.upper} and~\ref{prop:minimax.lower}
    are satisfied. The former gives the upper bound, attained by
    $\hat f_{(\llfloor s\rrfloor,h_n)}$, while the latter gives the
    matching minimax lower bound.
\end{proof}

\subsection{Proof of Theorem~\ref{thm:oracle}}
\label{sec:proof.adaptive}

The proof of Theorem~\ref{thm:oracle} relies on two concentration
inequalities and a summability lemma, which we establish first.

\begin{lemma}[Concentration of $\hat v_\gamma$]\label{lem:var.concentration}
    Set
    \[
        \kappa_\vartheta
        :=
        \frac{\vartheta^2}
             {2\bigl(\finfty+\vartheta/3\bigr)}
        >0.
    \]
    Then, for every $\gamma=(m,h)\in\Gn$,
    \begin{equation*}
        \PP\bigl(|\hat v_\gamma-v_\gamma|>\varepsilon_\gamma\bigr)
        \;\leq\; 2\exp\bigl\{-\kappa_\vartheta\,n\,\revision{V_h}\bigr\}
        \;\leq\; 2\exp\bigl\{-\kappa_\vartheta(\log n)^3\bigr\}.
    \end{equation*}
\end{lemma}

\begin{proof}
    Write $\hat v_\gamma=\sum_{i=1}^n Z_i$, where
    \[
        Z_i=n^{-2}\{\revision{k}_\gamma(0,X_i-t)\revision{\mu}_h(X_i-t)\}^2.
    \]
   \revision{The variables $Z_i$ are independent, and   using that $k_{\gamma}^2(0,u)\le k_{\gamma}(0,0)Q_\gamma$ for $u\in \Eh$ we derive that $0\leq Z_i\leq c_\gamma^2$.
    Hence
    \[
        |Z_i-\EE Z_i|\leq c_\gamma^2.
    \]
The second-moment calculation used to prove part~(b) of Proposition~\ref{prop:risk.bound} gives
    \[v_\gamma=n\EE Z_1\leq\frac{\finfty\kg}{nV_h}=v_\gamma^\star\],
    and then
    \[\sum_{i=1}^n\Var(Z_i)
        \leq c_\gamma^2v_\gamma\leq c_\gamma^2v_\gamma^\star.\]}
Moreover using that $\Qg\ge 1$, we have
    \[
        v_\gamma^\star=\frac{\finfty\kg}{n\revision{V_h}}
        \leq\finfty c_\gamma^2n\Vh.
    \]

    Bernstein's inequality \citep[Theorem~2.10]{BoucheronLugosiMassart2013},
    with $\varepsilon_\gamma=\vartheta c_\gamma^2n\Vh$, gives
    \begin{align*}
        \PP(|\hat v_\gamma-v_\gamma|>\varepsilon_\gamma)
        &\leq2\exp\left\{-\frac{\varepsilon_\gamma^2}
        {2c_\gamma^2(v_\gamma^\star+\varepsilon_\gamma/3)}\right\}\\
        &\leq2\exp\{-\kappa_\vartheta n\Vh\}.
    \end{align*}
    The grid restriction $n\Vh\geq(\log n)^3$ proves the second bound.
\end{proof}

\begin{lemma}[Bernstein concentration of $\hat f_\gamma(t)$]\label{lem:bernstein.fhat}
    For every $\gamma\in\Gn$ and $x>0$,
    \begin{equation*}
        \PP\Bigl(\bigl|\hat f_\gamma(t)-\EE\hat f_\gamma(t)\bigr|>r_\gamma(v_\gamma,x)\Bigr)\;\leq\;2e^{-x}.
    \end{equation*}
\end{lemma}

\begin{proof}
    Write $\hat f_\gamma(t)=\sum_{i=1}^nY_i$, with
    $Y_i=n^{-1}\revision{k}_\gamma(0,X_i-t)\revision{\mu}_h(X_i-t)$.
   We have
    $|Y_i|\leq c_\gamma$, so $|Y_i-\EE Y_i|\leq2c_\gamma$.
    Moreover,
    \[
        \sum_{i=1}^n\Var(Y_i)\leq n\EE Y_1^2=v_\gamma.
    \]
    Bernstein's inequality \citep[Theorem~2.10]{BoucheronLugosiMassart2013}
    therefore gives
    \[
        \PP\left(|\hat f_\gamma(t)-\EE\hat f_\gamma(t)|
        >\sqrt{2v_\gamma x}+2c_\gamma x/3\right)\leq2e^{-x}.
    \]
    The claimed bound follows since
    $r_\gamma(v_\gamma,x)\geq\sqrt{2v_\gamma x}+2c_\gamma x/3$.
\end{proof}

\begin{lemma}[Summability over the grid]\label{lem:sums.adaptive}
    Suppose that $A>d+\nu$. There exists a constant $C_{\mathrm{sum}}>0$,
    depending only on $\vartheta,A,d,\nu,\hmax,\finfty$ and the
    constant $c_\nu^{\mathrm{PL}}$ in $\PL[-](\nu)$, such that
    \[
        \sum_{\gamma\in\Gn}
        \bigl(v_\gamma^\star+\varepsilon_\gamma+c_\gamma^2\bigr)
        e^{-\mathrm{pen}(\gamma)}
        \leq
        \frac{C_{\mathrm{sum}}}{n}.
    \]
\end{lemma}

\begin{proof}
    Fix $\gamma=(m,h)\in\Gn$.
    Using
    $n\Vh\geq1$, we obtain
    \[
        v_\gamma^\star\leq\frac{\finfty Q_\gamma^2}{n\Vh},
        \qquad
        \varepsilon_\gamma\leq\frac{\vartheta Q_\gamma^2}{n\Vh},
        \qquad
        c_\gamma^2\leq\frac{Q_\gamma^2}{(n\Vh)^2}
        \leq\frac{Q_\gamma^2}{n\Vh}.
    \]
    Set $C_\vartheta=\finfty+\vartheta+1$.
    For $\gamma=\gamma_\ell$, the penalty satisfies
    \[
        e^{-\mathrm{pen}(\gamma)}
        \leq e^{-A\ell}(e\vee Q_\gamma)^{-2}
        \leq e^{-A\ell}Q_\gamma^{-2}.
    \]
    Thus $\PL[-](\nu)$ and $h_\ell=\hmax e^{-\ell}$ give
    \begin{align*}
        \sum_{\gamma\in\Gn}
        (v_\gamma^\star+\varepsilon_\gamma+c_\gamma^2)
        e^{-\mathrm{pen}(\gamma)}
        &\leq\frac{C_\vartheta}{n}
        \sum_{\ell\geq1}\frac{e^{-A\ell}}{V_{h_\ell}}\\
        &\leq\frac{C_\vartheta}{n c_\nu^{\mathrm{PL}}\hmax^{d+\nu}}
        \sum_{\ell\geq1}e^{-\ell(A-d-\nu)}.
    \end{align*}
    Since $A>d+\nu$, the series converges, proving the claim with
    \[
        C_{\mathrm{sum}}=\frac{C_\vartheta}
        {c_\nu^{\mathrm{PL}}\hmax^{d+\nu}(e^{A-d-\nu}-1)}.
    \]
\end{proof}

\Restate{rethm}{thm:oracle}{\StmtOracle}

\begin{proof}
    The proof follows the general selection argument of
    \citet[Section~3]{Lepski2022adaptive}.
    We retain it for completeness because our comparison thresholds
    differ slightly from those used in his formulation.
    The argument below accommodates these thresholds without requiring
    monotonicity along the candidate family.

    Recall that $v_\gamma=\EE[\hat v_\gamma]$. We shall compare the
    random threshold $\hat{\mathbb U}_\gamma$ with the two deterministic
    thresholds
    \[
        \mathbb{U}_\gamma
        :=
        r_\gamma\bigl(v_\gamma,\mathrm{pen}(\gamma)\bigr),
        \qquad
        \bar{\mathbb{U}}_\gamma
        :=
        r_\gamma\bigl(
            v_\gamma^\star+2\varepsilon_\gamma,
            \mathrm{pen}(\gamma)
        \bigr),
        \quad \gamma\in\Gn.
    \]
    The second one is precisely the deterministic term appearing
    in~\eqref{eq:oracle.adaptive}. Since $r_\gamma(v,x)$ is
    nondecreasing in $v$ and $v_\gamma\leq v_\gamma^\star$, we have
    \[
        \mathbb U_\gamma\leq\bar{\mathbb U}_\gamma,
        \qquad \gamma\in\Gn.
    \]
    \proofpart{Pathwise oracle bound}
    Introduce the event on which every empirical variance proxy is
    sufficiently close to its expectation:
    \[
        \mathcal E
        :=
        \bigcap_{\gamma\in\Gn}
        \bigl\{|\hat v_\gamma-v_\gamma|\leq\varepsilon_\gamma\bigr\}.
    \]
    On $\mathcal E$, for every $\gamma\in\Gn$,
    \[
        v_\gamma
        \leq
        \hat v_\gamma+\varepsilon_\gamma
        \leq
        v_\gamma+2\varepsilon_\gamma
        \leq
        v_\gamma^\star+2\varepsilon_\gamma.
    \]
    The monotonicity of $r_\gamma(\cdot,\mathrm{pen}(\gamma))$ therefore
    yields
    \begin{equation}\label{eq:proof.E.sandwich}
        \mathbb{U}_\gamma\;\leq\;\hat{\mathbb{U}}_\gamma\;\leq\;\bar{\mathbb{U}}_\gamma
        \qquad\text{on }\mathcal{E},\ \text{for every }\gamma\in\Gn.
    \end{equation}

    To control all stochastic fluctuations simultaneously, define
    \[
        T
        :=
        \max_{\gamma\in\Gn}
        \bigl\{
            |\hat f_\gamma(t)-\EE\hat f_\gamma(t)|
            -\mathbb U_\gamma
        \bigr\}_+.
    \]
    By construction, for every $\gamma\in\Gn$,
    \begin{equation}\label{eq:proof.T.def}
        |\hat f_\gamma(t)-\EE\hat f_\gamma(t)|
        \leq
        \mathbb{U}_\gamma+T
        \leq
        \hat{\mathbb{U}}_\gamma+T
        \qquad\text{on }\mathcal E,
    \end{equation}
    where the second inequality follows from~\eqref{eq:proof.E.sandwich}.

    \proofstep{Goldenshluger--Lepski decomposition}
    Fix $\gamma\in\Gn$. By the triangle inequality through $\hat\gamma\vee\gamma$,
    \[
        |\hat f_{\hat\gamma}(t)-f(t)|
        \;\leq\;
        \bigl|\hat f_{\hat\gamma}(t)-\hat f_{\hat\gamma\vee\gamma}(t)\bigr|
        +\bigl|\hat f_{\hat\gamma\vee\gamma}(t)-\hat f_\gamma(t)\bigr|
        +|\hat f_\gamma(t)-f(t)|.
    \]
    For any $\gamma_1,\gamma_2\in\Gn$, the definition
    of $\hat A_{\gamma_1}$ in~\eqref{eq:Agamma.new}, evaluated at
    $\gamma'=\gamma_2$, implies
    \[
        |\hat f_{\gamma_1\vee\gamma_2}(t)-\hat f_{\gamma_2}(t)|
        \leq
        \hat A_{\gamma_1}
        +\hat{\mathbb U}_{\gamma_1\vee\gamma_2}
        +\hat{\mathbb U}_{\gamma_2}.
    \]
    Apply this observation first to
    $(\gamma_1,\gamma_2)=(\gamma,\hat\gamma)$ and then to
    $(\hat\gamma,\gamma)$. Since
    $\gamma\vee\hat\gamma=\hat\gamma\vee\gamma$, we obtain
    \begin{align*}
        \bigl|\hat f_{\hat\gamma\vee\gamma}(t)-\hat f_{\hat\gamma}(t)\bigr|
        &\leq
        \hat A_\gamma
        +\hat{\mathbb{U}}_{\hat\gamma\vee\gamma}
        +\hat{\mathbb{U}}_{\hat\gamma},\\
        \bigl|\hat f_{\hat\gamma\vee\gamma}(t)-\hat f_\gamma(t)\bigr|
        &\leq
        \hat A_{\hat\gamma}
        +\hat{\mathbb{U}}_{\hat\gamma\vee\gamma}
        +\hat{\mathbb{U}}_\gamma.
    \end{align*}
    The selection rule~\eqref{eq:gammahat.new} gives
    \[
        \hat A_{\hat\gamma}+\hat{\mathbb U}_{\hat\gamma}
        \leq
        \hat A_\gamma+\hat{\mathbb U}_\gamma.
    \]
    Combining these bounds with the triangle inequality yields
    \begin{equation}\label{eq:proof.triangle}
        |\hat f_{\hat\gamma}(t)-f(t)|
        \;\leq\;
        2\hat A_\gamma+2\hat{\mathbb{U}}_\gamma+2\hat{\mathbb{U}}_{\hat\gamma\vee\gamma}+|\hat f_\gamma(t)-f(t)|.
    \end{equation}

    No monotonicity of the penalty is needed to control the middle
    threshold. Indeed, the definition of $\vee$ ensures that
    $\hat\gamma\vee\gamma$ is equal to either $\hat\gamma$ or $\gamma$.
    Consequently,
    \[
        \hat{\mathbb{U}}_{\hat\gamma\vee\gamma}\;\leq\;\hat{\mathbb{U}}_{\hat\gamma}+\hat{\mathbb{U}}_\gamma,
    \]
    while the selection rule and $\hat A_{\hat\gamma}\geq0$ imply
    \[
        \hat{\mathbb U}_{\hat\gamma}
        \leq
        \hat A_\gamma+\hat{\mathbb U}_\gamma.
    \]
    Substitution into~\eqref{eq:proof.triangle} gives the pathwise
    decomposition
    \begin{equation}\label{eq:proof.GL.decomp}
        |\hat f_{\hat\gamma}(t)-f(t)|
        \;\leq\;
        4\hat A_\gamma+6\hat{\mathbb{U}}_\gamma+|\hat f_\gamma(t)-f(t)|.
    \end{equation}

    \proofstep{Control of the empirical bias proxy}
    Consider first $\gamma'\in\Gn$ such that
    $\gamma\preceq\gamma'$ and $\gamma'\neq\gamma$. Then
    $\gamma\vee\gamma'=\gamma'$, so the corresponding term in the
    maximum defining $\hat A_\gamma$ is
    $\{0-2\hat{\mathbb U}_{\gamma'}\}_+=0$.
    Since the grid is totally ordered by $\preceq$, only the candidates
    $\gamma'\preceq\gamma$ can therefore contribute.
    For them, $\gamma\vee\gamma'=\gamma$, and~\eqref{eq:Agamma.new}
    reduces to
    \[
        \hat A_\gamma
        =\max_{\substack{\gamma'\in\Gn\\\gamma'\preceq\gamma}}
        \bigl\{|\hat f_\gamma(t)-\hat f_{\gamma'}(t)|-\hat{\mathbb{U}}_\gamma-\hat{\mathbb{U}}_{\gamma'}\bigr\}_+.
    \]

    Fix such a $\gamma'$. Using
    $\EE\hat f_\gamma(t)=f(t)+b_\gamma$ and
    $\EE\hat f_{\gamma'}(t)=f(t)+b_{\gamma'}$, the triangle inequality
    gives
    \begin{align*}
        |\hat f_\gamma(t)-\hat f_{\gamma'}(t)|
        &\leq
        |\hat f_\gamma(t)-\EE\hat f_\gamma(t)|
        +|\hat f_{\gamma'}(t)-\EE\hat f_{\gamma'}(t)|\\
        &\quad
        +|b_\gamma|+|b_{\gamma'}|.
    \end{align*}
    Since $\gamma'\preceq\gamma$, the definition of $B(\gamma)$ gives
    $|b_\gamma|\leq B(\gamma)$ and
    $|b_{\gamma'}|\leq B(\gamma)$.
    Applying~\eqref{eq:proof.T.def} to both stochastic terms, we obtain,
    on $\mathcal E$,
    \[
        |\hat f_\gamma(t)-\hat f_{\gamma'}(t)|
        \;\leq\;
        \hat{\mathbb{U}}_\gamma+\hat{\mathbb{U}}_{\gamma'}+2T+2B(\gamma),
    \]
    uniformly over $\gamma'\preceq\gamma$. Subtracting the two
    thresholds and taking the positive part therefore yields
    \begin{equation}\label{eq:proof.Ahat.bound}
        \hat A_\gamma\;\leq\;2T+2B(\gamma) \qquad\text{on }\mathcal{E}.
    \end{equation}

    \proofstep{Control of the reference estimator}
    The last term in~\eqref{eq:proof.GL.decomp} is controlled in the
    same way. By~\eqref{eq:proof.T.def} and
    $|b_\gamma|\leq B(\gamma)$,
    \begin{equation}\label{eq:proof.last.term}
        |\hat f_\gamma(t)-f(t)|
        \;\leq\;
        |\hat f_\gamma(t)-\EE\hat f_\gamma(t)|+|b_\gamma|
        \;\leq\;
        \hat{\mathbb{U}}_\gamma+T+B(\gamma)
        \qquad\text{on }\mathcal{E}.
    \end{equation}

    \proofstep{Conclusion of the pathwise argument}
    Substitute~\eqref{eq:proof.Ahat.bound}
    and~\eqref{eq:proof.last.term} into~\eqref{eq:proof.GL.decomp}, and
    use $\hat{\mathbb U}_\gamma\leq\bar{\mathbb U}_\gamma$
    from~\eqref{eq:proof.E.sandwich}. For every $\gamma\in\Gn$, on
    $\mathcal E$,
    \begin{align*}
        |\hat f(t)-f(t)|
        &=
        |\hat f_{\hat\gamma}(t)-f(t)|\\
        &\leq
        4\{2T+2B(\gamma)\}
        +6\bar{\mathbb U}_\gamma
        +\{\bar{\mathbb U}_\gamma+T+B(\gamma)\}\\
        &=
        9T+9B(\gamma)+7\bar{\mathbb U}_\gamma.
    \end{align*}
    Since the coefficient of each term is at most $9$, taking the
    minimum over $\gamma\in\Gn$ gives
    \begin{equation}\label{eq:proof.pathwise.final}
        |\hat f(t)-f(t)|\,\1_{\mathcal{E}}
        \;\leq\;
        9\Bigl(\min_{\gamma\in\Gn}\bigl\{B(\gamma)+\bar{\mathbb{U}}_\gamma\bigr\}+T\Bigr).
    \end{equation}

    \proofpart{Expectation bounds}
    \proofstep{Contribution of the exceptional event}
    Using~\eqref{eq:cgamma.def}, we deduce that
   \revision{ \[
        c_\gamma
        =
        \frac{\sqrt{\kg\Qg}}{nV_h}
        \leq
        \frac{Q_\gamma}{nV_h}.
    \]
}
    The definition of $\Gn$ gives $Q_\gamma\leq n$, while the
    restriction $\gamma\in\{\gamma_\ell:\ell\in\mathcal L_n\}$ gives
    $nV_h\geq(\log n)^3$. Hence, for $n\geq3$,
    \[
        c_\gamma
        \leq
        \frac{n}{(\log n)^3}
        \leq
        n.
    \]

    In the notation of Lemma~\ref{lem:bernstein.fhat},
    $|\hat f_\gamma(t)|\leq n\max_i|Y_i|\leq nc_\gamma$.
    Thus
    \[
        |\hat f_\gamma(t)|\leq n^2,
        \qquad \gamma\in\Gn.
    \]
    Since $n\geq3$,
    \[
        |\hat f(t)-f(t)|
        \leq
        n^2+\finfty
        \leq
        \left(1+\frac{\finfty}{9}\right)n^2.
    \]
    Therefore,
    \[
        \EE\bigl[|\hat f(t)-f(t)|^2\,\1_{\mathcal{E}^c}\bigr]
        \leq
        \left(1+\frac{\finfty}{9}\right)^2
        n^4\PP(\mathcal E^c).
    \]
    By Lemma~\ref{lem:var.concentration} and a union bound over the at
    most $\lfloor\log n\rfloor$ candidates in $\Gn$,
    \[
        \PP(\mathcal E^c)
        \leq
        2\lfloor\log n\rfloor
        e^{-\kappa_\vartheta(\log n)^3}.
    \]
    The sequence
    \[
        n^5\lfloor\log n\rfloor
        e^{-\kappa_\vartheta(\log n)^3},
        \qquad n\geq3,
    \]
    is bounded. Consequently, there exists $C_{\mathcal E}>0$,
    depending only on $\vartheta$ and $\finfty$, such that, for every
    $n\geq3$,
    \begin{equation}\label{eq:proof.Ec}
        \EE\bigl[
            |\hat f(t)-f(t)|^2\,\1_{\mathcal E^c}
        \bigr]
        \leq
        \frac{C_{\mathcal E}}{n}.
    \end{equation}

    \proofstep{Control of the stochastic remainder}
    Fix $\gamma\in\Gn$ and define the individual excess
    \[
        T_\gamma
        :=
        \bigl\{
            |\hat f_\gamma(t)-\EE\hat f_\gamma(t)|
            -\mathbb U_\gamma
        \bigr\}_+.
    \]
    For $x\geq0$, set
    \[
        u(x)
        :=
        r_\gamma(v_\gamma,x)
        =
        \sqrt{2v_\gamma x}+c_\gamma x.
    \]
    The inequality $\sqrt a+\sqrt b\geq\sqrt{a+b}$, applied to the
    square-root terms in $r_\gamma$, gives
    \begin{align*}
        \mathbb U_\gamma+u(x)
        &=
        r_\gamma\bigl(v_\gamma,\mathrm{pen}(\gamma)\bigr)
        +r_\gamma(v_\gamma,x)\\
        &\geq
        r_\gamma\bigl(v_\gamma,\mathrm{pen}(\gamma)+x\bigr).
    \end{align*}
    Lemma~\ref{lem:bernstein.fhat} therefore yields the tail bound
    \begin{align*}
        \PP\bigl(T_\gamma>u(x)\bigr)
        &=
        \PP\left(
            |\hat f_\gamma(t)-\EE\hat f_\gamma(t)|
            >
            \mathbb U_\gamma+u(x)
        \right)\\
        &\leq
        \PP\left(
            |\hat f_\gamma(t)-\EE\hat f_\gamma(t)|
            >
            r_\gamma\bigl(v_\gamma,\mathrm{pen}(\gamma)+x\bigr)
        \right)\\
        &\leq
        2e^{-\mathrm{pen}(\gamma)}e^{-x}.
    \end{align*}

    We now integrate this tail estimate. The change of variables
    $z=u(x)$ in the usual identity
    $\EE[T_\gamma^2]=\int_0^\infty2z\,\PP(T_\gamma>z)\,\rd z$
    gives
    \[
        \EE[T_\gamma^2]
        =
        \int_0^\infty
        2u(x)\PP\bigl(T_\gamma>u(x)\bigr)u'(x)\,\rd x.
    \]
    Using the preceding tail bound and integrating by parts, we obtain
    \[
        \EE[T_\gamma^2]
        \leq
        2e^{-\mathrm{pen}(\gamma)}
        \int_0^\infty u(x)^2e^{-x}\,\rd x.
    \]
    Here the boundary terms vanish because $u(0)=0$ and
    $u(x)^2e^{-x}\to0$ as $x\to\infty$.

    It remains to evaluate the last integral. Since
    $(a+b)^2\leq2a^2+2b^2$,
    \[
        u(x)^2
        \leq
        4v_\gamma x+2c_\gamma^2x^2.
    \]
    Together with
    \[
        \int_0^\infty xe^{-x}\,\rd x=1,
        \qquad
        \int_0^\infty x^2e^{-x}\,\rd x=2,
    \]
    this gives
    \[
        \EE[T_\gamma^2]
        \leq
        8e^{-\mathrm{pen}(\gamma)}
        \bigl(v_\gamma+c_\gamma^2\bigr).
    \]
    Since $v_\gamma\leq v_\gamma^\star$, we may add the non-negative
    term $\varepsilon_\gamma$ in order to apply
    Lemma~\ref{lem:sums.adaptive} directly:
    \[
        \EE[T_\gamma^2]
        \leq
        8
        \bigl(v_\gamma^\star+\varepsilon_\gamma+c_\gamma^2\bigr)
        e^{-\mathrm{pen}(\gamma)}.
    \]

    Finally, $T=\max_{\gamma\in\Gn}T_\gamma$, and hence
    \[
        T^2
        \leq
        \sum_{\gamma\in\Gn}T_\gamma^2.
    \]
    Taking expectations and applying Lemma~\ref{lem:sums.adaptive}
    yields
    \begin{equation}\label{eq:proof.ET2}
        \EE[T^2]
        \leq
        8\sum_{\gamma\in\Gn}
        \bigl(v_\gamma^\star+\varepsilon_\gamma+c_\gamma^2\bigr)
        e^{-\mathrm{pen}(\gamma)}
        \leq
        \frac{8C_{\mathrm{sum}}}{n}.
    \end{equation}

    \proofpart{Conclusion of the oracle inequality}
    To simplify notation, set
    \[
        \mathcal O_n
        :=
        \min_{\gamma\in\Gn}
        \bigl\{B(\gamma)+\bar{\mathbb U}_\gamma\bigr\}.
    \]
    The pathwise bound~\eqref{eq:proof.pathwise.final} and the inequality
    $(a+b)^2\leq2a^2+2b^2$ imply
    \[
        \EE\bigl[|\hat f(t)-f(t)|^2\,\1_{\mathcal{E}}\bigr]
        \leq
        162\mathcal O_n^2
        +162\EE[T^2].
    \]
    Using~\eqref{eq:proof.ET2}, we obtain
    \begin{equation}\label{eq:proof.E}
        \EE\bigl[|\hat f(t)-f(t)|^2\,\1_{\mathcal{E}}\bigr]
        \leq
        162\mathcal O_n^2
        +
        \frac{1296C_{\mathrm{sum}}}{n}.
    \end{equation}

    We now add the contribution of the exceptional event. Combining
    \eqref{eq:proof.E} with~\eqref{eq:proof.Ec} gives
    \[
        \pwrisk(\hat f,f)^2
        =
        \EE\bigl[|\hat f(t)-f(t)|^2\bigr]
        \leq
        162\mathcal O_n^2
        +
        \frac{1296C_{\mathrm{sum}}+C_{\mathcal E}}{n}.
    \]
    Taking square roots and using
    $\sqrt{a+b}\leq\sqrt a+\sqrt b$, we conclude that
    \[
        \pwrisk(\hat f,f)
        \leq
        9\sqrt{2}\,\mathcal O_n
        +\frac{\sqrt{1296C_{\mathrm{sum}}+C_{\mathcal E}}}{\sqrt n}.
    \]
    This is~\eqref{eq:oracle.adaptive} with
    \[
        C_1:=9\sqrt2,
        \qquad
        C_2:=\sqrt{1296C_{\mathrm{sum}}+C_{\mathcal E}}.
    \]
    This proves~\eqref{eq:oracle.adaptive} for every $n\geq3$ for which
    $\Gn$ is nonempty.
\end{proof}

\subsection{Proof of Theorem~\ref{thm:adaptive}}
\label{sec:proof.adaptive.rate}

The adaptive result follows by evaluating the oracle bound at a
suitably chosen degree--bandwidth pair.

\Restate{rethm}{thm:adaptive}{\StmtAdaptive}

\begin{proof}
    Fix $s,L>0$. For $f\in\Sigma(s,L)$, let $q$ be its approximating
    polynomial, and take the norms below on $\Ehmax$ with respect to the
    local probability measure $\mu_{\hmax}$. The Bernstein--Markov bound
    at degree $\llfloor s\rrfloor$ gives
    \[
        \|q\|_\infty^2
        \leq C_{\llfloor s\rrfloor}^{\mathrm{BM}}\|q\|_2^2
        \leq C_{\llfloor s\rrfloor}^{\mathrm{BM}}
             \|q\|_\infty\|q\|_1.
    \]
    Nonnegativity and unit mass of $f$ give
    $\|q\|_1\leq V_{\hmax}^{-1}+L(\hmax)^s$.
    By $\PL[-](\nu)$ this controls
    $\finfty\leq\|q\|_\infty+L(\hmax)^s$ uniformly by a constant
    $\FsL$ with the dependencies stated in the theorem.
    For all sufficiently large $n$, set
    \[
        \ell^\star=
        \left\lfloor\frac{\log(n/\log n)}{2s+d+\nu}\right\rfloor,
        \qquad \gamma^\star=\gamma_{\ell^\star}.
    \]

    \proofpart{The oracle candidate belongs to the grid}
    The definitions give
    \begin{equation}\label{eq:proof.hstar}
        h_{\ell^\star}\asymp
        \left(\frac{\log n}{n}\right)^{1/(2s+d+\nu)}\longrightarrow0.
    \end{equation}
    Also $1\leq\ell^\star\leq\lfloor\log n\rfloor$ for large $n$, and
    $\PL[-](\nu)$ yields
    \[
        nV_{h_{\ell^\star}}
        \gtrsim n^{2s/(2s+d+\nu)}
        (\log n)^{(d+\nu)/(2s+d+\nu)}.
    \]
    This exceeds $(\log n)^3$ eventually, so $\ell^\star\in\mathcal L_n$.

    Since
    \[
        \frac{\log n}{2\ell^\star}\longrightarrow
        s+\frac{d+\nu}{2}>s,
    \]
    the degree $m_{\ell^\star}=\lfloor\log n/(2\ell^\star)\rfloor$
    satisfies, for large $n$,
    \begin{equation}\label{eq:proof.mstar.range}
        \llfloor s\rrfloor\leq m_{\ell^\star}\leq\bar m,
        \qquad
        \bar m:=\left\lceil\frac{2s+d+\nu}{2}\right\rceil+1.
    \end{equation}
    Uniform regularity and Proposition~\ref{prop:BM} imply
    \[
        Q_{\gamma^\star}\leq
        \max_{\llfloor s\rrfloor\leq m\leq\bar m}C_m^{\mathrm{BM}}
        =:C_{\mathrm{BM}}^\star<\infty.
    \]
    Hence $Q_{\gamma^\star}\leq n$ eventually, proving
    $\gamma^\star\in\Gn$.

    \proofpart{Stochastic threshold}
    At these finitely many degrees, uniform regularity gives
    \[
        1\leq k_{\gamma^\star}(0,0)\leq Q_{\gamma^\star}\leq C_{\mathrm{BM}}^\star.
    \]
    By the definitions of the thresholds and $\finfty\leq\FsL$,
    \[
        v_{\gamma^\star}^\star+2\varepsilon_{\gamma^\star}
        \lesssim\frac{1}{nh_{\ell^\star}^{d+\nu}},
        \qquad
        c_{\gamma^\star}\lesssim\frac{1}{nh_{\ell^\star}^{d+\nu}}.
    \]
    The boundedness of $Q_{\gamma^\star}$ and $m_{\ell^\star}$, together
    with $\ell^\star\asymp\log n$, yields
    $\mathrm{pen}(\gamma^\star)\lesssim\log n$.
    Since~\eqref{eq:proof.hstar} implies
    $h_{\ell^\star}^{2s+d+\nu}\asymp\log n/n$, the deterministic
    threshold from the oracle inequality satisfies
    \begin{align*}
        \bar{\mathbb U}_{\gamma^\star}
        &=r_{\gamma^\star}
        (v_{\gamma^\star}^\star+2\varepsilon_{\gamma^\star},
        \mathrm{pen}(\gamma^\star))\\
        &\lesssim
        \sqrt{\frac{\log n}{nh_{\ell^\star}^{d+\nu}}}
        +\frac{\log n}{nh_{\ell^\star}^{d+\nu}}
        \lesssim h_{\ell^\star}^s.
    \end{align*}

    \proofpart{Bias over finer candidates}
    By definition,
    \[
        B(\gamma^\star)
        =\max_{\substack{\ell'\geq\ell^\star\\\gamma_{\ell'}\in\Gn}}
        |b_{\gamma_{\ell'}}|.
    \]
    For these candidates, $h_{\ell'}\leq h_{\ell^\star}$ and
    $0\leq m_{\ell'}\leq m_{\ell^\star}\leq\bar m$.
    The general bias bound and \PBM at these finitely many degrees give
    \begin{equation}\label{eq:proof.bias.common}
        |b_{\gamma_{\ell'}}|
        \leq C_{\mathrm{bias}}h_{\ell'}^{\min(m_{\ell'}+1,s)},
    \end{equation}
    uniformly over the class and the retained candidates.

    If $m_{\ell'}\geq\llfloor s\rrfloor$, the right-hand side is at
    most $C_{\mathrm{bias}}h_{\ell^\star}^s$.
    Otherwise, the degree definition gives
    $\ell'(m_{\ell'}+1)>\tfrac12\log n$, and hence
    \[
        h_{\ell'}^{m_{\ell'}+1}
        =\hmax^{m_{\ell'}+1}e^{-\ell'(m_{\ell'}+1)}
        \leq\max(1,\hmax)^{\bar m+1}n^{-1/2}.
    \]
    Since $s/(2s+d+\nu)<1/2$, we have
    $n^{-1/2}=o(h_{\ell^\star}^s)$.
    Both cases therefore yield
    \begin{equation}\label{eq:proof.Bgammastar}
        B(\gamma^\star)\lesssim h_{\ell^\star}^s.
    \end{equation}

    \proofpart{Conclusion}
    Applying Theorem~\ref{thm:oracle} at $\gamma^\star$ gives
    \[
        \pwrisk(\hat f,f)
        \lesssim B(\gamma^\star)+\bar{\mathbb U}_{\gamma^\star}+n^{-1/2}
        \lesssim\left(\frac{\log n}{n}\right)^{s/(2s+d+\nu)},
    \]
    where the bounds on the bias and threshold are uniform in $f$.
    For the finitely many remaining $n\geq3$, the estimator is zero
    when $\Gn$ is empty. Otherwise, the proof of the oracle inequality
    gives $|\hat f(t)|\leq n^2$.
    Thus $\pwrisk(\hat f,f)\leq n^2+\FsL$ in either case.
    Enlarging the constant proves~\eqref{eq:rate.adaptive} for all
    $n\geq3$, with the dependencies stated in the theorem.
\end{proof}

\section{Proofs of Section \emph{The exponential cusp}}\label{sec:supp.Dexp}

The calculations below use scalar bandwidths.
Section~\ref{sec:proof.Dexp.anisotropic} uses the vector $\bh(r)$.

\subsection{Proof of Theorem~\ref{thm:Dexp.best.fixed.degree}}
\label{sec:proof.single.Dexp}

The fixed-degree analysis rests on two properties of the equivalent
kernel.
We first identify its behaviour in the natural boundary-layer
coordinates of the exponential cusp.
We then quantify the resulting failure to reproduce monomials above the
fitted degree.
Combining these two ingredients with the variance gives the best risk
attainable at any fixed degree.

\begin{lemma}[$L^1$ blow-up of the equivalent kernel]\label{lem:Dexp.kernel.L1}
    Fix $m\in\N_0$ and set $\gamma=(m,h)$. There exist constants
    $\Xi_m,\Lambda_m^{\mathrm{exp}}\in(0,\infty)$ such that, as $h\to0$,
    \begin{align}
        h^{2m} k_\gamma(0,0)&\longrightarrow\Xi_m,
        \label{eq:Dexp.kernel.diagonal}\\
        h^m\int_{\Eh}| k_\gamma(0,u)|\,\rd\tmu(u)
        &\longrightarrow\Lambda_m^{\mathrm{exp}}.
        \label{eq:Dexp.kernel.L1}
    \end{align}
    Moreover, for every $s>0$,
    \begin{equation}\label{eq:Dexp.kernel.weighted.L1}
        h^{m-s}\int_{\Eh}| k_\gamma(0,u)|\,
        \|u\|_\infty^s\,\rd\tmu(u)\longrightarrow\Lambda_m^{\mathrm{exp}}.
    \end{equation}
\end{lemma}

\begin{proof}
    \proofpart{Natural boundary-layer coordinates}
    Let $(X_h,Y_h)$ be a vector with distribution $\tmu$, which is the uniform
    probability measure on $\Eh$.
    Define
    \begin{equation*}
        Z_h:=\frac{h-X_h}{h^2},
        \qquad
        T_h:=Y_h e^{1/X_h}.
    \end{equation*}
    The inverse coordinate map is
    \begin{equation}\label{eq:Dexp.coordinate.map}
        F_h(z,t)
        :=
        \left(h-h^2z,\,
        e^{-1/(h-h^2z)}t\right),
        \qquad
        (z,t)\in\Omega_h:=(0,1/h)\times(0,1).
    \end{equation}
    Thus $(X_h,Y_h)=F_h(Z_h,T_h)$ almost surely.
    
    \revision{We have for $(z,t)\in\Omega_h$
\begin{align*}
        \PP(Z_h\leq z,T_h\leq t)
        &=\frac{1}{V_h}\int_{h-zh^2}^h\int_0^{te^{-1/x}}\rd y\,\rd x\\
        &=\frac{t}{V_h}\int_{h-zh^2}^h e^{-1/x}\,\rd x
        =t\int_0^z q_h(a)\,\rd a.
    \end{align*}
Thus $T_h$ is uniform on $(0,1)$ and independent of $Z_h$. Moreover,
   their joint distribution with support on $\Omega_h$ has density $f_{Z_h,T_h}(z,t)=q_h(z)$ where
    \begin{equation}\label{eq:Dexp.zh.density}
        q_h(z)
        =\frac{h^2}{V_h}\exp\!\left\{-\frac{1}{h-h^2z}\right\}=\frac{h^2e^{-1/h}}{V_h}\exp\!\left\{-\frac{z}{1-hz}\right\},
    \end{equation}
   where we use that
    \begin{equation}\label{eq:Dexp.reciprocal.identity}
        \frac1{h-h^2z}=\frac1h+\frac{z}{1-hz}.
    \end{equation}
   Using that $\lim_{h\to0}\frac{h^2e^{-1/h}}{V_h}=1$ and that $z/(1-hz)\ge z$ on $(0,1/h)$, we obtain that
    \begin{equation}\label{eq:Dexp.zh.limit}
        \lim_{h\to 0}q_h(z)= e^{-z},
        \qquad
        q_h(z)\le C_q e^{-z}
    \end{equation}
    for all sufficiently small $h$, after extending $q_h$ by zero to
    $(0,\infty)$. } Thus $(Z_h,T_h)$ converges in distribution to an
    independent pair $(Z,T)$ with $Z\sim\operatorname{Exp}(1)$ and
    $T\sim\mathcal U[0,1]$. Relation~\eqref{eq:Dexp.zh.limit} also gives
    convergence of all polynomial moments of $Z_h$.

    \proofpart{Adapted basis and Gram limit}
    Let
    \begin{equation*}
        \mathcal I_m:=\{(r,j)\in\N_0^2:r+j\le m\}
    \end{equation*}
    and, for $(r,j)\in\mathcal I_m$, define the polynomial
    \begin{equation}\label{eq:Dexp.adapted.basis}
        p_{r,j}^{(h)}(x,y)
        :=h^{-2r}e^{j/h}(h-x)^r y^j.
    \end{equation}
    Translation in the first coordinate and multiplication by non-zero
    constants show that
    $\{p_{r,j}^{(h)}:(r,j)\in\mathcal I_m\}$ is a basis of $\Pm$.
    Write
    \[
        \mathbf p_h(x,y)
        :=\bigl(p_{r,j}^{(h)}(x,y)\bigr)_{(r,j)\in\mathcal I_m}.
    \]
    Pulling this basis back through $F_h$, define
    \begin{equation}\label{eq:Dexp.basis.transform}
        \phi_{r,j}^{(h)}(z,t)
        :=p_{r,j}^{(h)}\bigl(F_h(z,t)\bigr)
        =z^r\exp\!\left\{-\frac{jz}{1-hz}\right\}t^j
    \end{equation}
   where we use \eqref{eq:Dexp.reciprocal.identity}.
   Consider the polynomial basis
    \[
        \boldsymbol\phi_h(z,t)
        :=\bigl(\phi_{r,j}^{(h)}(z,t)\bigr)_{(r,j)\in\mathcal I_m}
    \]
   and the associated Gram matrix 
    \begin{equation}\label{eq:Dexp.Gram.pullback}
        \Gamma_h
        :=\EE\!\left[
        \mathbf p_h(X_h,Y_h)\mathbf p_h(X_h,Y_h)^\top\right]
        =
        \int_{\Omega_h}
        \boldsymbol\phi_h(z,t)\boldsymbol\phi_h(z,t)^\top
        q_h(z)\,\rd z\,\rd t.
    \end{equation}
    For every fixed $(z,t)\in(0,\infty)\times(0,1)$,
    $\boldsymbol\phi_h(z,t)$ converges componentwise to
    \begin{equation}\label{eq:Dexp.limit.family}
        \boldsymbol\phi_\star(z,t)
        :=\bigl(z^r e^{-jz}t^j\bigr)_{(r,j)\in\mathcal I_m}.
    \end{equation}
    Moreover, every component of $\boldsymbol\phi_h(z,t)$ is bounded by
    $z^r$ on $\Omega_h$.
    Relations~\eqref{eq:Dexp.zh.limit} and~\eqref{eq:Dexp.Gram.pullback}
    thus show, by dominated convergence, that $\Gamma_h$ converges
    entrywise to the Gram matrix $\Gamma_\star$ of
    $\boldsymbol\phi_\star$ in
    $L^2(e^{-z}\rd z\,\rd t)$ on $(0,\infty)\times(0,1)$.
    This limiting family is linearly independent. Indeed, if a linear
    combination vanishes almost everywhere, then, for almost every $z$,
    it is the zero polynomial in $t$. For each $j$, the resulting
    polynomial in $z$, multiplied by $e^{-jz}$, must therefore vanish
    identically.
    Hence $\Gamma_\star$ is positive definite, and as $h\to 0$
    \begin{equation}\label{eq:Dexp.Gram.limit}
        \Gamma_h^{-1}\longrightarrow\Gamma_\star^{-1}.
    \end{equation}

    \proofpart{Kernel limits}
    Recall that $k_{m,h}$ is the reproducing kernel of $\Pm$ in $L^2(\tmu)$.
    The coordinate map $F_h$ parametrises only the interior of $\Eh$ and
    does not provide coordinates for the target $(0,0)$.
    We therefore compute its evaluation vector directly from the
    polynomial basis.
    Set $\mathbf v_h:=\mathbf p_h(0,0)$.
    By~\eqref{eq:Dexp.adapted.basis},
    \begin{equation*}
        (\mathbf v_h)_{r,j}=h^{-r}\1_{\{j=0\}},
        \qquad
        \lim_{h\to 0 }h^m\mathbf v_h= e_{(m,0)},
    \end{equation*}
    where $e_{(m,0)}$ is the coordinate vector indexed by $(m,0)$.
    The standard kernel formula in the basis $\mathbf p_h$ gives
    \begin{equation}\label{eq:Dexp.kernel.original.coordinates}
        k_{m,h}(0,(x,y))
        =\mathbf v_h^\top\Gamma_h^{-1}\mathbf p_h(x,y).
    \end{equation}
    Pulling the second argument back through $F_h$ yields
    \begin{equation}\label{eq:Dexp.kernel.representation}
        k_{m,h}\bigl(0,F_h(z,t)\bigr)
        =\mathbf v_h^\top\Gamma_h^{-1}
        \boldsymbol\phi_h(z,t).
    \end{equation}
    Evaluating~\eqref{eq:Dexp.kernel.original.coordinates} at the origin
    gives
    \begin{equation*}
        h^{2m}k_{m,h}(0,0)
        =(h^m\mathbf v_h)^\top
        \Gamma_h^{-1}(h^m\mathbf v_h)
        \longrightarrow
        e_{(m,0)}^\top\Gamma_\star^{-1}e_{(m,0)}
        =:\Xi_m>0,
    \end{equation*}
    which proves~\eqref{eq:Dexp.kernel.diagonal}.

    For each fixed
    $(z,t)\in(0,\infty)\times(0,1)$,
    \begin{equation*}
        h^mk_{m,h}\bigl(0,F_h(z,t)\bigr)
        =
        (h^m\mathbf v_h)^\top\Gamma_h^{-1}
        \boldsymbol\phi_h(z,t)
        \longrightarrow
        k_{m,\star}(z,t)
        :=e_{(m,0)}^\top\Gamma_\star^{-1}
        \boldsymbol\phi_\star(z,t).
    \end{equation*}
    The vector $h^m\mathbf v_h$ and the matrices $\Gamma_h^{-1}$ are
    bounded for small $h$.
    Together with~\eqref{eq:Dexp.zh.limit}, this bounds the absolute value
    of the integrand by $C_m(1+z^m)e^{-z}$.
    Dominated convergence therefore gives
    \begin{align*}
        h^m\int_{\Eh}
        |k_{m,h}(0,u)|\,\rd\tmu(u)
        &=
        \int_{\Omega_h}
        \left|(h^m\mathbf v_h)^\top\Gamma_h^{-1}
        \boldsymbol\phi_h(z,t)\right|
        q_h(z)\,\rd z\,\rd t\\
        &\longrightarrow
        \int_0^\infty\int_0^1
        |k_{m,\star}(z,t)|e^{-z}\,\rd t\,\rd z
        =:\Lambda_m^{\mathrm{exp}}<\infty.
    \end{align*}
    Finally, $k_{m,\star}$ is not zero in
    $L^2(e^{-z}\rd z\,\rd t)$: otherwise the linear independence of
    $\boldsymbol\phi_\star$ would imply
    $\Gamma_\star^{-1}e_{(m,0)}=0$.
    Hence
    $\Lambda_m^{\mathrm{exp}}>0$, proving~\eqref{eq:Dexp.kernel.L1}.

    To obtain the weighted limit, note that $e^{-1/x}\le x$ on $(0,1)$, so
    $\|(x,y)\|_\infty=x$ on $\Dexp$. In the coordinates above,
    $X_h/h=1-hZ_h$. Therefore
    \begin{align*}
        &h^{m-s}\int_{\Eh}
        |k_{m,h}(0,u)|\,\|u\|_\infty^s\,\rd\tmu(u)\\
        &\qquad =
        \int_{\Omega_h}
        \left|(h^m\mathbf v_h)^\top\Gamma_h^{-1}
        \boldsymbol\phi_h(z,t)\right|
        (1-hz)^s q_h(z)\,\rd z\,\rd t.
    \end{align*}
    The additional factor $(1-hz)^s$ is bounded by one and converges
    pointwise to one. The same dominated-convergence argument thus yields
    the limit $\Lambda_m^{\mathrm{exp}}$, proving~\eqref{eq:Dexp.kernel.weighted.L1} and
    completing the proof.
\end{proof}

The kernel limit also quantifies the residual moments that remain after
projection onto $\Pm$.

\begin{lemma}[Polynomial leakage]\label{lem:Dexp.polynomial.leakage}
    Fix $m\in\N_0$.
    As $h\to0$, for every fixed integer $k\geq m+1$,
    \begin{equation}\label{eq:Dexp.xk.leakage}
        \int_{\Eh} k_{m,h}(0,u)x^k\,\rd\tmu(u)
        =
        (-1)^m\binom{k-1}{m}h^k+O(h^{k+1}).
    \end{equation}
    In particular,
    \begin{equation}\label{eq:Dexp.first.leakage}
        \int_{\Eh} k_{m,h}(0,u)x^{m+1}\,\rd\tmu(u)
        =
        (-1)^m h^{m+1}\{1+O(h)\}.
    \end{equation}
    Moreover, for every fixed $i,j\in\N_0$ with $j\geq1$,
    \begin{equation}\label{eq:Dexp.y.leakage}
        \int_{\Eh} k_{m,h}(0,u)x^iy^j\,\rd\tmu(u)
        =
        O\!\left(h^{i-m}e^{-j/h}\right).
    \end{equation}
\end{lemma}

\begin{proof}
    \proofpart{Kernel envelope}
    By~\eqref{eq:Dexp.kernel.representation},
    $k_{m,h}(0,F_h(z,t))$ is the scalar product of
    $\mathbf v_h^\top\Gamma_h^{-1}$ and
    $\boldsymbol\phi_h(z,t)$.
    The vector $\mathbf v_h$ has norm at most $C_mh^{-m}$,
    the matrices $\Gamma_h^{-1}$ are bounded by
    \eqref{eq:Dexp.Gram.limit}, and every component of
    $\boldsymbol\phi_h(z,t)$ is bounded by $z^r$.
    Consequently,
    \begin{equation}\label{eq:Dexp.kernel.pointwise.bound}
        \left|k_{m,h}\bigl(0,F_h(z,t)\bigr)\right|
        \le C_m h^{-m}(1+z^m).
    \end{equation}
    Relation~\eqref{eq:Dexp.zh.limit} shows that all polynomial moments
    under $q_h(z)\,\rd z$ are bounded uniformly for small $h$.

    \proofpart{Pure \(x\)-monomials}
    For $r\le m$, polynomial reproduction gives
    \begin{equation*}
        \int_{\Eh}
        k_{m,h}(0,u)(h-x)^r\,\rd\tmu(u)=h^r.
    \end{equation*}
    For $r>m$, the change of variables~\eqref{eq:Dexp.coordinate.map}
    and~\eqref{eq:Dexp.kernel.pointwise.bound} give instead
    \begin{align*}
        \left|\int_{\Eh}
        k_{m,h}(0,u)(h-x)^r\,\rd\tmu(u)\right|
        &\le
        C_mh^{2r-m}
        \int_{\Omega_h}(1+z^m)z^r q_h(z)\,\rd z\,\rd t\\
        \le C_{m,r}h^{2r-m}.
    \end{align*}
    Expanding $x^k=\{h-(h-x)\}^k$, we obtain, for $k\ge m+1$,
    \begin{align*}
        \int_{\Eh}k_{m,h}(0,u)x^k\,\rd\tmu(u)
        &=
        h^k\sum_{r=0}^m(-1)^r\binom{k}{r}+
        \sum_{r=m+1}^k(-1)^r\binom{k}{r}
        h^{k-r}O(h^{2r-m})\\
        &=h^k\sum_{r=0}^m(-1)^r\binom{k}{r}
        +O(h^{k+1})\\
        &=(-1)^m\binom{k-1}{m}h^k+O(h^{k+1}),
    \end{align*}
    proving~\eqref{eq:Dexp.xk.leakage} and
    \eqref{eq:Dexp.first.leakage}.

    \proofpart{Monomials involving \(y\)}
    Finally, if $(x,y)=F_h(z,t)$, then $x\le h$ and
    \[
        y=e^{-1/h}
        \exp\!\left\{-\frac{z}{1-hz}\right\}t.
    \]
    Hence~\eqref{eq:Dexp.kernel.pointwise.bound} gives
    \begin{align*}
        &\left|\int_{\Eh}
        k_{m,h}(0,u)x^iy^j\,\rd\tmu(u)\right|\\
        &\qquad\le
        C_mh^{i-m}e^{-j/h}
        \int_{\Omega_h}
        (1+z^m)
        \exp\!\left\{-\frac{jz}{1-hz}\right\}
        t^j q_h(z)\,\rd z\,\rd t\\
        &\qquad\le C_{m,i,j}h^{i-m}e^{-j/h},
    \end{align*}
    which proves~\eqref{eq:Dexp.y.leakage}.
\end{proof}

These two kernel estimates determine both sources of fixed-degree bias
and the corresponding variance cost.

\begin{lemma}[Risk at a fixed degree]\label{lem:Dexp.fixed.degree}
    Let $s,L>0$ and fix $m\in\N_0$.
    If $m\leq\llfloor s\rrfloor$, set
    \[
        \alpha_m(s):=\min\{m+1,s-m\},
    \]
    and denote by $h_{n,m}$ the unique positive solution of
    \begin{equation}\label{eq:hnm.Dexp}
        n h_{n,m}^{2\alpha_m(s)+2m+2}e^{-1/h_{n,m}}=1.
    \end{equation}
    Then
    \begin{equation}\label{eq:rate.Dexp.fixed.degree}
        \inf_{0<h\leq\hmax}\sup_{f\in\Sigma(s,L)}
        \pwrisk\bigl(\hat f_{(m,h)},f\bigr)
        \asymp
        h_{n,m}^{\alpha_m(s)}
        \asymp
        (\log n)^{-\alpha_m(s)}.
    \end{equation}
    If $m>\llfloor s\rrfloor$, instead,
    \begin{equation}\label{eq:rate.Dexp.over.degree}
        \inf_{0<h\leq\hmax}\sup_{f\in\Sigma(s,L)}
        \pwrisk\bigl(\hat f_{(m,h)},f\bigr)
        \asymp1.
    \end{equation}
\end{lemma}

\begin{proof}
    \proofpart{Small-degree case (\(m\le\llfloor s\rrfloor\))}
    \proofstep{Reference bandwidth}
    Put
    \[
        p_m(s):=2\alpha_m(s)+2m+2.
    \]
    The map $h\mapsto h^{p_m(s)}e^{-1/h}$ is continuous and strictly
    increasing from zero to infinity on $(0,\infty)$.
    Thus~\eqref{eq:hnm.Dexp} has a unique positive solution, which belongs
    to $(0,\hmax]$ for all sufficiently large $n$.
    With $x_{n,m}:=1/h_{n,m}$, the defining equation becomes
    \[
        x_{n,m}+p_m(s)\log x_{n,m}=\log n.
    \]
    In particular, $x_{n,m}\to\infty$ and $x_{n,m}\le\log n$ for all
    sufficiently large $n$.
    Thus $\log x_{n,m}=o(\log n)$, so $x_{n,m}\sim\log n$ and hence
    \begin{equation}\label{eq:Dexp.hnm.log}
        h_{n,m}\sim\frac1{\log n}.
    \end{equation}

    \proofstep{Upper bound}
    Let $f\in\Sigma(s,L)$ and let
    \[
        q(x,y)=
        \sum_{i+j\le\llfloor s\rrfloor}b_{i,j}x^iy^j
    \]
    be an approximating polynomial.
    Decompose it as
    \[
        q=q_m^-+q_m^+,
        \qquad
        q_m^-:=
        \sum_{i+j\le m}b_{i,j}x^iy^j.
    \]
    The kernel reproduces $q_m^-$.
    Among the monomials of $q_m^+$, Lemma~\ref{lem:Dexp.polynomial.leakage}
    bounds the pure powers of $x$ by a constant times $h^{m+1}$.
    Every monomial containing $y$ is exponentially smaller and hence is
    also $O(h^{m+1})$.
    As in the proof of Proposition~\ref{prop:risk.bound}, the coefficients
    of $q$ are bounded uniformly over $f\in\Sigma(s,L)$.
    Therefore the polynomial contribution to the bias is at most
    $C_{m,s,L}^{\mathrm{poly}}h^{m+1}$.
    This contribution vanishes when $m=\llfloor s\rrfloor$.

    The approximation remainder contributes at most $L \mathcal A_{m,s}(h)$.
    By~\eqref{eq:Dexp.kernel.weighted.L1}, for all sufficiently small $h$,
    \[
        \mathcal A_{m,s}(h)\le C_{m,s}^{\mathrm{rem}}h^{s-m}.
    \]
    Consequently,
    \begin{align*}
        \left|\EE_f\hat f_{(m,h)}(t)-f(t)\right|
        &\le C_{m,s,L}^{\mathrm{poly}}h^{m+1}
        +L \mathcal A_{m,s}(h)\\
        &\le C_{m,s,L}^{\mathrm{bias}}
        \{h^{m+1}+h^{s-m}\}\\
        &\le 2C_{m,s,L}^{\mathrm{bias}}h^{\alpha_m(s)}
    \end{align*}
    for all sufficiently small $h$, with the first term omitted when
    $m=\llfloor s\rrfloor$.
    Moreover, Lemma~\ref{lem:sup.bound},
    \eqref{eq:Dexp.kernel.diagonal}, and~\eqref{eq:volume.Dexp} yield
    \begin{equation*}
        \Var_f\bigl(\hat f_{(m,h)}(t)\bigr)
        \le C_{m,s,L}^{\mathrm{var}}
        \frac{e^{1/h}}{n h^{2m+2}}.
    \end{equation*}
    Taking $h=h_{n,m}$ and using~\eqref{eq:hnm.Dexp} proves the upper
    bound in~\eqref{eq:rate.Dexp.fixed.degree}.

    \proofstep{Bias lower bound for \(h\ge h_{n,m}\)}
    We lower-bound the risk in this range through the bias.
    The complementary range $h\le h_{n,m}$ will be handled through the
    variance.

    \emph{Remainder alternative.}
    We first align an admissible approximation remainder with the sign of the
    equivalent kernel.
    Its mass is compensated on a fixed set lying outside $\Ehmax$.
    Since $\hmax<1/4$, choose a Borel set
    \begin{equation*}
        B\subset\bigl\{(x,y)\in\Dexp:2\hmax<x<3\hmax\bigr\}
    \end{equation*}
    with $0<\mathrm{Leb}(B)<\infty$. In particular, we have
    $B\cap\Ehmax=\varnothing$. Let
    $f_0\equiv\kappa_0:=\mathrm{Leb}(\Dexp)^{-1}$ on $\Dexp$.
    For $h\in(0,\hmax]$, define
    \begin{equation*}
        g_h(u):=\|u\|_\infty^s
        \operatorname{sgn}\!\bigl(k_{m,h}(0,u)\bigr)\1_{\Eh}(u),
        \qquad
        a_h:=\frac{1}{\mathrm{Leb}(B)}\int_{\Eh}g_h(u)\,\rd u,
    \end{equation*}
    with $\operatorname{sgn}(0)=0$. For a sufficiently small constant
    $\varepsilon>0$, independent of $h$, let
    \begin{equation}\label{eq:Dexp.bias.alternative}
        f_h(u):=f_0(u)+\varepsilon
        \bigl\{g_h(u)-a_h\1_B(u)\bigr\},
        \qquad u\in\Dexp.
    \end{equation}
    The definition of $a_h$ gives $\int_{\Dexp}f_h=1$. Moreover,
    \begin{equation*}
        \|g_h\|_\infty\le\hmax^s,
        \qquad
        |a_h|\le
        \frac{\hmax^s\Vscale{\hmax}}{\mathrm{Leb}(B)},
    \end{equation*}
    so choosing $\varepsilon$ small enough ensures $f_h\ge0$ for every
    $h\le\hmax$. On $\Ehmax$, the correction supported on $B$ vanishes,
    and, with the constant approximating polynomial $q\equiv\kappa_0$,
    \begin{equation*}
        |f_h(u)-q(u)|=\varepsilon|g_h(u)|
        \le\varepsilon\|u\|_\infty^s.
    \end{equation*}
    Taking also $\varepsilon\le L$ shows that $f_h\in\Sigma(s,L)$, and
    $f_h(0)=f_0(0)=\kappa_0$.
    The alternative $f_h$ may depend on $h$: the supremum over
    $\Sigma(s,L)$ is taken separately for each estimator
    $\hat f_{(m,h)}$.

    Since $B\cap\Eh=\varnothing$ and the estimator reproduces constants,
    \begin{align}
        \left|\EE_{f_h}\hat f_{(m,h)}(t)-f_h(t)\right|
        &=\varepsilon\int_{\Eh}
        |k_{m,h}(0,u)|\,\mu_h(u)\,
        \|u\|_\infty^s\,\rd u\notag\\
        &=\varepsilon \mathcal A_{m,s}(h).
        \label{eq:Dexp.bias.lower}
    \end{align}
    By~\eqref{eq:Dexp.kernel.weighted.L1}, there are
    $c_{m,s,L}^{\mathrm{rem}}>0$ and $h_0\in(0,\hmax]$ such that the last
    expression is at least
    $c_{m,s,L}^{\mathrm{rem}}h^{s-m}$ for $h\le h_0$.
    On $[h_0,\hmax]$, the function $\mathcal A_{m,s}$ is positive and continuous:
    the local moment matrix is continuous and nonsingular, and the
    integrand in~\eqref{eq:bias.moment.def} is non-zero on a set of positive
    measure.
    Its infimum on this compact interval is therefore positive.

    \emph{Polynomial alternative.}
    When $m<\llfloor s\rrfloor$, a second alternative captures the
    possible leading bias $h^{m+1}$.
    Define $\pi_m(u):=x^{m+1}\1_{\Ehmax}(u)$ and
    \begin{equation*}
        \bar a_m:=\frac1{\mathrm{Leb}(B)}
        \int_{\Dexp}\pi_m(u)\,\rd u,
        \qquad
        \bar f_m:=f_0+\bar\varepsilon
        \{\pi_m-\bar a_m\1_B\},
    \end{equation*}
    where $\bar\varepsilon>0$ is fixed and sufficiently small.
    The same mass-compensation argument makes $\bar f_m$ a non-negative
    probability density.
    It belongs to $\Sigma(s,L)$ because, on $\Ehmax$, its approximating
    polynomial is
    $f_0+\bar\varepsilon x^{m+1}\in\mathcal P_{\llfloor s\rrfloor}$ and
    the approximation error is zero.
    Since $B\cap\Eh=\varnothing$,
    Lemma~\ref{lem:Dexp.polynomial.leakage} gives
    \begin{equation}\label{eq:Dexp.polynomial.bias.lower}
        \left|\EE_{\bar f_m}\hat f_{(m,h)}(t)-\bar f_m(t)\right|
        \ge c_{m,L}^{\mathrm{poly}}h^{m+1},
        \qquad 0<h\le h_0,
    \end{equation}
    after decreasing $h_0$ if necessary.

    The root mean squared risk is at least the absolute bias.
    Hence, when $m<\llfloor s\rrfloor$, the two alternatives give
    \[
        \sup_{f\in\Sigma(s,L)}
        \pwrisk\bigl(\hat f_{(m,h)},f\bigr)
        \ge
        \max\left\{
        c_{m,s,L}^{\mathrm{rem}}h^{s-m},
        c_{m,L}^{\mathrm{poly}}h^{m+1}
        \right\}
        \ge
        c_{m,s,L}^{\mathrm{bias}}h^{\alpha_m(s)}
    \]
    for $0<h\le h_0$.
    When $m=\llfloor s\rrfloor$, one has
    $\alpha_m(s)=s-m$, and the remainder alternative alone gives the same
    conclusion.
    The positive lower bound for $\mathcal A_{m,s}$ on $[h_0,\hmax]$ handles the
    remaining bandwidths.
    Since $h_{n,m}<h_0$ for all sufficiently large $n$, it follows that,
    for every $h\in[h_{n,m},\hmax]$,
    \begin{equation}\label{eq:Dexp.large.h.lower}
        \sup_{f\in\Sigma(s,L)}\pwrisk\bigl(\hat f_{(m,h)},f\bigr)
        \ge c_{m,s,L}^{\mathrm{large}}h_{n,m}^{\alpha_m(s)}.
    \end{equation}

    \proofstep{Variance lower bound for \(h\le h_{n,m}\)}
    Under the constant density $f_0$, the estimator is unbiased and
    \[
        \mu_h(u)^2=V_h^{-1}\mu_h(u).
    \]
    Combining this identity with the squared-norm identity for the
    reproducing kernel gives
    \begin{equation*}
        \Var_{f_0}\bigl(\hat f_{(m,h)}(t)\bigr)
        =
        \frac1n\left\{
        \frac{\kappa_0}{V_h}k_{m,h}(0,0)-\kappa_0^2
        \right\}.
    \end{equation*}
    Relations~\eqref{eq:volume.Dexp} and
    \eqref{eq:Dexp.kernel.diagonal} imply, for all sufficiently small $h$,
    \begin{equation*}
        \Var_{f_0}\bigl(\hat f_{(m,h)}(t)\bigr)
        \ge c_m^{\mathrm{var}}
        \frac{e^{1/h}}{n h^{2m+2}}-\frac{\kappa_0^2}{n}.
    \end{equation*}
    The function
    $h\mapsto e^{1/h}h^{-(2m+2)}$ is decreasing.
    Hence, for $h\le h_{n,m}$,
    \begin{equation*}
        \frac{e^{1/h}}{n h^{2m+2}}
        \ge\frac{e^{1/h_{n,m}}}{n h_{n,m}^{2m+2}}
        =h_{n,m}^{2\alpha_m(s)}.
    \end{equation*}
    Since $n^{-1}=o(h_{n,m}^{2\alpha_m(s)})$, for all large $n$ and every
    $h\le h_{n,m}$,
    \begin{equation}\label{eq:Dexp.small.h.lower}
        \sup_{f\in\Sigma(s,L)}\pwrisk\bigl(\hat f_{(m,h)},f\bigr)
        \ge \sqrt{\Var_{f_0}(\hat f_{(m,h)}(t))}
        \ge c_{m,s,L}^{\mathrm{small}}h_{n,m}^{\alpha_m(s)}.
    \end{equation}
    Combining~\eqref{eq:Dexp.large.h.lower} and
    \eqref{eq:Dexp.small.h.lower} proves the lower bound in
    \eqref{eq:rate.Dexp.fixed.degree}.
    Relation~\eqref{eq:Dexp.hnm.log} proves the last comparison in that
    display.

    \proofpart{Large-degree case (\(m>\llfloor s\rrfloor\))}
    The remainder alternative above is still valid.
    Since now $m\ge s$, its small-bandwidth lower bound
    $c_{m,s,L}^{\mathrm{rem}}h^{s-m}$ is bounded away from zero, while
    the positivity of $\mathcal A_{m,s}$ on compact bandwidth intervals covers the
    remaining values of $h$.
    This proves the lower bound in~\eqref{eq:rate.Dexp.over.degree}.

    For the reverse bound, fix any $h_\circ\in(0,\hmax]$.
    The kernel $k_{m,h_\circ}$ is then fixed, and
    Lemma~\ref{lem:sup.bound} bounds both $f(t)$ and the density throughout
    $\Escale{h_\circ}$ uniformly over $\Sigma(s,L)$.
    The bias and variance of $\hat f_{(m,h_\circ)}$ are therefore bounded
    uniformly over the class, completing the proof.
\end{proof}

Optimising the degree-wise bound yields the best performance attainable
within the fixed-degree family.

\Restate{rethm}{thm:Dexp.best.fixed.degree}{\StmtDexpBest}

\begin{proof}
    Apply Lemma~\ref{lem:Dexp.fixed.degree} to the finitely many
    degrees $0\le m\le M$. Degrees $m>\llfloor s\rrfloor$ have risk of
    constant order, while among $0\le m\le\llfloor s\rrfloor$ the
    smallest risk corresponds to the largest exponent $\alpha_m(s)$.
    Since $m\mapsto m+1$ is increasing and $m\mapsto s-m$ is decreasing,
    their minimum is maximized at an integer closest to their crossing
    point $(s-1)/2$. Finally, if $s>1$, then $\alpha_0(s)=1<s$, whereas
    for $m\ge1$ one has $\alpha_m(s)\le s-m<s$.
\end{proof}

\subsection{Proof of Proposition~\ref{prop:Dexp.anisotropic}}
\label{sec:proof.Dexp.anisotropic}

\Restate{reprop}{prop:Dexp.anisotropic}{\StmtDexpAnisotropic}

\begin{proof}
    \proofpart{Rescaled geometry}
    Take $0<r\leq\hmax/2$ and
    $A_{\bh(r)}=\operatorname{diag}(2r,e^{-1/r})$.
    The inequality $e^{-1/r}\leq2r$ ensures
    $\bh(r)\in(0,\hmax]^2$.
    Changing variables in the defining inequalities of $\Dexp$ gives
    the set $F_r$ displayed in the main text.
    For $1/2<v_1<1$, its vertical section has height one.
    On $0<v_1<1/2$, that height is
    $\exp\{(1-(2v_1)^{-1})/r\}$.
    The substitution $z=2v_1$ therefore gives
    \[
        \mathrm{Leb}(F_r\triangle F_0)
        =\frac12\int_0^1 e^{-(z^{-1}-1)/r}\,\rd z.
    \]
    Since $z^{-1}-1\geq1-z$ for $0<z\leq1$,
    this integral is at most $r/2$.
    The transformed domains lie in the fixed compact set $[0,1]^2$,
    and $\mathrm{Leb}(F_0)=1/2$.
    Corollary~\ref{cor:geometric.stability} applies and gives
    \eqref{eq:Dexp.anisotropic.geometry}.
    In particular, for sufficiently small $r$,
    \[
        V_{\bh(r)}\geq c r e^{-1/r},
        \qquad \delta_m(\bh(r))\geq c_m>0,
        \qquad \Lambda_m(\bh(r))\leq C_m.
    \]

    \proofpart{Risk bound}
    Since $m\geq\llfloor s\rrfloor$, Proposition~\ref{prop:risk.bound}
    applies with $\mathfrak L_{m,s,L}=L$ and $\beta_m(s)=s$.
    Also $\|\bh(r)\|_\infty=2r$, so
    \[
        \mathcal A_{m,s}(\bh(r))
        \leq(2r)^s\Lambda_m(\bh(r))\leq C_m(2r)^s.
    \]
    Lemma~\ref{lem:sup.bound} and the general risk bound yield
    \begin{align*}
        \sup_{f\in\Sigma(s,L)}
        \pwrisk\bigl(\hat f_{(m,\bh(r))},f\bigr)^2
        &\leq L^2\mathcal A_{m,s}(\bh(r))^2
        +\frac{\FsL}{nV_{\bh(r)}\delta_m(\bh(r))}\\
        &\leq C\left(r^{2s}+\frac{e^{1/r}}{nr}\right).
    \end{align*}
    This proves~\eqref{eq:Dexp.anisotropic.risk}.

    \proofpart{Bandwidth choice}
    The choice~\eqref{eq:Dexp.anisotropic.bandwidth} satisfies
    $r_n\to0$, so the risk bound applies for all sufficiently large $n$.
    Moreover,
    \[
        r_n^{2s}=2^{2s}(\log n)^{-2s},
        \qquad
        \frac{e^{1/r_n}}{nr_n}
        =\frac{\log n}{2\sqrt n}
        =o\bigl((\log n)^{-2s}\bigr).
    \]
    Substitution in~\eqref{eq:Dexp.anisotropic.risk} proves
    \eqref{eq:Dexp.anisotropic.rate}.
\end{proof}

\subsection{Proof of Theorem~\ref{thm:Dexp}}\label{sec:proof.Dexp}

\Restate{rethm}{thm:Dexp}{\StmtDexpMinimax}

\begin{proof}
    Proposition~\ref{prop:Dexp.anisotropic} gives the upper bound with
    $m=\llfloor s\rrfloor$. For $0<s\leq1$, the degree-zero estimator
    has bias at most $Lh_n^s$ and variance at most
    $\FsL/(n\Vscale{h_n})\lesssim h_n^{2s}$ by
    \eqref{eq:volume.Dexp} and~\eqref{eq:hn.Dexp}.
    We prove the lower bound by a two-point construction.

    \proofpart{Construction and admissibility}
    We specialise the bump-function construction used in the proof of
    Proposition~\ref{prop:minimax.lower}.
    Choose $\phi_0\in C^\infty_c(\Delta)$ such that
    $0\leq\phi_0\leq1$, $\phi_0(0)=1$, and $\phi_0$ is constant in a
    neighbourhood of the origin.
    Then there exists $C_\phi>0$ such that
    \[
        |\phi_0(v)-1|\leq C_\phi\|v\|_\infty^s,
        \qquad v\in\Delta.
    \]
    Set $\phi_{h_n}(u):=\phi_0(u/h_n)$, where $h_n$ is defined
    by~\eqref{eq:hn.Dexp}.

    Let $f_0\in\Sigma(s,L)$ be a density on $\Dexp$ that is constant on
    $\Ehmax$, with $f_0(u)=f_0(0)>0$ for $u\in\Ehmax$.
    Set
    \[
        c_n
        :=
        \frac{1}{\Vscale{\hmax}}
        \int_{\Escale{h_n}}\phi_{h_n}(u)\,\rd u.
    \]
    Since $0\leq\phi_0\leq1$,
    \[
        0\leq c_n
        \leq\frac{\Vscale{h_n}}{\Vscale{\hmax}}
        \longrightarrow0.
    \]

    By~\eqref{eq:volume.Dexp}, fix $C_{\mathrm{vol}}>0$ such that
    \[
        \Vscale{h}\leq C_{\mathrm{vol}}h^2e^{-1/h}
    \]
    for all sufficiently small $h$.
    Choose
    \[
        0<\varepsilon
        \leq
        \min\left\{
            \frac{L}{C_\phi},\,
            L,\,
            \left(
                \frac{f_0(0)}
                     {2C_{\mathrm{vol}}\|\phi_0\|_\infty^2}
            \right)^{1/2}
        \right\},
        \qquad
        \rho_n:=\varepsilon h_n^s,
    \]
    and define
    \begin{equation}\label{eq:f1.Dexp}
        f_1
        :=
        f_0
        +
        \rho_n
        \left(
            \phi_{h_n}\1_{\Escale{h_n}}
            -
            c_n\1_{\Ehmax}
        \right).
    \end{equation}
    The definition of $c_n$ gives $\int_{\Dexp}f_1=1$.
    Since $\rho_n\to0$, one has $f_1\geq0$ for all sufficiently large
    $n$.
    The verification of membership in $\Sigma(s,L)$ in the proof of Proposition~\ref{prop:minimax.lower}
    applies verbatim and gives $f_0,f_1\in\Sigma(s,L)$.
    Finally,
    \begin{equation}\label{eq:Dexp.lower.signal}
        f_1(0)-f_0(0)
        =
        \rho_n(1-c_n)
        \geq\frac{\rho_n}{2}
    \end{equation}
    for all sufficiently large $n$.

    \proofpart{Kullback--Leibler divergence}
    As in~\eqref{eq:KL.bound},
    \[
        D_{\mathrm{KL}}(f_1\|f_0)
        \leq
        \frac{\rho_n^2}{f_0(0)}
        \|\phi_0\|_\infty^2\Vscale{h_n}.
    \]
    The definition of $C_{\mathrm{vol}}$, the identity
    $\rho_n=\varepsilon h_n^s$, and~\eqref{eq:hn.Dexp} yield
    \begin{align*}
        nD_{\mathrm{KL}}(f_1\|f_0)
        &\leq
        \frac{C_{\mathrm{vol}}\|\phi_0\|_\infty^2}{f_0(0)}
        \varepsilon^2nh_n^{2s+2}e^{-1/h_n} \\
        &=
        \frac{C_{\mathrm{vol}}\|\phi_0\|_\infty^2}{f_0(0)}
        \varepsilon^2
        \leq\frac12.
    \end{align*}

    \proofpart{Conclusion}
    By Pinsker's inequality, Le~Cam's two-point lemma,
    and~\eqref{eq:Dexp.lower.signal},
    \begin{align*}
        \inf_{\hat f}\sup_{f\in\Sigma(s,L)}
        \pwrisk(\hat f,f)^2
        &\geq
        \frac{\rho_n^2}{16}
        \left(
            1-\sqrt{\frac n2D_{\mathrm{KL}}(f_1\|f_0)}
        \right)^2_+ \\
        &\geq
        \frac{\rho_n^2}{64}
        =
        \frac{\varepsilon^2}{64}h_n^{2s}
        \asymp
        (\log n)^{-2s},
    \end{align*}
    since $h_n\sim1/\log n$.
    Taking square roots proves the lower bound in~\eqref{eq:rate.Dexp}.
\end{proof}

\section{Proofs of Section \emph{The connected ladder}}\label{sec:proof.Dlad}

\subsection{Proof of Proposition~\ref{prop:Dlad.volume}}
\label{sec:proof.Dlad.volume}

\Restate{reprop}{prop:Dlad.volume}{\StmtDladVolume}

\begin{proof}
    Recall that
    \[
        k(h)
        :=
        \min\bigl\{
            j\geq2:\zeta_j+\omega_j\leq h
        \bigr\}.
    \]
    \proofpart{Volume decomposition}
    For all sufficiently small $h$, $\Eh$ contains
    $I_j\cap\{0\leq y\leq h\}$ for every $j\geq k(h)$,
    together with a possible slice of $I_{k(h)-1}$ and the residual
    part of the cusp $\Dexp$, up to a Lebesgue-null set.
    Since
    \[
        h
        <
        \zeta_{k(h)-1}+\omega_{k(h)-1},
    \]
    the slice of $I_{k(h)-1}$ has width
    $(h-\zeta_{k(h)-1})_+$.
    The contribution of the cusp satisfies
    \[
        \int_0^h e^{-1/x}\,\rd x
        =
        \int_{1/h}^{\infty}e^{-y}y^{-2}\,\rd y
        =
        O\bigl(h^2e^{-1/h}\bigr).
    \]
    Therefore,
    \[
        \Vh
        =
        h\left[
            \sum_{j\geq k(h)}\omega_j
            +
            \bigl(h-\zeta_{k(h)-1}\bigr)_+
        \right]
        +
        O\bigl(h^2e^{-1/h}\bigr).
    \]

    For $j\geq2$,
    \[
        \frac{\omega_{j+1}}{\omega_j}
        =
        \theta_{j+1}^{1+\nu}
        \leq
        e^{-5}.
    \]
    Hence the sum is dominated by its first term.
    Since $k(h)\to\infty$ as $h\to0$,
    \[
        \sum_{j\geq k(h)}\omega_j
        =
        \omega_{k(h)}
        \bigl(1+O(\theta_{k(h)+1}^{1+\nu})\bigr)
        =
        \omega_{k(h)}\bigl(1+o(1)\bigr).
    \]
    Moreover, the cusp contribution is negligible compared with
    $h\omega_{k(h)}$.
    We have thus obtained the asymptotic equivalent

   \revision{ \[
        V_h
        \sim
        h\omega_{k(h)}
        +
        h^2\left(1-\frac{\zeta_{k(h)-1}}{h}\right)_+.
    \]}

    \proofpart{Two-sided volume bounds}
    By the definition of $k(h)$, $\zeta_{k(h)}\leq h$, and therefore
    \revision{ $
        \omega_{k(h)}
        \leq
        h^{1+\nu}$.}
    If the slice of $I_{k(h)-1}$ is nonempty, then
    \[
        \zeta_{k(h)-1}
        \leq
        h
        <
        \zeta_{k(h)-1}+\omega_{k(h)-1}.
    \]
    Consequently,
    \[
        \left(1-\frac{\zeta_{k(h)-1}}{h}\right)_+
        \leq
        \frac{\omega_{k(h)-1}}{h}
        =
        \zeta_{k(h)-1}^\nu
        \frac{\zeta_{k(h)-1}}{h}
        \leq
        h^\nu.
    \]
    Together with the volume decomposition, these estimates give
    $V_h\leq Ch^{2+\nu}$ for all sufficiently small $h$.

    For the lower bound, the full bar $I_{k(h)}$ gives
   $
       V_h
        \geq
        h\omega_{k(h)}$.
    Since
    \[
        h
        \leq
        \zeta_{k(h)-1}+\omega_{k(h)-1}
        \leq
        2\zeta_{k(h)-1},
    \]
    it follows that
    \[
        \frac{V_h}{h^2}
        \geq
        \frac{\zeta_{k(h)}^{1+\nu}}
             {2\zeta_{k(h)-1}}
        =
        \frac12\,
        \theta_{k(h)}^{1+\nu}
        \zeta_{k(h)-1}^\nu
        \geq
        c\,h^\nu\theta_{k(h)}^{1+\nu}.
    \]
    This proves~\eqref{eq:volume.Dlad}.

    \proofpart{Oscillatory envelopes}
    Since $\log(1/\zeta_j)=j^2$ and
    \[
        \zeta_{k(h)}
        \leq
        h
        \leq
        2\zeta_{k(h)-1},
    \]
    we have
    \[
        \sqrt{\log(1/h)}\leq k(h)
        \leq 1+\sqrt{\log(1/h)+\log2}.
    \]
    Hence $k(h)=\sqrt{\log(1/h)}+O(1)$, with a remainder bounded
    uniformly for all sufficiently small $h$.
    Recalling that $\theta_j=e^{-(2j-1)}$, we obtain
    \[
        \theta_{k(h)}^{1+\nu}
        =
        \exp\left\{
            -2(1+\nu)\sqrt{\log(1/h)}+O(1)
        \right\}
        =
        h^{o(1)}.
    \]
    Substituting into~\eqref{eq:volume.Dlad} and absorbing the bounded
    exponential factor into the constants proves~\eqref{eq:volume.Dlad22}.
    The two-sided bounds also imply
    \revision{\[
        V_h=h^{2+\nu+o(1)}.
    \]
}
    Consider now $h=h_k=2\zeta_k$.
    Since $\omega_k\leq\zeta_k$,
    \[
        \zeta_k+\omega_k
        \leq
        h_k,
    \]
    whereas $\zeta_{k-1}>h_k$ for all sufficiently large $k$.
    Hence $k(h_k)=k$, and the contribution of the full bar $I_k$ gives
    \[
        \frac{V_{h_k}}{h_k^2}
        \geq
        \frac{\zeta_k^{1+\nu}}{2\zeta_k}
        =
        \frac12\zeta_k^\nu
        \asymp
        h_k^\nu.
    \]
    Combined with the general upper bound, this yields
    \begin{equation}\label{eq:Dlad.volume.hk}
        \Vscale{h_k}
        \asymp
        \zeta_k^{2+\nu}.
    \end{equation}
    These scales, free of the deficiency factor
    $\theta_k^{1+\nu}$ appearing in the general lower bound, will be
    used in the minimax arguments below.

    At the opposite envelope, set $\bar h_k=\zeta_{k-1}$.
    Then $k(\bar h_k)=k$ and the second term in the equivalent above
    vanishes, so
    \[
       \frac{ V_{\bar h_k}}{\bar h_k^{2+\nu}}
        \sim\frac{\omega_k}{\bar h_k^{1+\nu}}
        =\theta_k^{1+\nu}\longrightarrow0.
    \]
    Together with the general upper bound and the behaviour at $h_k$, this
    proves
    \[
        \liminf_{h\downarrow0}h^{-(2+\nu)}V_h=0,
        \qquad
        0<\limsup_{h\downarrow0}h^{-(2+\nu)}V_h<\infty.
    \]
    If $\PL(\rho)$ held for some $\rho$, this logarithmic limit would force
    $\rho=\nu$, while the preceding liminf would contradict the required
    lower bound.
\end{proof}

\subsection{Proof of Proposition~\ref{prop:Dlad.PBM}}
\label{sec:proof.Dlad.PBM}

The proposition follows from a two-sided Christoffel estimate and a
bound on the absolute kernel integral along the same sequence.
We first introduce the rescaled measure on which this estimate is based.

\begin{definition}[Rescaled bar--cusp decomposition]
\label{def:Dlad.rescaled.measure}
    Let
    \[
        \overline E_k
        :=
        \{v\in\R^2:h_kv\in\Escale{h_k}\},
    \]
    and let $\mu_k$ be the push-forward of
    $\mu_{h_k}$ under the map $u\mapsto u/h_k$.
    For $j\geq0$, set
    \[
        \rho_{j,k}:=\frac{\zeta_{k+j}}{2\zeta_k},
        \qquad
        \tau_{j,k}:=\frac{\omega_{k+j}}{2\zeta_k},
        \qquad
        R_{j,k}
        :=
        [\rho_{j,k},\rho_{j,k}+\tau_{j,k}]\times[0,1].
    \]
    We assign the overlap between the bars and the exponential cusp to
    the bar component by defining
    \[
        P_k:=\bigcup_{j\geq0}R_{j,k},
        \qquad
        C_k:=\overline E_k\setminus P_k.
    \]
    Thus $P_k$ contains the complete portions of the bars present in
    the local neighbourhood. Up to a Lebesgue-null set, $C_k$ is the
    residual part of the rescaled cusp not covered by a bar.
    Both sets have positive Lebesgue measure.

    Let
    \[
        \eta_k:=\mu_k(C_k),
    \]
    and, for every Borel set $A\subseteq\R^2$, define the conditional
    probability measures
    \[
        \mu_{P,k}(A)
        :=
        \frac{\mu_k(A\cap P_k)}{1-\eta_k},
        \qquad
        \mu_{\mathrm{cusp},k}(A)
        :=
        \frac{\mu_k(A\cap C_k)}{\eta_k}.
    \]
    By construction,
    \[
        \mu_k
        =
        (1-\eta_k)\mu_{P,k}
        +
        \eta_k\mu_{\mathrm{cusp},k}.
    \]
    Finally, set
    \[
        w_{j,k}
        :=
        \frac{\omega_{k+j}}{\sum_{i\geq0}\omega_{k+i}},
        \qquad
        \mu_{1,k}
        :=
        \sum_{j\geq0}w_{j,k}\,
        \mathcal U[\rho_{j,k},\rho_{j,k}+\tau_{j,k}].
    \]
\end{definition}

When no confusion can arise, we suppress the second index and write
$\rho_j$, $\tau_j$, $w_j$ and $\mu_1$ for
$\rho_{j,k}$, $\tau_{j,k}$, $w_{j,k}$ and $\mu_{1,k}$, respectively.
Here and below, a subscript $\nu$ on $\lesssim$ or $\asymp$ indicates
that the comparison constants may depend on $\nu$, but not on $k$.

\begin{lemma}[Geometry of the rescaled local measure at $h_k$]
\label{lem:Dlad.rescaled.measure}
    With the notation of
    Definition~\ref{def:Dlad.rescaled.measure}, for every $m\in\N_0$,
    \[
        \delta_m(h_k)
        =
        \inf_{p\in\Zm}\|1-p\|_{L^2(\mu_k)}^2.
    \]
    Moreover,
    \[
        \mu_{P,k}
        =
        \mu_{1,k}\otimes\mathcal U[0,1],
    \]
    and, for all sufficiently large $k$,
    \begin{equation}\label{eq:Dlad.eta}
        \eta_k\leq e^{-1/(2h_k)}.
    \end{equation}
    In addition,
    $\rho_{0,k}=1/2$, $\tau_{0,k}=\zeta_k^\nu/2$, and
    \begin{equation}\label{eq:Dlad.weights}
        w_{0,k}\geq1-2\theta_{k+1}^{1+\nu},
        \qquad
        w_{1,k}\asymp_\nu\theta_{k+1}^{1+\nu},
        \qquad
        \sum_{j\geq1}w_{j,k}\leq2\theta_{k+1}^{1+\nu}.
    \end{equation}
    Finally, for every $j\geq1$,
    $R_{j,k}\subseteq[0,\theta_{k+1}]\times[0,1]$, since
    \begin{equation}\label{eq:Dlad.bands}
        \rho_{1,k}+\tau_{1,k}
        =
        \frac{\zeta_{k+1}+\omega_{k+1}}{2\zeta_k}
        \leq\theta_{k+1}.
    \end{equation}
\end{lemma}

\begin{proof}
    \proofpart{Variational identity}
    The correspondence $p(u)\mapsto p(h_kv)$ is a bijection of $\Zm$
    onto itself.
    The variational definition of $\delta_m(h_k)$ therefore gives the
    first identity.

    \proofpart{Cusp mass}
    Up to a Lebesgue-null set, $C_k$ is contained in the rescaled
    exponential cusp.
    Before rescaling, this cusp has area at most $h_ke^{-1/h_k}$,
    whereas bar $I_k$ contributes $h_k\omega_k$ to the local volume.
    Hence
    \[
        \eta_k
        \leq
        \frac{e^{-1/h_k}}{\omega_k}
        =
        \exp\{-1/h_k+(1+\nu)k^2\}.
    \]
    Since $k^2=\log(2/h_k)$, the exponent is at most $-1/(2h_k)$ for
    all sufficiently large $k$, which proves~\eqref{eq:Dlad.eta}.

    \proofpart{Product decomposition}
    Rescaling $I_{k+j}\cap\Escale{h_k}$ by $h_k=2\zeta_k$ gives
    $R_{j,k}$.
    Its mass before normalisation is proportional to
    $h_k\omega_{k+j}$.
    This proves the product decomposition of $\mu_{P,k}$ and the
    formula for the weights.
    The identities for $\rho_{0,k}$ and $\tau_{0,k}$ follow directly from
    $h_k=2\zeta_k$ and $\omega_k=\zeta_k^{1+\nu}$.

    \proofpart{Bar weights and locations}
    Set $r_k:=\theta_{k+1}^{1+\nu}$.
    Since
    \[
        \frac{\omega_{i+1}}{\omega_i}
        =
        \theta_{i+1}^{1+\nu}
        \leq e^{-5},
        \qquad i\geq2,
    \]
    the tail of the bar masses satisfies
    \[
        \sum_{j\geq1}\frac{\omega_{k+j}}{\omega_k}
        \leq
        \frac{r_k}{1-e^{-5}}
        \leq
        2r_k.
    \]
    Consequently,
    $\sum_{j\geq1}w_{j,k}\leq2r_k$ and
    $w_{0,k}\geq1-2r_k$.
    Moreover, $w_{1,k}=r_kw_{0,k}\asymp_\nu r_k$, which
    proves~\eqref{eq:Dlad.weights}.

    The sequence $\zeta_j+\omega_j$ is decreasing.
    Thus every band with $j\geq1$ lies to the left of the band with
    $j=1$.
    Finally,
    \[
        \rho_{1,k}+\tau_{1,k}
        =
        \frac{\theta_{k+1}}{2}
        \bigl(1+\zeta_{k+1}^{\nu}\bigr)
        \leq
        \theta_{k+1},
    \]
    which proves~\eqref{eq:Dlad.bands}.
\end{proof}

\begin{lemma}[Collapse of $\delta_1$ at $h_k$]\label{lem:Dlad.delta1}
    Let $m_1$ and $\sigma_1^2$ be the mean and variance of the marginal
    measure $\mu_1$, which is the first-coordinate marginal of
    $\mu_{P,k}$ by Lemma~\ref{lem:Dlad.rescaled.measure}.
    Then
    \[
        m_1\longrightarrow\frac12,
        \qquad
        \sigma_1^2
        \asymp_\nu
        \zeta_k^{2\nu}+\theta_{k+1}^{1+\nu}
        \asymp_\nu
        \theta_{k+1}^{1+\nu}.
    \]
    Moreover, there exist constants
    $0<c_\nu^\delta\leq C_\nu^\delta<\infty$ such that, for all
    sufficiently large $k$,
    \begin{equation}\label{eq:Dlad.delta1}
        c_\nu^\delta\theta_{k+1}^{1+\nu}
        \leq
        \delta_1(h_k)
        \leq
        C_\nu^\delta\theta_{k+1}^{1+\nu}.
    \end{equation}
\end{lemma}

\begin{proof}
    \proofpart{Mean}
    Set $w_+:=\sum_{j\geq1}w_j$.
    The mean of the leading band is
    \[
        \overline v_0
        =
        \rho_0+\frac{\tau_0}{2}
        =
        \frac12+\frac14\zeta_k^\nu.
    \]
    Let $\overline v_+$ be the conditional mean of the remaining bands.
    Relation~\eqref{eq:Dlad.bands} gives
    $0\leq\overline v_+\leq\theta_{k+1}$.
    Hence~\eqref{eq:Dlad.weights} gives
    \[
        m_1
        =
        w_0\overline v_0+w_+\overline v_+
        =
        \frac12
        +
        O\bigl(\zeta_k^\nu+\theta_{k+1}^{1+\nu}\bigr),
    \]
    and therefore $m_1\to1/2$.

    \proofpart{Variance}
    To control the variance from above, centre at $1/2$ rather than at
    $m_1$.
    Since the leading band is uniform and all the remaining bands are
    contained in $[0,\theta_{k+1}]$,
    \[
        \sigma_1^2
        \leq
        \int\left(v_1-\frac12\right)^2\,\rd\mu_1(v_1)
        \leq
        \frac{\tau_0^2}{3}+\frac{w_+}{4}
        \lesssim_\nu
        \zeta_k^{2\nu}+\theta_{k+1}^{1+\nu}.
    \]
    Conversely, the law of total variance, first between the leading
    band and the union of the other bands, yields
    \[
        \sigma_1^2
        \geq
        \frac{w_0\tau_0^2}{12}
        +
        w_0w_+
        \bigl(\overline v_0-\overline v_+\bigr)^2.
    \]
    For all sufficiently large $k$, the first term is bounded below by
    a positive multiple of $\zeta_k^{2\nu}$.
    The second is bounded below by a positive multiple of
    $w_1\asymp_\nu\theta_{k+1}^{1+\nu}$, because
    $\overline v_0-\overline v_+\to1/2$.
    Thus
    \[
        \sigma_1^2
        \asymp_\nu
        \zeta_k^{2\nu}+\theta_{k+1}^{1+\nu}.
    \]

    \proofpart{One-dimensional Christoffel function}
    For a probability measure $\mu$ on $\R$, define
    \[
        \delta_1^{(1)}(\mu)
        :=
        \inf_{a\in\R}\int(1-av)^2\,\rd\mu(v).
    \]
    Direct minimisation in $a$ gives
    \[
        \delta_1^{(1)}(\mu_1)
        =
        \frac{\sigma_1^2}{m_1^2+\sigma_1^2}
        \asymp_\nu
        \sigma_1^2,
    \]
    where the last comparison follows from $m_1\to1/2$.

    \proofpart{Comparison with \(\delta_1(h_k)\)}
    It remains to compare this one-dimensional quantity with
    $\delta_1(h_k)$.
    Let $p(v)=av_1+bv_2\in\mathcal Z_1$.
    For each fixed $v_2$, the definition of
    $\delta_1^{(1)}(\mu_1)$ gives
    \[
        \int\{1-p(v_1,v_2)\}^2\,\rd\mu_1(v_1)
        \geq
        (1-bv_2)^2\delta_1^{(1)}(\mu_1).
    \]
    Indeed, when $1-bv_2\neq0$, one factors $1-bv_2$ from the residual
    and uses the admissible slope $a/(1-bv_2)$. When it vanishes, the
    lower bound is immediate.
    Integrating with respect to $v_2\sim\mathcal U[0,1]$ and using
    \[
        \inf_{b\in\R}\int_0^1(1-bv_2)^2\,\rd v_2=\frac14
    \]
    shows that
    \[
        \delta_1(h_k)
        \geq
        \frac{1-\eta_k}{4}\,
        \delta_1^{(1)}(\mu_1).
    \]
    For the reverse inequality, use the test polynomial $p(v)=2v_1$.
    Its squared residual is at most $\zeta_k^{2\nu}$ on the leading
    band and at most $1$ on every other bar.
    It is also at most $1$ on the rescaled cusp, where
    $0\leq v_1\leq1$.
    Therefore,
    \[
        \delta_1(h_k)
        \leq
        \zeta_k^{2\nu}
        +
        2\theta_{k+1}^{1+\nu}
        +
        \eta_k.
    \]
    Relation~\eqref{eq:Dlad.eta} makes the last term negligible.
    Finally,
    \[
        \zeta_k^{2\nu}
        =
        o\bigl(\theta_{k+1}^{1+\nu}\bigr),
    \]
    because
    $2\nu k^2-(1+\nu)(2k+1)\to\infty$.
    Combining the preceding bounds proves~\eqref{eq:Dlad.delta1} and
    the stated estimates for $\sigma_1^2$.
\end{proof}

\Restate{reprop}{prop:Dlad.PBM}{\StmtDladPBM}

\begin{proof}
    Fix $m\geq1$.
    For all sufficiently large $k$,
    Lemma~\ref{lem:Dlad.delta1} applies and gives the second inequality
    below.
    The first follows from
    $\mathcal Z_1\subseteq\mathcal Z_m$:
    \[
        \delta_m(h_k)
        \leq
        \delta_1(h_k)
        \leq
        C_\nu^\delta\theta_{k+1}^{1+\nu}.
    \]
    This proves~\eqref{eq:delta.Dlad}.
    Since $\theta_{k+1}^{1+\nu}\to0$,
    Proposition~\ref{prop:PBM} then shows that \PBM[m] fails.

    To control the absolute kernel integral, set
    $r_k=\theta_{k+1}^{1+\nu}$ and use the rescaled measures of
    Definition~\ref{def:Dlad.rescaled.measure}.
    The affine reproducing kernel at zero for the product measure
    $\mu_{P,k}=\mu_{1,k}\otimes\mathcal U[0,1]$ is
    \[
        q_{P,k}(v)
        =1+\frac{m_1(m_1-v_1)}{\sigma_1^2}+3-6v_2.
    \]
    The leading band has width $\tau_{0,k}=\zeta_k^\nu/2=o(r_k)$,
    while all other bands have total mass $O(r_k)$.
    Hence $|m_1-1/2|=O(r_k)$ and
    \[
        \int|v_1-m_1|\,\rd\mu_{1,k}(v_1)=O(r_k).
    \]
    Since $\sigma_1^2\asymp r_k$, this gives
    $\int|q_{P,k}|\,\rd\mu_{P,k}\lesssim1$.

    Let $G_{P,k}$ and $G_k$ be the Gram matrices in the basis
    $\Psi(v)=(1,v_1,v_2)^\top$ under $\mu_{P,k}$ and $\mu_k$.
    The product structure and Lemma~\ref{lem:Dlad.delta1} give
    $\lambda_{\min}(G_{P,k})\asymp r_k$.
    All measures are supported in $[0,1]^2$, so
    $\|G_k-G_{P,k}\|_2\lesssim\eta_k$.
    By~\eqref{eq:Dlad.eta}, $\eta_k=o(r_k^2)$, and the inverse identity
    yields
    \[
        \|G_k^{-1}-G_{P,k}^{-1}\|_2
        \lesssim\eta_k r_k^{-2}=o(1).
    \]
    The full affine kernel $q_k(v)=e_1^\top G_k^{-1}\Psi(v)$ therefore
    differs from $q_{P,k}$ by $o(1)$ uniformly on $[0,1]^2$.
    Moreover, $\sup_{[0,1]^2}|q_{P,k}|\lesssim r_k^{-1}$.
    The mixture decomposition of $\mu_k$ now gives
    \[
        \Lambda_1(h_k)=\int|q_k|\,\rd\mu_k
        \leq\int|q_{P,k}|\,\rd\mu_{P,k}
        +O(\eta_k/r_k)+o(1)\lesssim1.
    \]
\end{proof}

\subsection{Proof of Theorem~\ref{thm:Dlad.low}}
\label{sec:proof.Dlad.low}

\Restate{rethm}{thm:Dlad.low}{\StmtDladLow}

\begin{proof}
    \proofpart{The implicit bandwidth}
    The map $h\mapsto\Vh$ is continuous, positive and non-decreasing on
    $(0,\hmax]$.
    Consequently, $F(h):=h^{2s}\Vh$ is continuous and strictly
    increasing, tends to $0$ as $h\downarrow0$, and is positive at
    $\hmax$.
    Thus~\eqref{eq:hn.implicit.Dlad} has a unique solution whenever
    $n\hmax^{2s}\Vscale{\hmax}\geq1$.
    Otherwise, the theorem defines $h_n=\hmax$.
    The next two parts first treat the former case with $h_n$
    sufficiently small, and hence all sufficiently large $n$.

    \proofpart{Upper bound}
    Since $0<s\leq1$, the approximating polynomial in
    Definition~\ref{def:Sigma} is the constant $f(t)$.
    The degree-zero estimator is
    \[
        \hat f_{(0,h)}(t)
        =
        \frac{1}{n\Vh}
        \sum_{i=1}^n\1_{\Eh}(X_i-t).
    \]
    Its bias is at most $Lh^s$, while
    Lemma~\ref{lem:sup.bound} gives
    \[
        \Var\bigl(\hat f_{(0,h)}(t)\bigr)
        \leq
        \frac{\FsL}{n\Vh}.
    \]
    At $h=h_n$, Relation~\eqref{eq:hn.implicit.Dlad} therefore yields
    \[
        \sup_{f\in\Sigma(s,L)}
        \pwrisk\bigl(\hat f_{(0,h_n)},f\bigr)
        \leq
        C_{s,L}^{\mathrm{up}}h_n^s,
        \qquad
        C_{s,L}^{\mathrm{up}}
        :=
        \bigl(L^2+\FsL\bigr)^{1/2}.
    \]

    \proofpart{Lower bound}
    Let
    $f_0\equiv\kappa_0:=\mathrm{Leb}(\Dlad)^{-1}$.
    Choose a smooth function $\phi_0$ supported in $\Delta$ such that
    \[
        0\leq\phi_0\leq1,
        \qquad
        \phi_0(0)=1,
        \qquad
        |\phi_0(v)-1|
        \leq
        C_\phi\|v\|_\infty^s.
    \]
    Put $\phi_h(u):=\phi_0(u/h)$ and
    \begin{align*}
        c_h
        &:=
        \mathrm{Leb}(I_2)^{-1}
        \int_{\Eh}\phi_h(u)\,\rd u,\\
        f_1(t+u)
        &:=
        f_0(t+u)
        +
        \rho_h
        \bigl\{
            \phi_h(u)\1_{\Eh}(u)
            -
            c_h\1_{I_2}(t+u)
        \bigr\},
        \qquad
        \rho_h:=\varepsilon h^s,
    \end{align*}
    where $I_2$ is disjoint from $\Eh$ for all sufficiently small $h$.

    We verify that the two alternatives are admissible.
    The definition of $c_h$ gives $\int f_1=1$, and
    \[
        0
        \leq
        c_h
        \leq
        \frac{\Vh}{\mathrm{Leb}(I_2)}
        \longrightarrow0.
    \]
    Hence, for small $h$,
    $f_1\geq\kappa_0-\rho_hc_h\geq\kappa_0/2$.
    Next, take $q_h\equiv\kappa_0+\rho_h$.
    On $\Eh$,
    \[
        |f_1(t+u)-q_h|
        =
        \rho_h|\phi_0(u/h)-1|
        \leq
        \varepsilon C_\phi\|u\|_\infty^s.
    \]
    On $\Ehmax\setminus\Eh$, one has $\|u\|_\infty\geq h$ and
    \[
        |f_1(t+u)-q_h|
        \leq
        \rho_h(1+c_h).
    \]
    Thus, after fixing $\varepsilon>0$ sufficiently small in terms of
    $s$ and $L$, both $f_0$ and $f_1$ belong to $\Sigma(s,L)$.
    Moreover,
    \[
        f_1(t)-f_0(t)=\rho_h.
    \]

    Since $f_0=\kappa_0$ on the support of the perturbation and
    \[
        c_h^2\mathrm{Leb}(I_2)
        \leq
        \frac{\Vh^2}{\mathrm{Leb}(I_2)},
    \]
    there exists a finite constant $C_{\mathrm{KL}}$, independent of
    $h$ and $n$, such that
    \[
        nD_{\mathrm{KL}}(f_1\|f_0)
        \leq
        C_{\mathrm{KL}}n\rho_h^2\Vh.
    \]
    At $h=h_n$, the right-hand side is at most
    $C_{\mathrm{KL}}\varepsilon^2$.
    Decreasing $\varepsilon$ once more if necessary, Le~Cam's lemma
    gives a constant $c_{s,L}^{\mathrm{low}}>0$ such that
    \[
        \mathcal R_n(s,L)
        \geq
        c_{s,L}^{\mathrm{low}}h_n^s.
    \]
    This proves~\eqref{eq:rate.implicit.Dlad} for all sufficiently
    large $n$.

    \proofpart{The remaining sample sizes}
    Only finitely many integers $n\geq3$ remain.
    For each of them, the preceding bias and variance bounds applied
    at the bandwidth $h_n$ give a finite upper-bound constant.
    Their maximum therefore proves the upper bound in
    \eqref{eq:rate.implicit.Dlad} on this finite set, still using
    $\hat f_{(0,h_n)}$.
    For the lower bound, use the same two-point construction at a fixed
    sufficiently small bandwidth and choose its amplitude uniformly
    over the remaining sample sizes.
    Comparing the resulting positive lower bound with the finitely many
    values $h_n^s$ completes the proof of
    \eqref{eq:rate.implicit.Dlad} for every $n\geq3$.

    \proofpart{Oscillation envelope}
    Set
    \[
        x_n:=\log(1/h_n),
        \qquad
        N:=\log n,
        \qquad
        a:=2s+2+\nu,
        \qquad
        D(h):=-\log\{V_h{}/h^{2+\nu}\}.
    \]
    The defining equation for $h_n$ becomes
    \begin{equation}\label{eq:Dlad.Nx}
        N=ax_n+D(h_n).
    \end{equation}
    Relation~\eqref{eq:volume.Dlad} and the definition of $k(h)$ imply
    that there exists $C_\nu^D<\infty$ such that, uniformly as
    $h\downarrow0$,
    \[
        -C_\nu^D
        \leq
        D(h)
        \leq
        2(1+\nu)\sqrt{\log(1/h)}\,(1+o(1)).
    \]
    It follows from~\eqref{eq:Dlad.Nx} that
    \[
        x_n
        =
        \frac{N}{a}
        +
        O(\sqrt N),
        \qquad
        x_n\sim\frac{N}{a}.
    \]
    Since~\eqref{eq:rate.implicit.Dlad} gives
    $\log\mathcal R_n(s,L)=-sx_n+O(1)$,
    Relation~\eqref{eq:Dlad.Nx} yields
    \[
        \frac{
            \log\mathcal R_n(s,L)
            +
            \frac{s}{a}\log n
        }{\sqrt{\log n}}
        =
        \frac{sD(h_n)}{a\sqrt N}
        +
        o(1).
    \]
    The lower bound on $D$ shows that the liminf is at least $0$.
    The upper bound, together with $x_n\sim N/a$, shows that the limsup
    is at most
    \[
        \frac{2s(1+\nu)}{a^{3/2}}.
    \]

    \proofpart{Subsequences attaining the two envelopes}
    At $h_k^-:=2\zeta_k$,
    \[
        V_{h_k^-}\asymp(h_k^-)^{2+\nu},
        \qquad
        D(h_k^-)=O(1).
    \]
    At $h_k^+:=\zeta_{k-1}$, $\Escale{h_k^+}$ contains no portion of
    bar $k-1$ with positive area,
    and the proof of Proposition~\ref{prop:Dlad.volume} gives
    \[
       V_{h_k^+}
        \asymp
        \omega_k h_k^+,
        \qquad
        V_{h_k^+}
        \asymp (h_k^+)^{2+\nu}
        \theta_k^{1+\nu},
        \qquad
        D(h_k^+)
        =
        (1+\nu)(2k-1)+O(1).
    \]

    For either choice $h=h_k^\pm$, define
    \[
        n_k^\pm
        :=
        \left\lceil F(h_k^\pm)^{-1}\right\rceil.
    \]
    Then
    \[
        n_k^\pm F(h_k^\pm)=1+o(1),
        \qquad
        \frac{F(h_{n_k^\pm})}{F(h_k^\pm)}=1+o(1),
        \qquad
        h_{n_k^\pm}\leq h_k^\pm.
    \]
    We claim that
    \[
        \frac{h_{n_k^\pm}}{h_k^\pm}\longrightarrow1.
    \]
    Fix $\epsilon\in(0,1/2)$.
    For all sufficiently large $k$, neither interval
    $[(1-\epsilon)h_k^\pm,h_k^\pm]$ crosses a relevant bar boundary:
    bar $k$ remains fully included at the minus scales, while the plus
    interval remains in the gap preceding bar $k-1$.
    The volume calculation in the proof of
    Proposition~\ref{prop:Dlad.volume} then gives, uniformly for
    $r\in[1-\epsilon,1]$,
    \[
        \Vscale{r h_k^\pm}
        =
        (1+o(1))\,r h_k^\pm\omega_k.
    \]
    Consequently,
    \[
        \frac{F((1-\epsilon)h_k^\pm)}{F(h_k^\pm)}
        \longrightarrow
        (1-\epsilon)^{2s+1}
        <
        1.
    \]
    The monotonicity of $F$ and
    $F(h_{n_k^\pm})/F(h_k^\pm)\to1$ prove the claim.

    At the plus scales, $h_{n_k^+}\leq h_k^+$ also ensures that the
    solution stays before bar $k-1$ enters.
    Hence
    \[
        D(h_{n_k^-})=O(1),
        \qquad
        D(h_{n_k^+})
        =
        (1+\nu)(2k-1)+o(k).
    \]
    The same volume equivalents give
    \[
        \log n_k^\pm
        =
        ak^2+O(k).
    \]
    Substituting these relations into the normalised-risk identity
    above proves, respectively,
    \eqref{eq:oscillations.Dlad.inf} and
    \eqref{eq:oscillations.Dlad.sup}.
    Since $D(h_{n_k^-})=O(1)$, the identity
    $\log n_k^-=a\log(1/h_{n_k^-})+D(h_{n_k^-})$
    also gives $h_{n_k^-}^s\asymp(n_k^-)^{-s/a}$.
    Equation~\eqref{eq:rate.implicit.Dlad} proves the asserted baseline order.
\end{proof}

\subsection{Proof of Theorem~\ref{thm:Dlad.anisotropic}}
\label{sec:proof.Dlad.anisotropic}

\Restate{rethm}{thm:Dlad.anisotropic}{\StmtDladAnisotropic}

\begin{proof}
    Fix $1<s\leq2$ and $L>0$. Throughout the proof, write
    \[
        b=1+\nu,\qquad c=2s+1,\qquad a=b+c,
        \qquad B_k=B_k(s),\qquad v_k(n)=v_k(n,s).
    \]
    All comparison constants are independent of $n$ and $k$.

    \proofpart{The minimising index}
    Define the real thresholds
    \begin{equation}\label{eq:Dlad.high.thresholds}
        N_k=B_k^{-c/s}\omega_k^{-1},
        \qquad M_k=B_k^{-c/s}\omega_{k+1}^{-1}.
    \end{equation}
    The sequence $B_k$ is strictly decreasing and $v_k(n)$ is strictly
    increasing in $k$. Moreover,
    $M_{k-1}<N_k<M_k<N_{k+1}$. Direct comparison of the candidates gives
    \begin{equation}\label{eq:Dlad.high.phases}
        \mathfrak r_n(s)=
        \begin{cases}
            v_k(n),&M_{k-1}\leq n\leq N_k,\\
            B_k,&N_k\leq n\leq M_k,
        \end{cases}
    \end{equation}
    for all sufficiently large $n$.
    Indeed, on the first interval $B_k\leq v_k(n)\leq B_{k-1}$.
    Every smaller index has a bias term at least $B_{k-1}$, and every
    larger index has a variance term at least $v_k(n)$.
    On the second interval, the $k$th candidate equals $B_k$.
    Smaller indices have a larger bias term, while larger indices have
    a variance term at least $v_{k+1}(n)\geq B_k$.
    The infimum is attained because $v_k(n)\to\infty$ as $k\to\infty$.

    For $k=k_n$, these comparisons imply $M_{k-1}\leq n\leq M_k$
    once $n$ is sufficiently large. Since
    \begin{align}
        \log B_k&=-sk^2+2(s-1)k-(s-1),\label{eq:Dlad.high.logB}\\
        \log N_k&=ak^2-\frac{2c(s-1)}s k
                       +\frac{c(s-1)}s,\label{eq:Dlad.high.logN}\\
        \log M_k&=\log N_k+b(2k+1),\label{eq:Dlad.high.logM}
    \end{align}
    we have
    \begin{equation}\label{eq:Dlad.high.active.scales}
        k_n\longrightarrow\infty,\qquad
        \log n=ak_n^2+O(k_n),\qquad
        \log H_{k_n,n}=-k_n^2+O(k_n).
    \end{equation}
    In particular, both components of $\bh_{k_n,n}$ tend to zero.

    \proofpart{Affine projection with the anisotropic bandwidth}
    Take $k=k_n$ and abbreviate
    \[
        r=\zeta_k,\quad R=\zeta_{k-1},\quad
        w=\omega_k=r^b,\quad \theta=r/R,\quad
        \ell=\theta w,\quad H=H_{k,n}.
    \]
    For large $k$, $\ell<\omega_{k-1}$ and $R+\ell<2R$.
    The portion of $\Escale{\bh_{k,n}}$ supplied by the bars is the product
    \[
        P=A_k\times[0,H],\qquad
        A_k=[R,R+\ell]\cup
        \bigcup_{j\geq k}[\zeta_j,\zeta_j+\omega_j].
    \]
    Its horizontal length is $W_k=\ell+\sum_{j\geq k}\omega_j\sim w$, since $\ell=o(w)$ and $\sum_{j>k}\omega_j=o(w)$.
    Let $\mu_x$ and $\sigma_x^2$ be the mean and variance of the
    uniform probability measure on $A_k$.
    The $k$th bar has a mass bounded away from zero and lies at
    $x\asymp r$. The outer slice has mass of order $\theta$ and lies
    at $x\asymp R$. All remaining bars lie below $2r$.
    The outer slice contributes order $\theta R=r$ to the first
    moment and order $\theta R^2=rR$ to the second. The inner bars
    contribute $O(r^2)$ to the second moment. The variance between
    the outer slice and the dominant bar is bounded below by a
    constant times $\theta R^2$. Consequently,
    \begin{equation}\label{eq:Dlad.high.horizontal.moments}
        \mu_x\asymp r,\qquad
        \sigma_x^2\asymp rR,\qquad
        \frac{\mu_x^2}{\sigma_x^2}=O(\theta).
    \end{equation}
    Under the uniform probability measure on $P$, the two coordinates
    are independent. The affine reproducing kernel at zero is
    \[
        \mathcal K_P(x,y)
        =1+\frac{\mu_x(\mu_x-x)}{\sigma_x^2}+3-\frac{6y}{H}.
    \]
    Its coefficients in the basis $(1,x/R,y/H)$ are bounded uniformly,
    by~\eqref{eq:Dlad.high.horizontal.moments}.

    We next include the cusp. Its contribution outside $P$ has area
    at most $2R e^{-1/(2R)}$. If $\eta$ denotes its mass under the
    uniform probability measure on $\Escale{\bh_{k,n}}$, then
    \[
        \eta\leq\frac{2R e^{-1/(2R)}}{HW_k}=o(\theta),
    \]
    by~\eqref{eq:Dlad.high.active.scales}.
    To control the kernel, let $G_P$ and $G$ be the Gram matrices in
    the basis $\Psi=(1,x/R,y/H)^\top$ under the uniform probabilities
    on $P$ and on $\Escale{\bh_{k,n}}$, respectively.
    The coordinates of $\Psi$ are bounded on $\Escale{\bh_{k,n}}$, so
    $\|G-G_P\|_2\leq C\eta$.
    Independence and~\eqref{eq:Dlad.high.horizontal.moments} give
    $\lambda_{\min}(G_P)\asymp\theta$.
    Hence $\|G^{-1}\|_2\lesssim\theta^{-1}$, and
    \[
        \|G^{-1}e_1-G_P^{-1}e_1\|_2
        \leq\|G^{-1}\|_2\|G-G_P\|_2\|G_P^{-1}e_1\|_2
        \lesssim\eta/\theta=o(1).
    \]
    The full kernel $\mathcal K(x,y)=e_1^\top G^{-1}\Psi(x,y)$ is
    therefore uniformly bounded on $\Escale{\bh_{k,n}}$. It follows that
    \begin{equation}\label{eq:Dlad.high.kernel.bounds}
        V_{\bh_{k,n}}\asymp H\omega_k,\qquad
        \delta_1(\bh_{k,n})\asymp1,\qquad
        \Lambda_1(\bh_{k,n})\lesssim1.
    \end{equation}
    Here the lower bound for $\mathcal K(0,0)$ follows from reproduction
    of constants. Its upper bound follows from the bounded coefficients.

    Under the uniform probability measure on $\Escale{\bh_{k,n}}$, the bars with index at least
    $k$ have $x\leq2r$, while the outer slice has mass $O(\theta)$.
    The cusp has mass $o(\theta)$ and $x\leq2R$.
    The boundedness of $\mathcal K$ therefore gives
    \begin{align}
        \mathcal A_{1,s}(\bh_{k,n})
        &\lesssim(r\vee H)^s+\theta(R\vee H)^s\notag\\
        &\lesssim H^s+rR^{s-1}=H^s+B_k.
        \label{eq:Dlad.high.bias.moment}
    \end{align}
    Proposition~\ref{prop:risk.bound} and Lemma~\ref{lem:sup.bound} yield
    \[
        \sup_{f\in\Sigma(s,L)}\pwrisk(\hat f_{(1,\bh_{k,n})},f)
        \lesssim H^s+B_k+(nH\omega_k)^{-1/2}
        =B_k+2v_k(n)\lesssim\mathfrak r_n(s).
    \]
    This proves the upper bound for the estimator specified in the theorem.

    \proofpart{A lower bound between successive bias levels}
    Let $j$ be large and take
    \[
        B_{j+1}\leq T\leq B_j,\qquad
        r=\zeta_j,\qquad R=\zeta_{j-1},\qquad H=T^{1/s}.
    \]
    Choose a sufficiently small constant $\epsilon>0$, depending only
    on $s$ and $L$, and define on $\Dlad$
    \[
        p(x,y)=\epsilon T(1-x/r),\qquad
        g(x,y)=\operatorname{sign}(p(x,y))
        \bigl(|p(x,y)|-(L/2)\|(x,y)\|_\infty^s\bigr)_+.
    \]
    Then $g(0)=\epsilon T$ and
    $|g-p|\leq(L/2)\|(x,y)\|_\infty^s$.

    On every bar with index at least $j+1$, $0\leq x<r$ for large $j$,
    so $|p|\leq\epsilon T$. Taking $\epsilon\leq L/2$ ensures that
    $g$ vanishes there whenever $y>H$.
    The total area of these bars below height $H$ is of order
    $H\omega_{j+1}$. Their contribution to $\int g^2$ is at most
    $C\epsilon^2T^{2+1/s}\omega_{j+1}$.

    On bar $j$, $|p|\leq\epsilon T r^\nu$. Since $T\leq rR^{s-1}$,
    \[
        \frac{T r^\nu}{r^s}
        \leq r^\nu(R/r)^{s-1}\longrightarrow0.
    \]
    Thus $g$ vanishes on this entire bar for large $j$.
    On each outer bar, $x\geq R$ and
    \[
        |p(x,y)|\leq\epsilon T x/r
        \leq\epsilon R^{s-1}x\leq\epsilon x^s.
    \]
    Hence $g$ also vanishes on all these bars.

    It remains to bound the cusp outside the bars. For $x\geq2r$,
    $g(x,y)\ne0$ implies
    \[
        x^{s-1}<\frac{2\epsilon T}{Lr}
        \leq\frac{2\epsilon}{L}R^{s-1}.
    \]
    Choose $\epsilon$ so that $(2\epsilon/L)^{1/(s-1)}<1/2$.
    Since $2r<R/2$ for large $j$, the cusp contribution is supported
    in $x\leq R/2$. On this set,
    $|g|\leq C\epsilon T R/r$, and its area is at most
    $(R/2)e^{-2/R}$. Therefore
    \[
        \int_{\Dexp\setminus\cup_i I_i}g^2
        \leq C\epsilon^2T^2(R/r)^2R e^{-2/R}
        =o\bigl(\epsilon^2T^{2+1/s}\omega_{j+1}\bigr).
    \]
    The last comparison is uniform for $B_{j+1}\leq T\leq B_j$.
    Indeed, all factors other than $e^{-2/R}$ have logarithms
    of order $j^2$, whereas $R^{-1}=e^{(j-1)^2}$.
    We have proved
    \begin{equation}\label{eq:Dlad.high.perturbation.norm}
        \int_{\Dlad}g^2
        \leq C\epsilon^2T^{2+1/s}\omega_{j+1}.
    \end{equation}
    The same support bounds give
    $\|g\|_\infty\leq C\epsilon T R/r\leq C\epsilon R^s\to0$.

    Choose a fixed bounded probability density $\psi$ supported in
    the interior of $I_2$, at positive distance from the origin, and set
    \[
        \widetilde g=g-\left(\int_{\Dlad}g\right)\psi,
        \qquad f_0=\mathrm{Leb}(\Dlad)^{-1},\qquad
        f_1=f_0+\widetilde g.
    \]
    The correction has integral $\int g$, vanishes at zero, and tends
    uniformly to zero. On its support $\|(x,y)\|_\infty^s$ is bounded
    away from zero. Thus, for large $j$,
    $|\widetilde g-p|\leq L\|(x,y)\|_\infty^s$ on the local neighbourhood.
    Also $\|\widetilde g\|_\infty\to0$, so $f_1$ is a positive density
    in $\Sigma(s,L)$ with approximating polynomial $f_0+p$.
    Cauchy--Schwarz and the finiteness of $\mathrm{Leb}(\Dlad)$ give
    \[
        \|\widetilde g\|_2^2
        \leq 2\|g\|_2^2+2\left(\int_{\Dlad}g\right)^2\|\psi\|_2^2
        \leq\bigl\{2+2\mathrm{Leb}(\Dlad)\|\psi\|_2^2\bigr\}\|g\|_2^2.
    \]
    Consequently,
    \[
        nD_{\mathrm{KL}}(f_1\|f_0)
        \leq C\epsilon^2 nT^{(2s+1)/s}\omega_{j+1}.
    \]
    Since $f_1(0)-f_0(0)=\epsilon T$, Le~Cam's lemma proves
    \begin{equation}\label{eq:Dlad.high.lower.amplitude}
        \mathcal R_n(s,L)\gtrsim T
        \quad\text{if}\quad
        B_{j+1}\leq T\leq B_j,
        \qquad nT^{c/s}\omega_{j+1}\leq1.
    \end{equation}
    The constant $\epsilon$ can be fixed small enough for the
    product-law divergence to be bounded by a constant less than one.

    \proofpart{Minimax comparison for every sample size}
    If $M_{k-1}\leq n\leq N_k$, take $j=k-1$ and $T=v_k(n)$.
    Then $B_k\leq T\leq B_{k-1}$ and
    $nT^{c/s}\omega_k=1$.
    If $N_k\leq n\leq M_k$, take $j=k$ and $T=B_k$.
    Then $nT^{c/s}\omega_{k+1}\leq1$.
    Applying~\eqref{eq:Dlad.high.lower.amplitude} in these two cases
    and using~\eqref{eq:Dlad.high.phases} yields
    $\mathcal R_n(s,L)\gtrsim\mathfrak r_n(s)$.
    Together with the upper bound and the definition of minimax risk,
    this proves~\eqref{eq:Dlad.high.minimax}.

    \proofpart{Fluctuation envelopes}
    At both thresholds $N_k$ and $M_k$, the profile equals $B_k$.
    Equations~\eqref{eq:Dlad.high.logB}--\eqref{eq:Dlad.high.logM} give
    \begin{align*}
        \log B_k+(s/a)\log N_k
        &=\frac{2b(s-1)}a k+O(1),\\
        \log B_k+(s/a)\log M_k
        &=\frac{2b(2s-1)}a k+O(1),
    \end{align*}
    while both $\sqrt{\log N_k}$ and $\sqrt{\log M_k}$ equal
    $\sqrt a\,k+O(1)$.
    To control intermediate sample sizes, extend the profile to real
    $n$ and put $J(x)=\log\mathfrak r_{e^x}(s)+(s/a)x$.
    On a variance segment its slope is $s/a-s/c<0$.
    On a constant segment its slope is $s/a>0$.
    The endpoint expressions imply $J(x)=O(\sqrt x)$ throughout.
    Hence, on each segment,
    \[
        \frac{\rd}{\rd x}\frac{J(x)}{\sqrt x}
        =\frac{J'(x)}{\sqrt x}-\frac{J(x)}{2x^{3/2}}
    \]
    has the sign of $J'(x)$ for all sufficiently large $x$.
    The lower and upper limits are therefore attained at $N_k$ and
    $M_k$, respectively, and have the values asserted in the theorem.
    Rounding either threshold up to an integer changes its logarithm
    by $o(1)$. The function $\log\mathfrak r_{e^x}(s)$ has bounded
    slopes, so this rounding does not change the limits.
    Finally,~\eqref{eq:Dlad.high.minimax} changes the logarithm of the
    profile by only $O(1)$, for both the minimax risk and the maximal
    risk of the stated estimator.
\end{proof}

\else
\begin{supplement}
\stitle{Supplement to ``A minimax theory for pointwise density estimation on geometrically complex domains''}
\sdescription{Proofs of all results stated in the paper, together with the auxiliary lemmas on tangent measures, Gram factorisations and the singular examples.}
\end{supplement}
\fi
\section*{Acknowledgment}

\noindent K.\ Bertin is supported by FONDECYT regular grants 1221373 and 1230807 from ANID-Chile, Mathamsud AMSUD250012 and the Centro de Modelamiento Matemático (CMM) BASAL fund FB210005 for centers of excellence from ANID-Chile. N.\ Klutchnikoff is supported by FONDECYT regular grant 1221373 from ANID-Chile and Mathamsud AMSUD250012.

\bibliographystyle{imsart-nameyear}
\bibliography{bib}

\end{document}